\documentclass[a4paper,11pt]{book}
\usepackage[T1]{fontenc}
\usepackage[utf8]{inputenc}
\usepackage{lmodern}
\usepackage[pdftex]{graphicx}
\usepackage[brazilian]{babel}
\usepackage{graphicx}	
\usepackage{amsmath}							
\usepackage{amssymb}
\usepackage{authblk}
\usepackage[T1]{fontenc}
\usepackage[utf8]{inputenc}
\usepackage{lmodern} 
\usepackage{amsthm}
\usepackage{mathrsfs}
\usepackage{tikz}
\usetikzlibrary{decorations.markings}
\usetikzlibrary{patterns}
\usepackage{makeidx}
\usepackage{listings}
\usepackage{color}
\usepackage{mathrsfs}
\usepackage{textcomp}

\makeindex

\definecolor{codegreen}{rgb}{0,0.6,0}
\definecolor{codegray}{rgb}{0.5,0.5,0.5}
\definecolor{codepurple}{rgb}{0.58,0,0.82}
\definecolor{backcolour}{rgb}{0.95,0.95,0.92}

\lstdefinestyle{mystyle}{
    backgroundcolor=\color{backcolour},   
    commentstyle=\color{red!85!green},
    keywordstyle=\color{green!35!blue},
    numberstyle=\tiny\color{codegray},
    stringstyle=\color{green!85!blue},
    basicstyle=\footnotesize,
    breakatwhitespace=false,         
    breaklines=true,                 
    captionpos=b,                    
    keepspaces=true,                 
    numbers=left,                    
    numbersep=5pt,                  
    showspaces=false,                
    showstringspaces=false,
    showtabs=false,                  
    tabsize=2,
	emph={Compara_tipos,Busca_Coloracao_Legitima, MatrizIncidencia, Posicao, MatrizA, Oq, OQ, PSE_MI},  
	emphstyle=\color{brown}
}
 
\theoremstyle{plain}
\newtheorem{teo}{Teorema}[chapter] 
\newtheorem{lema}[teo]{Lema} 
\newtheorem{cor}[teo]{Corolário}
\newtheorem{prop}[teo]{Proposição}

\theoremstyle{definition}
\newtheorem{df}[teo]{Definição} 
\newtheorem{ex}[teo]{Exemplo} 
\newtheorem{ob}[teo]{Observação}

\DeclareMathOperator*{\argmin}{arg\,min}

\newcommand{\rnum}[1]{\lowercase\expandafter{\romannumeral #1\relax}}
\newcommand*{\QEDA}{\hfill\ensuremath{\blacksquare}}%
\newcommand*{\QEDB}{\hfill\ensuremath{\square}}%

\begin{document}

\frontmatter 

\thispagestyle{empty}
\begin{center}
    \vspace*{1.5 cm}
    \LARGE{\textbf{Compressão de Entropia e Colorações Legítimas em Planos Projetivos}}\\
    
    \vspace*{1.4cm}
    \Large{Luís Doin}
    
    \vskip 2.5cm
    \large{\textsc{
    Trabalho de conclusão de curso apresentado\\
    ao\\
    Instituto de Matemática e Estatística\\
    da\\
    Universidade de São Paulo\\
    para\\
    obtenção do título\\
    de\\
    Bacharel em Matemática Aplicada}}
    
    \vskip 2.5cm
    Curso: Matemática Aplicada \\ Habilitação em Métodos Matemáticos\\
    Orientador: Prof. Rodrigo Bissacot\\

    
    \vskip 4 cm
    \normalsize{São Paulo, julho de 2017}
\end{center}
\pagebreak
\thispagestyle{empty}
\pagebreak

      


\newpage
\thispagestyle{empty}
    \begin{center}
        \vspace*{1.5 cm}
    \LARGE{\textbf{Compressão de Entropia e Colorações Legítimas em Planos Projetivos}}\\
        \vspace*{2 cm}
    \end{center}

    \vskip 3.5cm

    \begin{flushright}
       Esta versão do trabalho de formatura contém as correções e alterações sugeridas
	pela Comissão Julgadora durante a defesa da versão original do trabalho,
	realizada em 11/07/2017. Uma cópia da versão original está disponível no
	Instituto de Matemática e Estatística da Universidade de São Paulo.

 \vskip 2cm

    \end{flushright}
    \vskip 4.2cm

    \begin{quote}
    \noindent Comissão Julgadora:
    
    \begin{itemize}
		\item Prof. Dr. Rodrigo Bissacot - IME-USP
		\item Prof. Dr. Guilherme Oliveira Mota - CMCC-UFABC
		\item Prof. Dr. Fernando Mário de Oliveira Filho	 - IME-USP
  \end{itemize}
      
    \end{quote}
\pagebreak

\pagenumbering{roman}     


\chapter*{Agradecimentos} 

\

Esse trabalho nasceu de uma ideia ousada de como orientar uma iniciação científica em matemática: ao invés de colocar livros nas mãos de graduandos, coloque artigos! O professor Bissacot vai ainda um passo adiante incentivando seus alunos a atacar conjecturas e desenvolver ideias originais. Essa busca por novidades, mesmo que como objetivo secundário, coloca o aluno verdadeiramente em contato com a vida acadêmica, experienciando problemas muitos distintos daqueles encontrados em sala de aula onde o ambiente e as questões propostas são mais controlados. Dessa forma, ao longo do curso o aluno ganha motivação e ânimo e consegue se lembrar mais facilmente porque mesmo decidiu vencer toda a burocracia que é um curso de graduação. Ao menos essa foi minha experiência e agradeço muito ao professor Bissacot pela orientação nesse caminho. Essa é a parte acadêmica, pois outro agradecimento vai por todo apoio que recebi dele de forma geral. Quem o conhece sabe a pessoa humana que é. 

\

Além de meu professor, também recebi a ajuda de alguns amigos. Dentre eles, quero agradecer especialmente à La$\ddot{\text{i}}$s pelo incentivo em geral e pela ajuda no apêndice de planos projetivos, ao Marcelo, por aguçar minha curiosidade quanto a teoria de compressão de entropia, ao Eric, por sua predisposição natural em ajudar e dividir conhecimento e ao João, por me ajudar a passar em análise estatística! 

\

Gostaria também de ressaltar a contribuição dos professores Fernando Oliveira Filho e Guilherme Mota. Esse texto passou por alterações e correções importantes após a avaliação de ambos.   

\

Não posso deixar de agradecer a CNPq e a FAPESP pelo apoio financeiro a essa pesquisa.

\

Finalmente, agradeço à minha família: minha mãe Márcia, minha irmã Larissa e meu irmão Leonardo, por me incentivarem desde sempre.

\chapter*{Resumo}

\noindent Doin, L. \textbf{Colorindo Planos Projetivos Finitos via Compressão de Entropia}. 
2017.
Trabalho de Conclusão de Curso - Instituto de Matemática e Estatística,
Universidade de São Paulo, São Paulo, 2017.
\\

Esse estudo teve como objetivo fornecer ao aluno experiência em pesquisa em matemática. O caminho escolhido para tanto foi pensado de forma a apresentar ao aluno ferramentas e técnicas importantes não abordadas na graduação. O aluno então estudou o artigo \cite{paper} de Noga Alon e Zoltan F$\ddot{\text{u}}$redi sobre colorações legítimas em planos projetivos finitos onde faz-se uso de diferentes técnicas contidas no Método Probabilístico dentre elas o Lema Local de Lovász (LLL) exposto em \cite{livrinho}. De posse de versões melhoradas do LLL (\cite{artigo Bissacot}  e \cite{LLLL}) e do método de entropy compression sistematizado por Louis Esperet e Aline Parreau em \cite{Entropy}, buscou-se uma melhora dos resultados obtidos por Noga Alon e Zoltan F$\ddot{\text{u}}$redi. Para tanto, desenvolveu-se uma variação do método de entropy compression provando-se que tal variação pode ser aplicada em todos os problemas formulados dentro da linguagem \textit{variable version} do Lema Local de Lovász -- isso é feito na seção 3.2.3. Os resultados obtidos aplicando essa variação são expostos no Capítulo 4. Os resultados desse trabalho também estão expostos em \cite{BL}.
\

\

\noindent \textbf{Palavras-chave:} Lema Local de Lovász, Entropy Compression, Planos Projetivos Finitos.

\tableofcontents

\mainmatter

\chapter*{Introdução}

\

Em meados de 1940 surgia uma ideia para se obter provas de existência de objetos matemáticos por meio da construção de espaços de probabilidade adequados. Muitas técnicas foram desenvolvidas com esse espírito e esse conjunto de ideias, que mais tarde ficaria conhecido como Método Probabilístico, se mostrou poderoso e aplicável em diversas áreas da matemática. Nas palavras de Béla Bollobás \cite{Bollobas}: ``\emph{Existence results based on probabilistic ideas can now be found in many branches of mathematics, especially in analysis, the geometry of Banach spaces, number theory, graph theory, combinatorics and computer science. Probabilistic methods have become an important part of the arsenal of a great many mathematicians. Nevertheless, this is only a beginning: in the next decade or two (escrito em 1985) probabilistic methods are likely to become even more prominent.}''. E posteriormente, em 2001, nas de Niranjan Balachandran \cite{Balachandran}: \emph{``The Probabilistic Method has now become one of the most important and indispensable tools for the Combinatorist. There have been several hundred papers written which employ probabilistic ideas and some wonderful monographs. ... Over the past 2 decades, the explosion of research material, along with the wide array of very impressive results demonstrates another important aspect of the Probabilistic Method; some of the techniques involved are subtle, one needs to know how to use those tools, more so than simply understand the theoretical underpinnings.''}. Pensando no dito por Balachandran, nesse trabalho cada capítulo terá como uma de suas preocupações fornecer ao leitor novas maneiras de aplicar as ideias apresentadas.   
\

Paul Erd$\H{o}$s, um dos mais prolíficos matemáticos de todos os tempos, foi talvez o maior responsável por popularizar o Método Probabilístico publicando diversos trabalhos a partir de 1947 usando esse método como em \cite{artigo 1, artigo 2, Erdos, LLLL}. De fato, Noga Alon e Joel Spencer usam uma foto de Erd$\ddot{\text{o}}$s como capa do seu livro \textit{The Probabilistic Method} \cite{AMP} que é tido como referência por compilar diferentes abordagens e aplicações do método. Justiça seja feita, um dos primeiros trabalhos de que se tem notícia onde tal método é utilizado foi escrito em 1943 por Szele em \cite{artigo 5} enquanto estudava caminhos Hamiltonianos. A abordagem do Método Probabilístico pode ser resumida nas seguintes palavras: para assegurar a existência de um determinado objeto (um grafo com determinadas propriedades pré-estabelecidas, por exemplo), constrói-se um espaço de probabilidade adequado ao problema onde o objeto de interesse tem probabilidade positiva de ocorrer concluindo assim que tal objeto existe mesmo que não saibamos exatamente como construí-lo. Exemplos do método sendo aplicado para provar a existência de procedimentos de rotinas eficientes podem ser encontrados em \cite{MP3} e \cite{MP4} e problemas de coloração em grafos em \cite{MP1} e \cite{MP2}. 
\
Uma das ferramentas mais afiadas do Método Probabilístico é o Lema Local de Lovász (LLL) capaz de provar a existência de objetos onde outras abordagens falham. Estudaremos esse lema no próximo capítulo. O seu ponto fraco, assim como o de todas as técnicas do Método Probabilístico, é não construir o objeto para o qual prova a existência. Dada as pequenas probabilidades com as quais o LLL pode trabalhar, desenvolver um algoritmo que encontra o objeto apontado por Lovász se mostrou uma tarefa desafiadora. 
\

Paralelamente a tudo isso, em 1948 Claude E. Shannon fundava a área de teoria da informação com seu artigo \textit{A Mathematical Theory of Communication} \cite{Shannon}. No ano seguinte o artigo se transformava em livro ganhando o nome \textit{The Mathematical Theory of Communication}, uma pequena porém relevante mudança dada a percepção da generalidade de seu trabalho. 
\

Em 2010 Robin Moser e Gárbor Tardos em seu celebrado artigo \cite{MT} desenvolveram, para um versão do LLL chamada de \textit{variable version}\footnote{Esse nome é devido a Kolipaka e Szegedy por \cite{MT2}.}-- a priori mais restritiva -- um algoritmo que constrói  o objeto cuja existência conhecemos pelo LLL  provando que tal algoritmo pára em tempo finito usando uma bonita ideia baseada em conceitos de compressão de dados estudada na teoria da informação. Essa ideia ficou conhecida como \textit{entropy compression method}\footnote{Esse nome é em devido a Tao por \cite{Tao}.}.  Em 2013 Jaros Grytczuc, Jakub Kozik e Piotr Micek, em \cite{GKeM}, perceberam que essa abordagem poderia ser aplicada diretamente nos problemas obtendo dessa forma novos resultados em sequências não repetitivas usando um argumento de compressão de entropia. Também em 2013, em \cite{Entropy}, Louis Esperet e Aline Parreau sistematizaram a abordagem proposta em \cite{GKeM} para um determinado conjunto de problemas.  Muitos resultados foram obtidos dessa forma como em \cite{Entropy3, Entropy4, Entropy1, Entropy, Entropy2, Entropy5, Entropy6}.  O método pode ser resumido da seguinte forma: desenvolve-se um algoritmo que constrói o objeto passo a passo aleatoriamente. Se em um determinado passo a construção produz uma característica que o objeto não possui então esse passo é desfeito juntamente com alguns mais. A condição para o algoritmo parar é achar um objeto com as características desejadas. A ideia então é usar um argumento de compressão de entropia para provar que o algoritmo pára. Estudamos o artigo \cite{Entropy} no capítulo 3 juntamente com o aparato teórico que o sustenta. Na seção 3.2.3 provamos que essa abordagem pode ser aplicada em todos os problemas que podem ser atacados utilizando-se variable version do LLL. Como aplicação, trabalhamos no problema de existência de colorações legítimas em planos projetivos finitos proposto por  Noga Alon e Zoltan F\"uredi em \cite{paper} resultando em \cite{BL}. O artigo de Alon e F\"uredi é exposto no capítulo 2 e os resultados que obtivemos são expostos no capítulo 4. 

\

Em tempo, no desenvolvimento desse trabalho teve-se como norte escrever um texto, se não totalmente, o mais autocontido possível. Dessa forma nos apêndices A, B e C o leitor encontrará o conteúdo necessário de planos projetivos finitos, funções geradoras e entropia de Shannon respectivamente.

\chapter{O Lema Local de Lovász}

\

\begin{flushright}

\textit{``This course could be called The Probabilistic Method, \\ The Erd$\H{o}$s Method or, my favorite, Erd$\H{o}$s Magic'' - $ \ $ \\ Joel Spencer }

\end{flushright}

\

A frase acima, introduzindo o curso de Grafos Aleatórios na New York University em 2012, descreve bem o atual status que o Método Probabilístico obteve dentro da teoria combinatória. A palavra \textit{magic} não é usada na frase somente por ter sido Paul Erd$\H{o}$s quem popularizou a abordagem, mas pelo poder que este conjunto de técnicas e ideias acabou exibindo quando aplicado aos mais variados problemas. 
\

Por vezes, a probabilidade de um objeto escolhido aleatoriamente ter as propriedades desejadas é razoavelmente alta e então muitas abordagens dentro do método probabilístico obtém êxito. Ferramentas mais sofisticadas, tais como o Lema Local de Lovász (LLL), nos permite provar a existência de objetos que ocorrem com probabilidades muito pequenas. Com o passar dos anos, versões do LLL mais poderosas e/ou mais adequadas a determinados problemas foram desenvolvidas. Nesse capítulo estudamos diferentes enunciados do Lema Local de Lovász (seção 1.1) e ganhamos familiaridade com algumas técnicas para aplicá-los (seção 1.2).   

\

\section{Diferentes Versões do Lema Local de Lovász}

\

Para provar a existência de um objeto matemático, o LLL o define pelo seu complementar. Ao invés de definir as características de objetos de determinada classe, define uma coleção $(A_x)_{x \in X}$  de \textit{eventos ruins} -- num espaço de probabilidade $(\Omega,\mathcal{A},P)$ -- que englobam as características cuja negação é condição necessária e suficiente para pertencer à classe buscada. Dentro da filosofia do Método Probabilístico, quer-se descobrir condições, as mais gerais possíveis, de modo a garantir que com probabilidade positiva nenhum destes eventos ruins ocorre. Em outras palavras, busca-se o mínimo de restrições possíveis para que 

$$\mathbf{P}(\bigcap_{x \in X}\overline{A_x} ) > 0.$$

Daqui em diante, $X$ denotará um conjunto finito e $G$ um grafo tal que $V (G) = X$.

\begin{ex} Seja $(A_x)_{x \in X} $ uma família de eventos independentes em $(\Omega,\mathcal{A},P)$. Basta exigir que $P(A_x) = p_x < 1,\forall x \in X$, pois nesse caso: 

$$ \mathbf{P}(\bigcap_{x \in X}\overline{A_x}) = \prod_{x \in X} \mathbf{P}(\overline{A_x}) =  \prod_{x \in X} (1 - p_x) > 0.$$\QEDB
\end{ex}

\begin{ex} Seja $(A_x)_{x \in X}$ uma família de eventos em $(\Omega,\mathcal{A},P)$ tal que $\mathbf{P}(A_x) = p_x$ sobre a qual não se sabe quais são independentes. Como $\mathbf{P}(\bigcup_{x \in X} A_x) \leq  \sum_{x \in X} \mathbf{P}(A_x) = \sum_{x \in X} p_x $, impondo que $\sum_{x \in X} p_x < 1$ obtemos o resultado. De fato

$$ \mathbf{P}(\bigcap_{x \in X}\overline{A_x}) = \mathbf{P}(\overline{\bigcup_{x \in X}A_x}) = 1 - \mathbf{P}(\bigcup_{x \in X}A_x) \geq 1 - \sum_{x \in X} p_x > 0 .$$\QEDB

\end{ex}

O último exemplo não levou em consideração a dependência entre os eventos. A ideia então é que, caso seja possível usar esta informação, obteremos maiores valores para as probabilidades $p_x$ de forma a ainda a garantir $\mathbf{P}(\bigcap_{x \in X}\overline{A_x}) > 0$. O primeiro resultado nesta direção foi obtido em \cite{Erdos} por Paul Erd$\H{o}$s e László Lovász em 1973, o qual continha a primeira versão do então chamado Lema Local de Lovász. Antes de enunciar o teorema, precisamos da

\

\begin{df} Dizemos que $G$ é um grafo de dependência para a família de eventos $(A_x)_{x \in X}$ em um espaço de probabilidade quando, para cada $x \in X$, $A_x$ é mutuamente independente da família $\lbrace A_y : y \in X \setminus \Gamma^{\ast} (x) \rbrace$, onde $\Gamma^{\ast} (x)$ se refere à vizinhança de $x$ unida com $x$. 
\end{df}

\

É importante ressaltar que o grafo de dependência não é único -- basta observar que dado um grafo de dependência $G$ para uma família de eventos $(A_x)_{x \in X}$ qualquer outro grafo obtido a partir deste adicionando elos a $G$ será um grafo de dependência para a família. O resultado obtido por Erd$\H{o}$s e Lovász pode então ser enunciado da seguinte forma:

\

\begin{teo}[Erd$\H{o}$s e Lovász \cite{Erdos}] Seja $G$ um grafo de dependência para uma família de eventos $(A_x)_{x \in X}$ com grau máximo $\Delta$ e $\mathbf{P}(A_x) = p_x$. se $4p_x\Delta \leq 1$, para todo $x \in$ X. Então $\mathbf{P}(\bigcap_{x \in X}\overline{A_x}) > 0$.
\end{teo}

\

Uma versão mais geral foi obtida por Joel Spencer em 1977 em \cite{Spencer}:

\

\begin{teo}[Spencer \cite{Spencer}] Seja $G$ um grafo de dependência para a família de eventos $(A_x)_{x \in X}$ e suponha que existam $(r_x)_{x \in X}$ números em [0,1) tais que, para cada $x \in X$,

$$ \mathbf{P}(A_x) = p_x \leq r_x \prod_{y \in \Gamma(x)}(1 - r_y).$$

 Então  $ \ \ \mathbf{P}(\bigcap_{x \in X}\overline{A_x}) \geq \prod_{x \in X}(1 - r_x) > 0$.
\end{teo}

\

O artigo \cite{paper} de Noga Alon e Zoltan F$\ddot{\text{u}}$redi mencionado usa a seguinte versão do LLL, cuja prova pode ser encontrada em \cite{livrinho}:

\

\begin{teo}[Lema Local Lovász - caso simétrico] Seja G um grafo de dependência para uma família de eventos $(A_i)_{i \in [n]}$ com grau máximo $\Delta$ e $\mathbf{P}(A_i) = p_i$. Suponha $p(\Delta + 1)e \leq 1$ e considere $p = sup_{i \in [n]} p_i$. Então $\mathbf{P}(\bigcap_{x \in X}\overline{A_x}) > 0$.
\end{teo}

\

Pela praticidade na verificação das hipóteses, esta versão do teorema é muitas vezes preferida pelos autores. Neste caso não é necessário ter uma informação mais refinada a respeito da estrutura do grafo de dependência, apenas o valor da probabilidade dos eventos (ou o supremo destas) e o grau máximo do grafo de dependência.

\

O próximo resultado a ser apresentado é o Lopsided Lema Local de Lovász que leva em  consideração que os eventos ruins não precisam ser independentes. No que concerne ao LLLL, eventos positivamente correlacionados (em uma certa maneira) podem ser considerados como eventos independentes. Isso possibilita a construção de um ``grafo de dependência'' com um grau menor para a mesma família de eventos. 

\

\begin{df} Dizemos que $G$ é um grafo de dependência assimétrico\footnote{A nomenclatura usual em inglês para este tipo de grafo é \textit{lopsidependency graph}.} para a família de eventos $(A_x)_{x\in X}$ em um espaço de probabilidade quando, para cada $x \in X$, 

$$\textbf{P}(A_x\vert\ \bigcap_{y\in Y} A_y) \leq \textbf{P}(A_x)$$

\noindent para todo $Y \subseteq X \setminus \Gamma^*(x)$. 
\end{df}

\

Este fato foi observado por Spencer e Erd$\H{o}$s em 1991 em \cite{LLLL}, onde perceberam a importância de trabalhar com grafos mais gerais. O novo conceito foi usado para uma nova aplicação do Lema Local de Lovász no estudo dos Latin Transversals. Estudaremos essa aplicação na próxima seção. O resultado obtido por Spencer e Erd$\H{o}$s é o seguinte:

\

\begin{teo}[Spencer e Erd$\H{o}$s \cite{LLLL}]  Suponha que G é um grafo de dependência assimétrico para a família de eventos $(A_x)_{x\in X}$ e suponha que existam $(r_x)_{x\in X}$ números em $[0,1)$ tais que, para cada $x$, 

$$\textbf{P}(A_x) = p_x \leq r_x \prod_{y \in \Gamma(x)}(1-r_y).$$ 

Então $\ \ \textbf{P}(\bigcap_{x\in X}\overline{A_x}) \geq \prod_{x \in X}(1-r_x) > 0$.
\end{teo}

\

A seguir, apresentamos a versão do lema fornecida por Bissacot \textit{et al}. em \cite{artigo Bissacot} em 2011. Recentemente, o LLL foi relacionado a um resultado da mecânica estatística. Essa ligação surpreendente foi apontada por Scott e Sokal em \cite{Sokal} em 2005. Nesse artigo, Scott e Sokal concluíram que o LLL é uma reformulação do critério de Dobrushin (\cite{Dobrushin}) para a convergência da pressão da gás de caroço duro em $G$. Em 2007, Fernández e Procacci, em \cite{FP}, forneceram um novo critério para a convergência da pressão do gás de caroço duro em $G$, e mostraram que esse novo critério é sempre mais eficiente que o critério de Dobrushin. Bissacot \textit{et  al.}, em \cite{artigo Bissacot}, usaram o critério de Fernández-Procacci e os resultados obtidos em \cite{Sokal} para melhorar o Lema Local de Lovász. Essa versão do LLL foi usada, por exemplo, para melhorar o resultado mencionado anteriormente em Latin Transversal (\cite{artigo Bissacot}) e para obter novos resultados sobre colorações dos elos do grafo completo $K_n$ (\cite{Kn}). Cabe mencionar que Pegden, em \cite{MT1}, mostrou que o algoritmo de Moser-Tardos (\cite{MT}) vale para esse novo lema. Para enunciar o lema, precisamos da seguinte nomenclatura:

\

 Sejam $X$ um conjunto finito, $\lbrace \mu_x \rbrace_{x \in X}$ uma família de números não negativos e $G$ um grafo tal que $V (G) = X$. 
Para cada $x \in X$ definimos a partição restrita aos vizinhos de $x$ por:

$$ Z_{\Gamma^{\ast}(x)}(\mu) = \sum_{ \stackrel{T \in I(G)}{T \subset \Gamma^{\ast}_{G} (x)}} \prod_{x \in T} \mu_x,$$

\noindent onde $I(G)$ denota a família de subconjuntos independentes em relação ao grafo $G$ do conjunto $X = V (G)$. Definimos também a seguinte função auxiliar:

$$ Z_{\Gamma(x)}(\mu) = \sum_{ \stackrel{T \in I(G)}{T \subset \Gamma_{G} (x)}} \prod_{x \in T} \mu_x.$$

\begin{teo}[Bissacot \textit{et al}. \cite{artigo Bissacot}] Suponha que G é um grafo de dependência para a família de eventos $(A_x)_{x \in X}$ em um espaço de probabilidade $(\Omega,\mathcal{A},P)$ e que existem $\mu = \lbrace \mu_x\rbrace_{x \in X}$ números reais em $[0,\infty)$ tais que, para cada $x \in X$, $ \mathbf{P}(A_x) = p_x$ satisfaz

$$ p_x \leq \frac{\mu(x)}{Z_{\Gamma^{\ast}(x)}(\mu)}.$$

Então $\ \ \mathbf{P}(\bigcap_{x \in X}\overline{A_x}) \geq \prod_{x \in X} (1 - p_x)^{Z_{\Gamma(x)}(\mu)} > 0.$
\end{teo}

\

Agora, reenunciaremos o Teorema 1.5 para esclarecer a relações entre essas duas versões do LLL.

\

\begin{teo} Seja $H$ um grafo de dependência para a família de eventos $\lbrace A_x\rbrace_{x\in X}$ cada um com \textbf{P}($A_x) = p_x$.  Suponha que existam $\lbrace\mu_x\rbrace_{x\in X}$ números reais em $[0,\infty)$ tais que, para cada $x \in X$, 

$$p_x \leq \frac{\mu_x}{\phi_x(\mu)}, \quad \text{onde} \quad \phi_x(\mu) = \sum_{R \subseteq \Gamma^*_H(x)}\prod_{x \in R}\mu_x.$$

Então  $\ \ \textbf{P}(\bigcap_{x\in X}\overline{A}_x) > 0.$
\end{teo}

\

A única diferença entre os dois teoremas acima é que $\phi_x(\mu)$ soma sobre todos os subconjuntos de $\Gamma^*_H(x)$ enquanto $Z_{\Gamma(x)}(\mu)$ soma apenas sobre os subconjuntos independentes de $\Gamma^*_H(x)$. Disso resulta que 

$$ Z_{\Gamma(x)}(\mu) \leq \phi_x(\mu).$$ 

Como $\phi_x(\mu)$ não leva em consideração a estrutura de $\Gamma^*_H(x)$, o Teorema 1.10 fornece o mesmo resultado se $\Gamma^*_H(x)$ for um conjunto independente ou um clique. Já o Teorema 1.9 fornece sua maior melhora se $\Gamma^*_H(x)$ for um clique e praticamente não fornece melhora se $\Gamma^*_H(x)$ for um conjunto independente. 
\

A desigualdade a seguir é de grande utilidade quando pretende-se aplicar o Teorema 1.9: seja $x$ um vértice do grafo de dependência $H$ para a família de eventos $\lbrace A_x\rbrace_{x\in X}$ e suponha que $\Gamma^*_H(x)$ é a união de $c_1,\ldots,c_k$ cliques então, pela definição de $Z_{\Gamma(x)}(\mu)$,

$$ Z_{\Gamma(x)}(\mu) \leq \prod_{i=1}^k(1 + \sum_{y\in c_i}\mu_y).$$   

\

Assim como anteriormente, o seguinte resultado também é válido. 

\

\begin{teo} Seja $G$ um grafo de dependência assimétrico para a família de eventos $(A_x)_{x\in X}$. Suponha que existam $\mu = \{\mu_x\}_{x\in X}$ números reais em $[0,\infty)$ tais que, para cada $x\in X$ temos que 

$$ \textbf{P}(A_x) = p_x  \leq \frac{\mu_x}{Z_{\Gamma(x)}(\mu)}.$$

\noindent Então $\ \ \textbf{P}(\bigcap_{x\in X} \overline{A}_x)  > 0.$
\end{teo}

\section{Aplicações do Lema Local de Lovász}

\

O objetivo dessa seção é fornecer ao leitor alguma familiaridade com as diferentes formas de se aplicar o LLL. Para tanto, selecionamos problemas que requerem o Lema em suas diferentes versões e organizamos os problemas em uma ordem didática, i.e., as proposições evoluem em nível de dificuldade e quando uma proposição usa uma versão mais antiga do LLL a proposição seguinte mostrará como o mesmo problema pode ser atacado com uma versão mais atual, obtendo assim melhores resultados\footnote{Pode parecer então que não vale a pena investir tempo aprendendo a usar uma versão antiga do LLL. Acontece que existem problemas que não conseguem ser resolvidos com versões mais atuais mas aceitam versões mais antigas.}. Vamos a elas:

\begin{prop} Sejam $G = (V, E)$ um grafo com grau máximo $\Delta$ e 
$V = V_1 \cup V_2 \cup \ldots \cup V_n$ uma partição de $V$. Suponha ainda que para cada conjunto $V_i$ temos $\vert V_i
\vert \geq 2e\Delta$, para todo $1 \leq i \leq n$. Então existe um
conjunto independente $W \subseteq V$ (i.e. não existem $a,b \in W$ tais que $(a,b) \in E$) de cardinalidade $n$ que contém exatamente um vértice de cada $V_i$.
\end{prop}

\textbf{Prova:} Aplicaremos o caso simétrico do LLL: seja $k = min\lbrace\vert V_i\vert : 0 \leq i \leq n \rbrace$. Assumiremos que $\vert V_i\vert = k$ para todo $i$ (o caso geral segue deste considerando a união de $n$ subconjuntos de cardinalidade $k$ onde cada um deles é subconjunto de algum $V_i$). Sejam, para cada $i$, $\psi_i$ variáveis aleatórias uniformes independentes cujos domínios são os respectivos $V_i$. Seja $W$ o conjunto formado pela união dos elementos gerados por uma amostragem de cada $\psi_i$.
\

Para cada elo $ab \in E$, seja $W_{ab}$ o evento ``$W$ contém os vértices $a$ e $b$''. Logo
$P(W_{ab}) = 0$, se $a$ e $b$ são elementos do mesmo $V_i$, e $P(W_{ab}) = 1/k^2 = p$, se $a \in V_i$ e $b \in V_j$
, com
$i  \neq j$. Note que, se os vértices $a$ e $b$ pertencem a $V_i \cup V_j$, temos que o evento $W_{ab}$ é mutuamente independente da família de eventos composta por todos os eventos $W_{cd}$ tais que nem $c$ e nem $d$ pertencem a $V_i \cup V_j$. Logo, como eventos de probabilidade zero são independentes em relação a quaisquer outros eventos, podemos construir um grafo de dependência $H = (\lbrace ab\rbrace_{ab\in E(G)},E(H))$ para a família de eventos ($W_{ab})_{ab\in E(G)}$ com grau máximo $2k\Delta - 1$. Com efeito, sejam $a \in V_i$ e $b \in V_j$, com $i  \neq j$. Como existem $k$ vértices em $V_i$ e cada um deles tem grau máximo $\Delta$ então existem no máximo $k\Delta$ eventos $W_{cd}$ com $c \in V_i$. Da mesma forma existem no máximo $k\Delta$ eventos $W_{cd}$ com $c \in V_j$. Logo, desconsiderando o próprio elo $ab$, a afirmação segue. Definamos $H$ como o grafo de dependência para a família $(W_{ab})_{ab\in E(G)}$ com menor número de elos possível ( i.e. existirá um elo de $H$ entre dois elos $ab$ e $cd$ se, e somente se, $W_{ab}$ e $W_{cd}$ forem dependentes).
\

Seja $\lbrace a, b\rbrace \subset V_i \cup V_j$. Pela definição do grafo $H$ temos que se $cd \in \Gamma_H(ab)$, onde $\Gamma_H(ab)$ é a vizinhança do elo $ab$ no grafo $H$, então $c \in V_i \cup V_j$ ou $d \in V_i \cup V_j$. Logo $\vert \Gamma_H(ab)\vert \leq 2k\Delta - 1$.

\

Pelo Lema Local de Lovász (caso simétrico), se $p(2k\Delta - 1 + 1)e \leq 1$ então

$$P(\cap_{ab\in E(G)}\overline{W}_{ab}) > 0.$$

A proposição segue pois, como $p = 1/k^2$ e $k = \vert V_i\vert$, $p(2k\Delta - 1 + 1)e \leq 1$ se, e somente se $\vert V_i \vert \geq 2e\Delta$. \textcolor{white}{a}\textcolor{white}{a}\QEDA

\

Conseguiremos agora uma melhora do resultado anterior aplicando uma versão mais refinada do LLL:

\

\begin{prop} Sejam $G = (V,E)$ um grafo com grau máximo $\Delta$ e $V = V_1 \cup V_2 \cup \ldots \cup V_n$ uma partição de $V$. Suponha ainda que para cada conjunto $V_i$ temos $\vert V_i\vert \geq 4\Delta$, para todo $ 1 \leq i \leq n$. Então existe um conjunto independente $W \subseteq V$ de cardinalidade $n$ que contém exatamente um vértice de cada $V_i$.
\end{prop}

\textbf{Prova:} Seja $H$ o grafo de dependência construído na proposição 1.2. Para obter essa melhora em relação à proposição 1.12, aplicaremos o teorema 1.9 usando o grafo de dependência $H$. Para tanto, precisamos encontrar uma cota superior para $Z_{\Gamma^*(ab)}(\mathbf{\mu})$, 
onde $ab$ é um elo arbitrário de $G$. A ideia é usar números  \textbf{$\mu$} = $(\mu_{ab})_{ab \in E(G)}$ com $\mu_{ab} = \mu$ para um mesmo $\mu > 0$. Dessa forma será fácil fatorar $Z_{\Gamma^*(ab)}(\mathbf{\mu})$ para posterior análise. Lembrando que, dado um vértice $v$ de um grafo de dependência, o teorema 1.9 usa como informação os conjuntos independentes formados por eventos em $\Gamma(v)$,  precisamos de uma cota superior para o número de pares (pois trincas não são possíveis) de eventos $\lbrace W_{cd},W_{fg} \rbrace$ independentes (em $H$) tais que $cd$ e $fg$ sejam adjacentes a $ab$ no grafo $H$. Claramente o número de pares não excede $k^2\Delta^2$, pois temos no máximo $k\Delta$ eventos $W_{cd}$ tais que $cd$ é adjacente a $ab$ onde ou $c$ ou $d$ não pertence a $V_i \cup V_j$. Como mencionado, não existem trincas $\lbrace W_{cd},W_{ef},W_{gh}\rbrace$ de eventos cujos respectivos elos $\lbrace cd,ef,gh\rbrace$ sejam adjacentes a $ab$ e independentes dois a dois em $H$. De fato, sejam $V_i$ e $V_j$ tais que  $\lbrace a,b\rbrace \subset V_i \cup V_j$. Se $cd \in \Gamma_H(ab)$, então $c \in V_i \cup V_j$ ou $d \in V_i \cup V_j$. O mesmo vale para $ef$ e $gh$. Assim, cada um dos eventos do conjunto $\lbrace W_{cd},W_{ef},W_{gh}\rbrace$ deveria ter seus índices com pelo menos um vértice em $V_i \cup V_j$ e os eventos deveriam ser independentes dois a dois, o que é impossível. Assim:

$$ Z_{\Gamma^*(ab)}(\mathbf{\mu}) \leq 1 + 2k\Delta\mu + k^2\Delta^2\mu^2 = (1 + k\Delta\mu)^2.$$

Seja 

$$ f(\mu) = \frac{\mu}{Z_{\Gamma^*(ab)}(\mathbf{\mu})} \geq \frac{\mu}{(1 + k\Delta\mu)^2}.$$

Como o lado direito da desigualdade acima assume seu valor máximo em $\mu_0 = \frac{1}{k\Delta}$, pelo teorema 1.9, se $p \leq \frac{1}{4k\Delta} \leq f(\mu_0)$ então a probabilidade de nenhum dos eventos da família $(W_{ab})_{ab\in E(G)}$ ocorrer é positiva. Isso conclui a prova, pois $p = 1/k^2$ e então $p \leq 1/4k\Delta$ se, e somente se, $ k \geq 4\Delta$. \textcolor{white}{a}\textcolor{white}{a}\QEDA

\

Sobre este problema, cabe citar o fato de que é conhecida sua constante ótima. Já provamos que o resultado funciona para $\vert V_i \vert \geq c\Delta$ para $c = 2e$ e para $c = 4$ (em particular para qualquer valor maior que este). Já se sabia que $c$ não poderia ser menor que $2$ através de contra-exemplos. Em \cite{Reed} Bruce Reed conjecturou que $2$ era o menor valor de $c$. A conjectura foi provada por Haxell em \cite{Haxell} em 2001  usando uma abordagem alheia ao método probabilístico. 

\

\begin{df} Seja $(a_{ij}) \in M_n(\mathbb{Z}$). Uma permutação $\sigma \colon \lbrace 1,\ldots,n\rbrace \rightarrow \lbrace 1,\ldots,n\rbrace$ é chamada Latin Transversal de $(a_{ij})$ se as entradas $a_{i\sigma(i)}$ são todas distintas para $i \in \lbrace 1,\ldots,n\rbrace$.
\end{df}

\

\begin{prop} Sejam $(a_{ij}) \in M_n(\mathbb{Z})$ e $k \leq (n-1)/4e$. Suponha que nenhum inteiro aparece em mais do que $k$ entradas de $(a_{ij})$. Então $(a_{ij})$ tem uma Latin Transversal.
\end{prop}

\textbf{Prova:} Nesse problema usaremos o teorema 1.8. Precisamos, portanto, como em todas as aplicações do LLL, criar um espaço de probabilidade apropriado e reconhecer quais são nossos eventos ruins (os quais, se forem todos evitados,  teremos a situação desejada). Para tanto, considere o espaço de probabilidade $(\Omega, \mathcal{P}(\Omega), \textbf{P})$, onde $\Omega$ é o conjunto de todas as permutações de $\lbrace 1,\ldots,n\rbrace$ em $\lbrace 1,\ldots,n\rbrace$, $\mathcal{P}(\Omega)$ é o conjunto das partes de $\Omega$ e $\textbf{P}$ é a medida de probabilidade homogênea em $\Omega$. Seja $T$ o conjunto das quatro-uplas ordenadas $(i,j,i^{'},j^{'})$ tais que $i < i^{'}$, $j \neq j^{'}$ e $a_{ij} = a_{i^{'}j^{'}}$. Para cada $(i,j,i^{'},j^{'}) \in T$, seja $A_{iji^{'}j^{'}}$ o evento onde $\sigma(i) = j$ e $\sigma(i^{'}) = j^{'}$. Claramente esses são os nossos eventos ruins. Nosso objetivo é mostrar que  

$$\textbf{P}(\bigcap_{i<i^{'}, j\neq j^{'}}\overline{A}_{iji^{'}j^{'}}) > 0.$$
\

Comecemos pela probabilidade dos eventos ruins:

$$ \textbf{P}(A_{iji^{'}j^{'}}) = \textbf{P}(``\sigma(i) = j" \cap ``\sigma(i^{'}) = j^{'}") = $$
$$ = \textbf{P}(``\sigma(i) = j")\textbf{P}(``\sigma(i^{'}) = j^{'}" \vert ``\sigma(i) = j") = \frac{1}{n} \cdot \frac{1}{n-1} = \frac{1}{n(n-1)}.  $$
\

Agora, seja $G$ o grafo cujo conjunto de vértices é $T$ onde dois vértices $(i,j,i^{'},j^{'})$ e $(p,q,p^{'},q^{'})$ são adjacentes se, e somente se, $\lbrace i,i^{'}\rbrace \cap \lbrace p,p^{'}\rbrace \neq \emptyset$ ou $\lbrace j,j^{'}\rbrace \cap \lbrace q,q^{'}\rbrace \neq \emptyset$. $G$ é um grafo de dependência assimétrico para a família de eventos $A_{iji^{'}j^{'}}$. Com efeito, seja $A_{iji^{'}j^{'}}$ fixado e $B$ o evento ``\textit{existe} $\lbrace x,x^{'}\rbrace$ tal que $ \lbrace x,x^{'}\rbrace \cap \lbrace i,i^{'}\rbrace = \emptyset$ e $\sigma(x) = j$ ou  $\sigma(x^{'}) = j^{'}$''. Seja $Y$ o conjunto dos membros de $T$ que não são adjacentes em $G$ a $(i,j,i^{'},j^{'})$. Logo, se $(p,q,p^{'},q^{'}) \in Y$ então $\overline{A}_{pqp^{'}q^{'}}$ é o conjunto das permutações onde $\sigma(p) \neq q$ ou $\sigma(p^{'}) \neq q^{'}$. Logo, dado que $\overline{A}_{pqp^{'}q^{'}}$ ocorre, aumenta a probabilidade de que $\sigma(p)$ ou $\sigma(p^{'})$ pertença a $\lbrace j,j^{'}\rbrace$. Logo

$$ \textbf{P}( A_{iji^{'}j^{'}} \vert  \overline{A}_{pqp^{'}q^{'}}) = \textbf{P}( A_{iji^{'}j^{'}} \cap B\vert  \overline{A}_{pqp^{'}q^{'}}) + \textbf{P}( A_{iji^{'}j^{'}}  \cap \overline{B}\vert \overline{A}_{pqp^{'}q^{'}}) =$$

$$ = 0 + \textbf{P}( A_{iji^{'}j^{'}}  \cap \overline{B}\vert \overline{A}_{pqp^{'}q^{'}}) = \textbf{P}(\overline{B} \vert  \overline{A}_{pqp^{'}q^{'}})\cdot\textbf{P}( A_{iji^{'}j^{'}} \vert  \overline{B} \cap \overline{A}_{pqp^{'}q^{'}}). $$    

Mas, pelo que foi discutido acima, temos que 

$$ \textbf{P}(\overline{B} \vert  \overline{A}_{pqp^{'}q^{'}}) < \textbf{P}(\overline{B}) \quad \text{e} \quad \textbf{P}( A_{iji^{'}j^{'}} \vert  \overline{B} \cap \overline{A}_{pqp^{'}q^{'}}) = \textbf{P}( A_{iji^{'}j^{'}} \vert  \overline{B}).   $$

Como

$$ \textbf{P}(A_{iji^{'}j^{'}}) = \textbf{P}(A_{iji^{'}j^{'}} \cap B) + \textbf{P}(A_{iji^{'}j^{'}} \cap \overline{B}) = $$

$$ = 0 + \textbf{P}(A_{iji^{'}j^{'}} \cap \overline{B}) = \textbf{P}(\overline{B})\textbf{P}(A_{iji^{'}j^{'}} \vert \overline{B})$$

temos então que

$$ \textbf{P}(A_{iji^{'}j^{'}}  \vert \overline{A}_{pqp^{'}q^{'}}) \leq \textbf{P}(A_{iji^{'}j^{'}}).$$

Pela extensão desse raciocínio temos que

$$ \textbf{P}(A_{iji^{'}j^{'}}  \vert  \ \bigcap_{(p,q,p^{'},q^{'}) \in Y}\overline{A}_{pqp^{'}q^{'}}) \leq \textbf{P}(A_{iji^{'}j^{'}})$$

\noindent provando então que $G$ é um grafo de dependência assimétrico. Ainda, $G$ tem grau máximo $4n(k-1)$. Com efeito, para $(i, j, i^{'}, j^{'})$ fixo, podemos escolher $(s, t)$ de $4n$ maneiras diferentes com $s \in \lbrace i, i^{'}\rbrace$ ou $t \in \lbrace j, j^{'}\rbrace$. Uma vez escolhido $(s, t)$ , temos menos do que $k$ escolhas para $(s^{'}, t^{'})$ diferentes de $(s, t)$, tal que $a_{st} = a_{s^{'}t^{'}}$, pois, por hipótese, existem no máximo $k$ entradas de ($a_{ij})$ com o mesmo valor. Assim, temos menos do que $4nk$ quatro-uplas $(s, t, s^{'}, t^{'})$ tais que $s = i$ ou $s = i^{'}$ ou $t = j$ ou $t = j^{'}$, com $a_{st} = a_{s^{'}t^{'}}$. Agora, para cada uma das quatro-uplas $(s, t, s^{'}, t^{'})$, podemos associar as quatro-uplas $(p, q, p^{'}, q^{'}) = (s, t, s^{'}, t^{'})$ se $s < s^{'}$ ou as quatro-uplas $(p, q, p^{'}, q^{'}) = (s^{'}, t^{'}, s, t)$ se $ s^{'} < s$, o que completa a prova, pois a condição $p(\Delta + 1)e \leq 1$ exigida pelo caso simétrico do teorema 1.8 é equivalente a $k \leq (n-1)/4e$. \textcolor{white}{a}\QEDA

\begin{prop} Sejam $(a_{ij}) \in M_n(\mathbb{Z})$ e $k \leq (n-1)/(256/27)$. Suponha que nenhum inteiro aparece em mas do que $k$ entradas de $(a_{ij})$. Então $(a_{ij})$ tem uma Latin Transversal.
\end{prop}

\textbf{Prova:} Conseguiremos essa melhora em relação à proposição 1.15 aplicando o teorema 1.11. Seja $G$ o grafo definido na proposição 1.15. Como na proposição 1.4, para facilitar a fatoração de $Z_{\Gamma^*(iji^{'}j^{'})}(\mathbf{\mu})$, usaremos \textbf{$\mu$} = $(\mu_{iji^{'}j^{'}})_{iji^{'}j^{'} \in V(G)}$ com $\mu_{iji^{'}j^{'}} = \mu$ para um mesmo $\mu > 0$.
\

O próximo passo é dar uma cota superior para $Z_{\Gamma^*(iji^{'}j^{'})}(\mathbf{\mu})$. Observe que, pela construção de $G$, para um vértice $iji^{'}j^{'}$ fixo, temos que $\Gamma^*(iji^{'}j^{'})$ é a união de 4 subconjuntos, cada um com cardinalidade máxima $nk$ tais que todos os vértices em cada um desses quatro subconjuntos são adjacentes. Logo um conjunto independente de vértices de $\Gamma^*(iji^{'}j^{'})$ tem cardinalidade menor ou igual a 4. Ainda,  

\

$\bullet$ Existem no máximo $4nk$ (a cardinalidade máxima de $\Gamma^*(iji^{'}j^{'})$) conjuntos independentes de vértices de $\Gamma^*(iji^{'}j^{'})$ com cardinalidade 1.

\

$\bullet$ Existem no máximo $\binom{4}{2}(nk)^2$ conjuntos independentes de vértices de $\Gamma^*(iji^{'}j^{'})$ com cardinalidade 2. Com efeito, os vértices deverão ser escolhidos de dois subconjuntos diferentes dos quatro subconjuntos que compõe $\Gamma^*(iji^{'}j^{'})$. Feito isso, existem $nk$ formas de se escolher um vértice em um subconjunto e outras $nk$ formas de se escolher um vértice no outro subconjunto. 

\

$\bullet$ Existem no máximo $\binom{4}{3}(nk)^3$ conjuntos independentes de vértices de $\Gamma^*(iji^{'}j^{'})$ com cardinalidade 3. A prova é muito similar a do caso anterior.   

\

$\bullet$ Existem no máximo $\binom{4}{4}(nk)^4$ conjuntos independentes de vértices de $\Gamma^*(iji^{'}j^{'})$ com cardinalidade 4. A prova é muito similar a do caso anterior.   

\

Logo 

$$ Z_{\Gamma^*(iji^{'}j^{'})}(\mathbf{\mu}) = \sum_{\stackrel{H\subset\Gamma^*(x)}{H \in I(G)} }\prod_{x\in H}\mu_x \leq  (1 + nk\mu)^4$$ 

e

$$ \frac{\mu}{Z_{\Gamma^*(iji^{'}j^{'})}(\mathbf{\mu})} \geq \frac{\mu}{(1+nk\mu)^4} \equiv f(\mu).$$

Como o lado direito assume seu valor máximo em $\mu_0 = 1/3nk$ , podemos usar o Teorema 1.11 na
região $p = 1 / n(n-1) \leq 27/256nk = f(\mu_0)$, o que é equivalente a dizer que $k \leq (n-1)/(256/27)$.\textcolor{white}{a}\QEDA

\

A seguir, definiremos o problema de coloração acíclica de elos e apresentaremos dois resultados obtidos pelo LLL. No capítulo 3, mostramos como Louis Esperet e Aline Parreau, em \cite{Entropy}, melhoraram esses resultados usando outra abordagem.  

\

\begin{df} Uma $\lambda$-coloração dos elos de um grafo $G = (V,E)$ é uma função $f\colon E \rightarrow \lbrace 1,2,\ldots,\lambda\rbrace$. Uma $\lambda$-coloração é dita acíclica se dois elos adjacentes não possuem a mesma cor e se todo ciclo tem pelo menos três cores. 
\end{df}

\

\begin{df} $a$'($G$) é o menor número de cores que podem formar uma coloração acíclica de $G$.
\end{df}

\

\begin{prop} Seja $G = (V,E)$ um grafo com grau máximo $\Delta \geq 3$. Então $a$'($G$) $ \leq  16\Delta $.
\end{prop}

\textbf{Prova:} Nesse proposição aplicaremos o teorema 1.5: sejam $T$ o conjunto dos pares $\lbrace ab,bc\rbrace$, onde $a,b,c \in V$ e $ab,bc \in E$, e, para $k\geq 2$, $C_{2k}(G)$ o conjuntos de todos os ciclos de $G$ de tamanho $2k$. Seja $X = T \cup (\cup_{k\geq 2} C_{2k})$.
Considerando os elementos de $C_{2k}$ como o conjunto de elos que formam cada ciclo temos que os elementos de $X$ são subconjuntos de $E$.  Seja $N = 16\Delta$ o número de cores que usaremos para colorir o grafo. Seguindo a mesma linha das proposições anteriores, definiremos um espaço de probabilidade (dessa vez de forma implícita) e reconheceremos os eventos ruins.
\

Para cada elo de $G$, associe uma cor escolhida uniformemente dentre as $N$ cores disponíveis. Para cada $\lbrace e_i,e_j\rbrace \in T$ seja $A_{\lbrace e_i,e_j\rbrace}$ o evento em que $e_i$ e $e_j$ recebem a mesma cor. Para cada $c_{2k} \in C_{2k}$, seja $A_{c_{2k}}$ o evento em que o ciclo $c_{2k}$ recebe apenas duas cores sem que dois elos adjacentes recebam a mesma cor. Repare que estamos considerando apenas os ciclos de tamanho par. Se um ciclo de tamanho ímpar for colorido com apenas duas cores obrigatoriamente teremos dois elos adjacentes com a mesma cor. Claramente se todos os eventos descritos acima forem evitados então teremos uma coloração acíclica dos elos de $G$. 
\

Como antes, investiguemos a probabilidade de tais eventos ocorrerem: sejam $e_i,e_j \in T$. Se uma cor for atribuída ao elo $e_i$, a probabilidade de a mesma cor ser atribuída ao elo $e_j$ é igual a $1/N$. Logo 

$$ \textbf{P}(A_{\lbrace e_i,e_j\rbrace}) = \frac{1}{N}.$$

Agora, seja $e_1,\ldots,e_{2k}, e_1 \in C_{2k}$. Se uma cor é atribuída a $e_1$ e uma cor diferente é atribuída a $e_2$ há apenas uma maneira de se estender essa coloração de forma a construir um \textit{evento ruim} $A_{c_{2k}}$. Como as escolhas de cores são independentes, tal extensão tem probabilidade $1/N^{2k-2}$ de ocorrer. Logo 

$$\textbf{P}(A_{c_{2k}}) \leq \frac{1}{N^{2k-2}}. $$  

O próximo passo é criar um grafo de dependência para a família de eventos ruins. Como cada evento ruim está associado a um e apenas um elemento de $X$, construiremos um grafo onde os vértices são os elementos de $X$. Como as escolhas de cores são independentes, é claro que o evento $A_{\lbrace e_i,e_j\rbrace}$ é independente do evento $A_{\lbrace e_r,e_s\rbrace}$ se, e somente se, $\lbrace e_i,e_j\rbrace \cap \lbrace e_r,e_s\rbrace = \emptyset$. Assim como $A_{\lbrace e_i,e_j\rbrace}$ é independente de $A_{c_{2k}}$ se, e somente se, $\lbrace e_i,e_j\rbrace\cap c_{2k} = \emptyset$. Da mesma forma, o evento $A_{c_{2k}}$ é independente dos eventos $A_{\lbrace e_i,e_j\rbrace}$ e $A_{c_{2m}}$ se, e somente se, $c_{2k} \cap \lbrace e_i,e_j\rbrace = \emptyset$ e $c_{2k} \cap c_{2m} = \emptyset$, respectivamente. Logo, seja $H(X,F)$ o grafo de dependência para a família de eventos ruins onde $\lbrace x,x^{'}\rbrace \in F$ se, somente se, $x \cap x^{'} \neq \emptyset$. 
\

Observe que cada elo $e$ está contido em no máximo $2(\Delta - 1)$ pares $\lbrace e_i,e_j\rbrace \in T$ e no máximo $(\Delta -1)^{2k-2}$  ciclos $c_{2k} \in C_{2k}$, para todo $k \geq 2$. Logo

\

1) para cada vértice $x = \lbrace e_i,e_j\rbrace \in T$ de $H$, $\Gamma^*_H(x)$ é a união de dois conjuntos $\Gamma^*_1(x)$ e $\Gamma^*_2(x)$ tais que, para $i \in \lbrace 1,2\rbrace$

$$ \vert\Gamma^*_i(x)\vert \leq 2(\Delta - 1) + \sum_{k\geq 2}(\Delta - 1)^{2k-2}$$

e cada elemento $z$ de $\Gamma^*_i(x)$ contém $e_i$ logo o subgrafo de $H$ induzido por $\Gamma^*_i(x)$ é um clique (esse fato só será aproveitado na próxima proposição).

\

2) para $k \geq 2$ e para cada vértice $y = c_{2k} = \lbrace e_1,\ldots, e_{2k},e_1\rbrace \in C_{2k}$ de $H$, $\Gamma^*_H(y)$ é a união de $2k$ conjuntos $\Gamma^*_1(y),\ldots,\Gamma^*_{2k}(y)$ tais que, para $j \in \lbrace 1,\ldots,2k\rbrace$,
 	
$$ \vert\Gamma^*_j(y)\vert \leq 2(\Delta - 1) + \sum_{s\geq 2}(\Delta - 1)^{2s-2}$$

e cada elemento $z$ de $\Gamma^*_j(y)$ contém $e_j$ logo o subgrafo de $H$ induzido por $\Gamma^*_j(y)$ é um clique (de novo, esse fato só será aproveitado na próxima proposição).

\

Dessa vez, não conseguiremos usar o caso simétrico do Lema de Lovász. Sejam, para cada $\lbrace e_i,e_j\rbrace \in T$ e cada $c_{2k} \in C_{2k}$, para $ k \geq 2$, $r_{\lbrace e_i,e_j\rbrace} = 2/N$ e $r_{c_{2k}} = (2/N)^{2k-2}$. Logo a condição do teorema 1.5 está satisfeita. Com efeito, para um $\overline{c}_{2k} \in C_{2k}$, consideremos $A_{\overline{c}_{2k}}$:

$$ r_{\overline{c}_{2k}}\prod_{A_j \in \Gamma(A_{\overline{c}_{2k}})}(1 - r_j) \geq (\frac{2}{N})^{2k-2}\cdot (1 - \frac{2}{N})^{4k\Delta}\cdot \prod_{l\geq 2}\big(1 - (\frac{2}{N})^{2l-2}\big)^{2k\Delta^{2l-2}}.$$

Usando o fato de que, para $z \geq 2$, $(1 - 1/z)^z \geq \frac{1}{4}$, temos que, lembrando que $N = 16\Delta$,

$$ \prod_{l\geq 2}(1 - (\frac{2}{N})^{2l-2})^{2k\Delta^{2l-2}} = \prod_{l\geq 2}((1 - \frac{1}{(N/2)^{2l-2}})^{(N/2)^{2l-2}})^{2k/8^{2l-2}} \geq $$

$$ \geq \prod_{l\geq 2}(\frac{1}{4})^{2k/8^{2l-2}} = (\frac{1}{4})^{2k\sum_{l\geq 2}\frac{1}{64^{l-1}}} = (\frac{1}{4})^{\frac{2k}{63}}.$$

Usando a mesma abordagem podemos mostrar que 

$$ (1 - \frac{2}{N})^{4k\Delta} \geq \big(\frac{1}{4}\big)^{\frac{k}{2}}.$$

Logo

$$r_{\overline{c}_{2k}}\prod_{A_j \in \Gamma(A_{\overline{c}_{2k}})}(1 - r_j) \geq (\frac{2}{N})^{2k-2}\cdot (\frac{1}{4})^{\frac{k}{2}}\cdot (\frac{1}{4})^{\frac{2k}{63}}.$$

Como

$$\frac{(\frac{2}{N})^{2k-2}\cdot (\frac{1}{4})^{\frac{k}{2}}\cdot (\frac{1}{4})^{\frac{2k}{63}}}{{\frac{1}{N^{2k-2}}}} =  {2}^{2k-2}\cdot (\frac{1}{4})^{\frac{k}{2}}\cdot (\frac{1}{4})^{\frac{2k}{63}} > 1 $$

temos o resultado desejado. O mesmo raciocínio pode ser usado para os eventos da forma $A_{\lbrace e_i,e_j\rbrace}$.\textcolor{white}{a}\QEDA

\
 
\begin{prop} Seja $G = (V,E)$ um grafo com grau máximo $\Delta \geq 3$. Então $a$'($G$) $ \leq  \lceil 9.62(\Delta - 1)\rceil$. 
\end{prop}

\textbf{Prova:} Usaremos o mesmo grafo $H$ da proposição anterior dessa vez com $N = c(\Delta - 1)$, para algum $c > 0$ a determinar. Como mencionado, a ideia é usar o teorema 1.9 para melhorar o resultado anterior. Dessa vez, como temos eventos ruins com probabilidades muito diferentes, não poderemos colocar $\mathbf{\mu_x} = \mu$ para todo $x \in X$. Faremos algo similar: para cada $x \in T$, seja $\mu_x = \mu_1 > 0$. Para cada $y \in C_{2k}$, seja $\mu_y = \mu_{k} > 0$. Então, pela desigualdade citada após o teorema 1.10, temos que

$$Z_{\Gamma^*(x)}(\mathbf{\mu}) \leq ( 1 + 2(\Delta - 1)\mu_1 + \sum_{s\geq 2}(\Delta-1)^{2s-2}\mu_s)^2$$
 
e
 
$$Z_{\Gamma^*(y)}(\mathbf{\mu}) \leq ( 1 + 2(\Delta - 1)\mu_1 + \sum_{s\geq 2}(\Delta-1)^{2s-2}\mu_s)^{2k}.$$

Logo devemos achar $\mathbf{\mu}$ de forma que

$$ \frac{1}{N} \leq \frac{\mu_1}{( 1 + 2(\Delta - 1)\mu_1 + \sum_{s\geq 2}(\Delta-1)^{2s-2}\mu_s)^2} $$

e 

$$ \frac{1}{N^{2k-2}} \leq \frac{\mu_k}{( 1 + 2(\Delta - 1)\mu_1 + \sum_{s\geq 2}(\Delta-1)^{2s-2}\mu_s)^{2k}}.$$

Seja $\mu_1 = q = \frac{\alpha}{(\Delta - 1)}$, com $0 < \alpha < 1$, e $\mu_k = q^{2k-2}$. Lembrando que $N \geq c(\Delta - 1)$, as inequações acima tomam a forma

$$\frac{1}{c} \leq \frac{\alpha}{(1+2\alpha+\sum_{s\geq 2}\alpha^{2s-2})^2}$$

e

$$\frac{1}{c} \leq \frac{\alpha}{(1+2\alpha+\sum_{s\geq 2}\alpha^{2s-2})^{\frac{2k}{2k-2}}}.$$

Como $k\geq 2$, a primeira inequação garante a segunda. Logo a condição que garante que nenhum \textit{evento ruim} ocorre é

$$ \frac{1}{c} \leq \frac{\alpha}{(1+2\alpha+\sum_{s\geq 2}\alpha^{2s-2})^2}$$

ou seja

$$ c \geq \frac{(1 + 2\alpha + \frac{\alpha^2}{1-\alpha^2})^2}{\alpha}.$$

Minimizando a função acima chega-se ao resultado que se $c \geq 9.62$ então todo grafo $G$ com grau máximo $\Delta$ admite uma coloração acíclica de elos usando $N = c(\Delta - 1)$ cores. \textcolor{white}{a}\QEDA

\

$$ \diamond \ \diamond \ \diamond $$

\chapter{O Artigo de Alon e F$\ddot{\text{u}}$redi}

\




\begin{flushright}

\textit{``Na realidade trabalha-se com poucas cores. O que dá a $ \ $ \\ ilusão do seu numero é serem postas em seu justo lugar.'' -  \\ Pablo Picasso}
\end{flushright}

\

Nesse capítulo $(P,L,R)$ denotará um plano projetivo finito\footnote{Visite o apêndice A para definições e resultados básicos de planos projetivos finitos.} de ordem $n$ com conjunto de pontos $P$, de linhas $L$ e relação $R\colon P \rightarrow L$.  

\

Dado $c \in \mathbb{N}$, uma $c$-coloração dos ponto de $P$ é uma função $f \colon P \rightarrow$ $\lbrace 1,\ldots,c \rbrace$. A mesma $c$-coloração também pode ser representada pela $c$-upla $(  f^{-1}({1}) , \ldots,  f^{-1}({c}) )$ onde  $f^{-1}(i)$ é o conjunto de todos os pontos que $f$ colore com a cor $i$.
Em relação a uma $c$-coloração, define-se o tipo de uma linha de $L$ como sendo $t_{l,f} := ( \vert l \cap f^{-1}(1) \vert, \ldots , \vert l \cap f^{-1}(c) \vert  )$ (aqui $l$ representa o conjunto que contém todos os pontos $p$ tais que $(p,l)$ $\in$ $R$. A mesma representação será usada mais adiante). 
Se todas as linhas de um $(P,L,R)$ tiverem tipos diferentes em relação a uma $c$-coloração, dá-se a essa coloração a alcunha de 
\textit{coloração legítima}. Se, para $l_i,l_j$ $\in$ $L$, $t_{li,f}$ = $t_{lj,f}$ então $\lbrace l_i,l_j \rbrace$ é dito ser um par ruim. 

O objeto de estudo aqui é o número mínimo de cores necessárias -- denotado por $c(P_n)$ -- para que exista uma $c$-coloração legítima de $(P,L,R)$.
O artigo \cite{paper} de Noga Alon e Zoltan F$\ddot{\text{u}}$redi afirma que se a ordem de $(P,L,R)$ for suficientemente grande (maior que $ 10^{250}$), então $5 \leqslant c(P_n) \leqslant 8$. 
Em outras palavras: para um $(P,L,R)$ de ordem suficientemente grande, se $c < 5$ então em toda $c$-coloração de $(P,L,R)$ existem pelo menos duas linhas com o mesmo tipo e
se $c = 8$ então existe uma $c$-coloração legítima de $(P,L,R)$, qualquer que seja sua ordem. 
Que $c(P_n)$ tem uma cota inferior já era esperado. Porém, surpreende $c(P_n)$ ter uma cota superior, i.e., com 8 cores é possível colorir de forma legítima planos projetivos de ordem tão grande quanto se queira. Isso é contra-intuitivo pois quanto maior a ordem do plano mais linhas temos para atribuir tipos diferente. Porém, também ganhamos mais pontos em cada linha o que faz com que consigamos produzir mais tipos com o mesmo número de cores. A dúvida seria se ganhamos mais tipos do que linhas. A resposta foi dada por Alon F$\ddot{\text{u}}$redi e é explicada a seguir. 
\

\begin{samepage}

\section{8 Cores São Suficientes}

\

A seguir provaremos o seguinte resultado mostrado em \cite{paper}:

\begin{teo} Para um $(P,L,R)$ de qualquer ordem, $c(P_n) \leqslant 8$.
\end{teo}

\textbf{Prova:} Durante toda a prova, $n$ representará a ordem do $(P,L,R)$ a ser colorido. $n$ será considerado tão grande quanto necessário em todas as passagens.\\ 

\end{samepage}
	
	Para essa prova usaremos os conhecimentos obtidos no capítulo 1 e aplicaremos o caso simétrico do Lema Local de Lovász: seja $C$ uma 8-coloração aleatória uniforme de $P$ (i.e., para cada $p$ $\in$ $P$, $\textbf{P}( C(p) = i ) = \frac{1}{8}$ , para todo i $\in$ $\lbrace 1,\ldots,8\rbrace$). 
	Sejam $A_{li,lj}$ os eventos ``$l_i$ e $l_j$ tem o mesmo tipo em relação a $C$'', para todo $l_i, l_j$ $\in$ $L$ com $i < j$. 
	Claramente esses $\binom{n^2 + n + 1}{2}$ eventos não são mutuamente independentes. Se o fossem, o teorema já estaria provado (ver \cite{livrinho}, pág. 12). Teremos mais trabalho que isso para usar o Lema de Lovász . Buscando maior fluidez e clareza, apenas enunciarei dois lemas - que sustentarão a prova - para só depois prová-los:

\

\begin{lema} Existe $S \subset P$ tal que, para todo $l$ $\in$ $L$,  $\ln n$ $\leqslant$ $\vert l \cap S \vert$ $\leqslant$ $20\ln n$.
\end{lema}	

\

	De agora em diante $S$ denotará tal conjunto. Seja $f \colon P\setminus S$ $\rightarrow$ $\lbrace 1,\ldots,8 \rbrace$ uma 8-coloração parcial de $(P,L,R)$. Para todo $l_i,l_j$ $\in$ $L$, define-se $\lbrace l_i,l_j\rbrace $ como sendo um par perigoso se 
	$d_1(t_{li,f} , t_{lj,f})$ $\leqslant$ $40\ln n$ (sejam $x = (x_1,\ldots,x_n)$ e $y = (y_1,\ldots,y_n)$ tem-se $d_1(x,y) = \sum_{i=0}^n \vert x_i - y_i\vert $). Entendamos essa definição: como, pelo lema 2.2, para todo $l$ $\in$ $L$, $\ln n \leqslant \vert l \cap S|$ $\leqslant$ $20\ln n$ e $\vert l\vert =  n + 1$ então $n + 1 - 20\ln n \leqslant \vert (P\setminus S) \cap l\vert = \vert  l \setminus S\cap l\vert \leqslant n + 1 - \ln n$.
	Logo $f \colon P\setminus S \rightarrow \lbrace 1,\ldots,8 \rbrace$ colore, no mínimo, $n + 1 - 20\ln n$ pontos de $l$. Portanto deixa de colorir, no máximo, $20\ln n$ de $l$ (que são os pontos de $l$ que estão em $S$). 	
	Logo, para cada par $l_i,l_j$ $\in$ $L$, $f$ deixa de colorir, no máximo, $40\ln n$ pontos de $l_i \cup l_j$.
	Suponhamos que, para algum $c \in \lbrace 1,\ldots,8 \rbrace$ e $l_i,l_j$ $\in$ $L$, $\ \vert \ ( \vert f^{-1}(c) \cap l_i \vert - \vert f^{-1}(c) \cap l_j \vert  ) \ \vert= k$. Isso quer dizer que $f$ coloriu, com a cor $c$, $k$ pontos a mais 
	em uma linha do que na outra. Quando essa coloração parcial for estendida a uma coloração total $g \colon P \rightarrow \lbrace 1,\ldots,8 \rbrace$ (colorindo-se os pontos restantes em $S$)
	pode acontecer dessa diferença ser equilibrada ( basta que a linha que recebeu menos cores $c$ tenha, no processo de extensão da coloração parcial para a total,
	mais $k$ pontos coloridos com a cor $c$ e a outra linha não tenha mais nenhum, por exemplo). Desse modo  $ \vert g^{-1}(c) \cap l_i \vert - \vert g^{-1}(c) \cap l_j \vert $ seria igual a zero.
	Se, durante a extensão, o mesmo acontecer para todos os $ \vert f^{-1}(c) \cap l_i\vert - \vert f^{-1}(c) \cap l_j \vert \neq 0 $ teremos $t_{li,g}$ = $t_{lj,g}$.
	Como não queremos que isso aconteça, chamamos de pares perigosos (em relação a uma coloração parcial $f$) os pares $\lbrace l_i,l_j \rbrace$ que podem se tornar ruins em alguma 
	extensão de $f$. Como uma extensão de $f$ só colorirá os pontos que $f$ ainda não coloriu - já vimos que são no máximo $40\ln n$ em cada par de linhas de $L$ - 
	essa extensão poderá diminuir $d_1(t_{li,f} , t_{lj,f})$ em, no máximo, $40\ln n$. Logo, se $d_1(t_{li,f} , t_{lj,f}) > 40\ln n$, $\lbrace l_i, l_j \rbrace$ não se tornará um par ruim
	em nenhuma extensão possível de $f$.

\

\begin{lema} Existe uma 8-coloração $f \colon P\setminus S \rightarrow \lbrace 1,\ldots,8 \rbrace$ na qual nenhum ponto $p$ $\in$ $P$ pertence a mais do que 4 pares perigosos. 
\end{lema}

\

Admitindo esses dois lemas, estamos aptos a aplicar o Lema Local de Lovász (caso simétrico):
	
\

	Seja $f \colon P\setminus S \rightarrow \lbrace 1,\ldots,8 \rbrace$ uma coloração parcial de $(P,L,R)$ que obedeça o Lema 2.3. Seja $C$ uma extensão aleatória uniforme de $f$.
	Para todo par perigoso $\lbrace l_i, l_j \rbrace$ (em respeito a $f$) seja $A_{l_i,l_j}$ o evento ``$\lbrace l_i, l_j \rbrace$ torna-se ruim com respeito à extensão aleatória $C$ de $f$''. 
	Sejam $ \lbrace l_i, l_j \rbrace$ um par perigoso, $S_i$ = $S \cap l_i$ e $S_j = S \cap l_j$ ($S_i$ e $S_j$ são os conjuntos que contém os pontos de $l_i$ e $l_j$, respectivamente, 
	que não foram coloridos por $f$).
	Pelo  lema 2.2 temos que:

$$\ln n \leqslant \vert S_i \vert \leqslant 20\ln n \quad \text{e}\quad \ln n \leqslant \vert S_j \vert \leqslant 20\ln n.$$		

Na prova do Lema 2.4, veremos que isso implica que:		 

$$ \textbf{P}(A_{l_i,l_j}) \leqslant \frac{100}{(\ln n)^{\frac{7}{2}}}.$$							
								
	Não é difícil constatar que cada evento $A_{li,lj}$ é mutuamente independente a todos os eventos $A_{lm,ln}$ que satisfaçam ($S_i \cup S_j) \cap (S_m \cup S_n) = \emptyset$ 
	(basta observar que qualquer informação sobre a coloração dos pontos de $S_m \cup S_n$ não fornece nenhuma informação sobre a coloração dos pontos de $S_i \cup S_j$).
	
	Pelo lema 2.3 temos que não existe $p$ $\in$ $P$ tal que $p$ pertença a mais do que quatro pares perigosos. Disso segue que passam no máximo
	quatro pares perigosos por cada $p$ $\in$ $S_i$ $\cup$ $S_j$. Como $ \vert S_i \cup S_j \vert   \leqslant  40\ln n$, temos que o número de
	pares perigosos $\{l_m, l_n\}$ com $ i \neq  m \neq j \neq n \neq i $ tais que $(S_i \cup S_j) \cap (S_m \cup S_n) \neq 0$ é no máximo $ \vert S_i \cup S_j \vert  4 \leqslant 160\ln n $. Finalmente estamos aptos a aplicar o Lema Local de Lovász (com $\Delta$ = 160 $\ln n$, e $p$ $\leqslant$ $\frac{100}{(\ln n)^{\frac{7}{2}}}$ ) para provar que que 
	$\textbf{P}( \bigcap \overline{A_{li,lj}}) > 0$ pois, como estamos supondo $n$ grande, ($160\ln n + 1 )e\frac{100}{(\ln n)^{\frac{7}{2}}} < 1. $ \textcolor{white}{a}\QEDA
	
	\

\begin{samepage}	
	
\section{Provas dos Lemas}

\

\textbf{Lema 2.2.} Existe $S \subset P$ tal que, para todo $l$ $\in$ $L$,  $\ln n$ $\leqslant$ $\vert l \cap S \vert$ $\leqslant$ $20\ln n$.

\

\end{samepage}

\textbf{Prova:} Para provar esse lema, usaremos o seguinte resultado provado em \cite{McDiarmid} e reenunciado em \cite{Matas}:

\begin{lema} Sejam $X_1,\ldots, X_n$ variáveis aleatórias i.i.d com distribuição Bernoulli($p$) e $f \colon X_1 \times \ldots \times X_n \rightarrow \mathbb{R}$ uma função que satisfaz $\vert f(x) - f(y) \vert = 1$ sempre que $x, y \in \{0,1\}^{n}$ diferirem em apenas uma coordenada e $f(x) \geq f(y)$ sempre $x, y \in \{0,1\}^{n}$ obedecerem $x_i \geq y_i$, para todo $i \in [n]$. Logo, para todo $a > 0$ valem as desigualdades:

$$ \textbf{P} ( \ f \ \geqslant \ \mathbb{E}(f) + a \ ) \leqslant  \exp \{ - \frac{a^2}{2(np + a/3)} \};$$

$$  \textbf{P} ( \ f \ \leqslant \ \mathbb{E}(f) - a \ ) \leqslant  \exp \{ - \frac{a^2}{2np} \}.$$

\end{lema}

Nosso desafio agora é construir um conjunto $S \subset P$ que satisfaz o lema 2.2. Para tanto, considere uma moeda com probabilidade de sair cara $c = \frac{10 \ln n}{n + 1}$. Jogaremos essa moeda para cada $ p \in P$. Se der cara, $p \in S$. Caso contrário, $p \notin S$ . Assim, $\vert S \vert$ é uma soma de Bernoullis com parâmetro $c$. Como $\vert P \vert = n^2 + n + 1$ então $|S| \sim$ Bin($n^2 + n + 1$, $c$). Seja $A_l$ o evento ``$\vert S \cap l \vert $ não satisfaz o lema 2.2''. Como $\vert l \vert = n + 1$ temos que $ \vert S \cap l \vert \sim Bin(n + 1, c) $ com $\mathbb{E}(\vert S \cap l \vert) = 10 \ln n \ $. 
Considere a função $f: \{1,0\}^{n+1} \to \mathbb{R}$ definida por $f(x) = \sum x_i$. Sejam $\widehat{x}, \overline{x} \in \{0,1\}^{n+1}$ vetores que diferem em apenas uma coordenada. Logo

$$ \vert f(\widehat{x}) - f(\overline{x}) \vert = 1. $$

Sejam $\widehat{y}, \overline{y} \in \{0,1\}^{n+1}$ tais que $\widehat{y}_i \geq \overline{y}_i$, para todo $i \in [n+1]$. Logo $f(\widehat{y}) \geq f(\overline{y})$. Portanto, tendo visto que $f$ satisfaz tais propriedades podemos usar as cotas enunciadas no lema 2.4. Como $ f(x_l) = \vert S \cap l\vert$, onde $x_l$ é o vetor de amostragens da linha $l$, temos então que

$$ \textbf{P}( \ \vert S \cap l \vert \geqslant 20\ln n \ ) \leqslant exp\{ - \frac{300 \ln^2 n}{80 \ln n}\} < \frac{1}{n^{3}}$$

e

$$ \textbf{P}( \ \vert S \cap l \vert \leqslant \ln n \ ) \leqslant exp\{ - \frac{81 \ln^2 n}{20 \ln n}\} < \frac{1}{n^{4}}.$$

Logo $$\textbf{P}(A_l) < \frac{1}{n^{4}} + \frac{1}{n^{3}}. $$ 

Sejam $X_i$ e $Y$ variáveis aleatórias com $X_i(\overline{A_{li}}) = 0$ e $X_i(A_{li}) = 1$, para todo $i \in \lbrace 1,\ldots, n^2 + n + 1 \rbrace$ e $Y = \sum_{i=1}^{n^2+n+1} X_i$. Logo temos 

$$ \mathbb{E}(Y) = \mathbb{E}(\sum_{i=1}^{n^2+n+1} X_i) = \sum_{i=1}^{n^2+n+1} \mathbb{E}(X_i) < \sum_{i=1}^{n^2+n+1} \frac{1}{n^{4}} + \frac{1}{n^{3}} = $$

$$ = (n^2 + n + 1)(\frac{1}{n^{4}} + \frac{1}{n^{3}}) < 1. $$

Acima, foi usado a linearidades da esperança e que $n$ é grande. Com isso concluímos que o evento ``para $S$, todas as linhas satisfazem o lema 2.2'' tem probabilidade maior que zero. Vejamos porque:

Seja $(\Omega,\mathscr{P_\Omega}, \textbf{P})$ um espaço de probabilidade com $\Omega = \lbrace A_0, \ldots, A_{n^2 + n + 1}\rbrace$, onde $A_i$ é o evento ``apenas $i$ linhas não satisfazem o lema 2.2'', para todo $i \in \lbrace 0, \ldots, n^2 + n + 1\rbrace$, $\mathscr{P_\Omega}$ é o conjuntos das partes de $\Omega$ e $\textbf{P}(A_i)$ é a probabilidade de $A_i$ sendo que $\textbf{P}(A_{li}) < \frac{1}{n^4}$. 
Por construção, a variável aleatória $Y$ conta quantas linhas não satisfazem o lema 2.2. Sendo assim temos que $\mathbb{E}(Y) = \sum_{i=0}^{n^2+n+1} i \textbf{P}(A_i)$. Faremos um argumento por absurdo:
Supondo que $\textbf{P}(A_0) = 0$ temos que: 

$$1 = \textbf{P}(\Omega) = \textbf{P}(\bigcup_{i=0}^{n^2+n+1} A_i) = \textbf{P}(\bigcup_{i=1}^{n^2+n+1} A_i) = \sum_{i=1}^{n^2+n+1} \textbf{P}(A_i) = 1  $$

Agora, se $\sum_{i=1}^{n^2+n+1} \textbf{P}(A_i) = 1  $ é evidente que $\sum_{i=1}^{n^2+n+1} i\textbf{P}(A_i) > 1  $.
Logo 

$$\mathbb{E}(Y) = \sum_{i=0}^{n^2+n+1} i\textbf{P}(A_i) =  \sum_{i=1}^{n^2+n+1} i\textbf{P}(A_i) > 1 . $$

Absurdo! Pois já estabelecemos acima que $\mathbb{E}(Y) < 1$. Logo $\textbf{P}(A_0) > 0$. Portanto existe um conjunto S que satisfaz o lema 2.2.\textcolor{white}{a}\QEDA

\

\begin{ob} Com a mesma argumentação acima é possível mostrar que existe $S \subset P$ tal que, para todo $l$ $\in$ $L$,  $\ln n$ $\leqslant$ $\vert l \cap S \vert$ $\leqslant$ $11\ln n$. Esse resultado será usado no capítulo 4.

\end{ob}

\textbf{Lema 2.3.} Existe uma 8-coloração $f$: $P\setminus S \rightarrow \lbrace 1,\ldots,8 \rbrace$ na qual nenhum ponto $p$ $\in$ $P$ pertence a mais do que quatro pares perigosos. 

\

\textbf{Prova:} Para provar esse lema, precisaremos estabelecer alguns outros resultados antes. De agora em diante $F$ denotará o conjunto $P\setminus S$, onde $S$ é um conjunto que obedece o lema 2.2.\\ 

\begin{lema} Sejam $ l \in L$ e $T \subset l$ tal que $\vert T\vert = k$ e $\widehat{t} = (t_1,\ldots,t_8)$ um vetor com entradas inteiras não negativas.  Então, para qualquer $g: T \rightarrow \lbrace 1,\ldots,8\rbrace$ e $f: F \rightarrow \lbrace 1,\ldots,8\rbrace$, temos que

$$ \textbf{P}(\ t_{l,f} = \widehat{t} \ \vert \ f(p) = g(p), \ \forall p \in T\ ) <  \frac{\binom {m - k} {\lfloor \frac{m-k}{8}\rfloor,\lfloor \frac{m-k+1}{8}\rfloor, \ldots, \lfloor \frac{m-k+7}{8}\rfloor}}{8^{m-k}},$$

onde $m = \vert l \cap F \vert$. Em particular se $ k \leqslant \sqrt{n}$ a probabilidade condicional acima é menor que $\frac{100}{n^\frac{7}{2}}$.\\
\end{lema}

\textbf{Prova:} Seja $s_i = \vert g^{-1}(i)\vert$. É evidente que $t_{l,f} = \widehat{t} \ $ se, e somente se, o número de pontos coloridos com $i$ de $L\cap F\setminus T$ (i.e., os pontos de $L\cap F$ que sobram para $f$ colorir após $g$ colorir os pontos em $L\cap T$) for igual a $t_i - s_i$. Para esses $\vert L\cap F \setminus T\vert$ pontos, há $8^{m-k}$ possíveis colorações  igualmente prováveis. Dessas, exatamente $\binom {m - k} {t_1 - s_1, t_2 - s_2,\ldots, t_8 - s_8}$ fazem com que $t_{l,f} = \widehat{t}$. Logo a probabilidade condicional acima é calculada por $\frac{\binom {m - k} {t_1 - s_1, t_2 - s_2,\ldots, t_8 - s_8}}{8^{m-k}}$. Esse multinomial assume valor máximo quando os valores de $t_i - s_i$ são os mais próximos possíveis. Com efeito, estabeleceremos, primeiro, o seguinte resultado: sejam $a_1$ e $a_2$ tais que $a_1 + a_2 = c$, com $c$ constante, então $\binom {c} {a_1, a_2}$ assume seu valor máximo quando $a_1$ e $a_2$ estão o mais próximo possível. Com efeito, se $c$ é par considere $a_1 = \frac{c}{2} = a_2$.
Sendo assim, temos que $(a_1 + i)! (a_2 - i)! > a_1! a_2!$, para todo $i \in \lbrace 1, \ldots, \frac{c}{2}\rbrace$, pois 

$$ (a_1 + i)! (a_2 - i)! = a_1!(a_1 + i)\ldots(a_1 + 1)(a_2 - i)! > $$

$$ > a_1!(a_2)\ldots (a_2 - i + 1)(a_2 - i)! = a_1!a_2!.$$ 

Portanto $\binom {c} {a_1, a_2} > \binom {c} {a_1 + i, a_2 - i}$, para todo $i \in \lbrace 1, \ldots, \frac{c}{2}\rbrace$. Se $c$ é ímpar o resultado segue com um raciocínio muito similar. O resultado obtido é facilmente generalizado para $\binom {c} {a_1, a_2, \ldots, a_8}$. E com isso o lema está provado, pois $\lfloor \frac{m-k}{8}\rfloor,\lfloor \frac{m-k+1}{8}\rfloor, \ldots, \lfloor \frac{m-k+7}{8}\rfloor$ são os números mais próximos cuja soma é $m - k$.
Repare que, pelo lema 2.2, $n + 1 - 20\ln n \leqslant m \leqslant n + 1 - \ln n$. Como $$\frac{\binom {m - k} {\lfloor \frac{m-k}{8}\rfloor,\lfloor \frac{m-k+1}{8}\rfloor, \ldots, \lfloor \frac{m-k+7}{8}\rfloor}}{8^{m-k}} $$ aumenta a medida que $m$ diminui (verifique!), temos que, sendo $N = n + 1 - 20\ln n - k$, o maior valor que essa probabilidade condicional pode assumir é 

\

$$\frac{\binom {N} {\lfloor \frac{N}{8}\rfloor,\lfloor \frac{N+1}{8}\rfloor, \ldots, \lfloor \frac{N+7}{8}\rfloor}}{8^{N}}. $$ 

Usando que $a! \sim (\frac{a}{e})^a \sqrt{2\pi a}$ (\cite{Bollobas}, pág. 4) e supondo, sem perda de generalidade,  $n \equiv 0(mod \ 8)$, temos o resultado particular para $k \leqslant \sqrt{n}$:

\

$$\frac{\binom {N} {\lfloor \frac{N}{8}\rfloor,\lfloor \frac{N+1}{8}\rfloor, \ldots, \lfloor \frac{N+7}{8}\rfloor}}{8^{N}} =
\frac{N!}{8^N(\frac{N}{8}!)^8} \sim \frac{(\frac{N}{e})^N \sqrt{2\pi N}}{8^N (\frac{(\frac{N}{8})^{\frac{N}{8}}}{e^{\frac{N}{8}}})^8 (\sqrt{2\pi \frac{N}{8}})^8}$$

Mas 

\

$$\frac{(\frac{N}{e})^N \sqrt{2\pi N}}{8^N (\frac{(\frac{N}{8})^{\frac{N}{8}}}{e^{\frac{N}{8}}})^8 (\sqrt{2\pi \frac{N}{8}})^8} = \frac{N^N 8^4\sqrt{2\pi N}}{8^N(\frac{N}{8})^N (2\pi N)^4 } = \frac{(\frac{4}{\pi})^4 \sqrt{2\pi}}{N^{\frac{7}{2}}} = $$

$$ = \frac{(\frac{4}{\pi})^4 \sqrt{2\pi}}{(n + 1 - 20\ln n - \sqrt{n})^\frac{7}{2}} < \frac{100}{n^{\frac{7}{2}}}.$$

Conforme mencionado anteriormente, esse resultado pode ser usado para provar que $ \textbf{P}(A_{l_i,l_j}) \leqslant \frac{100}{(\ln n)^{7/2}}$. Com efeito, para qualquer coloração de $l_j$, a probabilidade de $l_i$ ter o mesmo tipo é menor que 

$$\frac{\binom {\vert S_i \vert - 1} {\lfloor \frac{\vert S_i \vert - 1}{8}\rfloor,\lfloor \frac{\vert S_i \vert }{8}\rfloor, \ldots, \lfloor \frac{\vert S_i \vert + 6}{8}\rfloor}}{8^{\vert S_i \vert - 1}}.$$ 

Como $\vert S_i \vert > \ln n$, por um raciocínio similar ao desenvolvido acima, temos que essa probabilidade é menor que $\frac{100}{(\ln n)^{7/2}}$.\textcolor{white}{a}\QEDA

\

\begin{cor} Sejam $l_1,l_2 \in L$, $ \ T \subset P \ $ tais que $\vert l_2 \cap T\vert \leqslant \sqrt{n}$. Para qualquer $g:T \rightarrow \lbrace 1, \ldots, 8\rbrace$ tem-se que, em relação a $f$,  

$$\textbf{P}( \ \lbrace l_1,l_2\rbrace \  é \ perigoso \ \vert \ f(p) = g(p), \ \forall p \in T \ ) < \frac{(\ln n) ^ 9}{n^{\frac{7}{2}}}.$$ 
\end{cor}

\textbf{Prova:} Fixada uma coloração de $l_1 \cup T$ temos que, para cada $\widehat{t}$ com $d_1(t_{l_1,f},\widehat{t}) \leqslant 40\ln n$, pelo lema 2.5, a probabilidade de $t_{l_2,f} = \widehat{t}$ é menor que $\frac{100}{n^{7/2}}$. Há menos de $(100\ln n) ^8$  vetores $\widehat{t}$ que distam não mais que $40\ln n$ de $t_{l_1,f}$. Com efeito, seja $ t_{l_1, f} = ((l_1)_1, \ldots, (l_1)_8)$. Considere a seguinte notação: $r$ símbolos $\bullet$ na $i$-ésima coordena de um vetor $\widehat{t}$ significa que $\vert (l_1)_i - \widehat{t}_i\vert = r$. Por exemplo, o vetor   

$$ ( \ \bullet \ \bullet \ \bullet \ \bullet \ \bullet \ ; \ \bullet \ ; \ ; \ ; \ ; \ ; \ \bullet \ \bullet \ ; \ \bullet \ )$$

indica que $\vert (l_1)_1 - t_1 \vert = 5$, $\vert (l_1)_2 - t_2 \vert = 1$, $\vert (l_1)_3 - t_3 \vert = 0$, $\ldots$, $\vert (l_1)_7 - t_7 \vert = 2$ e $\vert (l_1)_8 - t_8 \vert = 1$. Com essa notação e levando em consideração que se $\vert (l_1)_i - \widehat{t}_i\vert = r$ temos que $\widehat{t}_i$ pode ter mais ou menos $r$ pontos coloridos com a cor $i$ é fácil perceber que o número de vetores que distam $D$ de $l_1$ é

$$2^8 \binom{D + 7}{7}.$$ 

\noindent Logo, o número de vetores que distam menos que $40 \ln n$ de $l_1$ é

$$ 2^8\sum_{i=0}^{40 \ln n} \binom{i + 7}{7} = 2^8 \binom{40 \ln n + 8}{8} < (100\ln n)^8.$$

\noindent Acima foi usado que $\binom{n}{n} + \binom{n+1}{n} + \ldots + \binom{n + k}{n} = \binom{n + k + 1}{n + 1}$ e que $n$ é grande. 





Portanto, 

$$\textbf{P}( \ \lbrace l_1,l_2\rbrace \  é \ perigoso \ \vert \ f(p) = g(p), \ \forall p \in T \ ) < (100\ln n)^8 \frac{100}{n^{7/2}} < \frac{(\ln n)^9}{n^{7/2}}.$$\textcolor{white}{a}\QEDA

\

\begin{lema} A probabilidade de existirem $l_1,l_2,l_3 \in L \ $ tais que $\lbrace l_1,l_2 \rbrace$ e $\lbrace l_1,l_3 \rbrace$
sejam perigosos (com respeito a $f$) é menor que $\frac{(\ln n)^{18}}{n}$.
\end{lema}

\textbf{Prova:} Sejam $l_1,l_2,l_3 \in L$. Pelo corolário 2.5 temos que 

$$\textbf{P}( \ \lbrace l_1,l_2 \rbrace \ é \ perigoso \ ) < \frac{(\ln n)^9}{n^{\frac{7}{2}}}.$$ 

\noindent De novo pelo corolário temos que

$$\textbf{P}( \ \lbrace l_1,l_3 \rbrace \ é \ perigoso \ \vert  \ \lbrace l_1,l_2 \rbrace \ é \ perigoso \ ) < \frac{(\ln n)^9}{n^{\frac{7}{2}}}.$$

\noindent Sendo $A_{l_i,l_j}$ o evento ``$\lbrace l_i,l_j \rbrace$ é perigoso'' temos que 

$$\textbf{P}( \ A_{l_1,l_2} \cap A_{l_1,l_3} \ ) = 
\textbf{P}( \ A_{l_1,l_2} \ ) \ \textbf{P}( \ A_{l_1,l_3} \ \vert \ A_{l_1,l_2} \ ) < \frac{(\ln n)^{18}}{n^7}.$$

Como existem $(n^2 + n + 1)\binom {n^2 + n}{2} < n^6$ diferentes formas de escolher o trio $l_1,l_2,l_3 \in L$ então, sendo $A_{l_1,l_2,l_3}$ o evento ``$\lbrace l_1,l_2\rbrace$ e $\lbrace l_1,l_3 \rbrace$ são perigosos'': 

$$ \textbf{P}( \ \bigcup _{ 1 \leqslant i \leqslant n^2 + n + 1, \ i \neq j < k \neq i } A_{l_i,l_j,l_k} \ ) <$$

$$< \sum_{ 1 \leqslant i \leqslant n^2 + n + 1, \ i \neq j < k \neq i } \ \textbf{P}( \ A_{l_i,l_j,l_k}) \  < n^6 \frac{(\ln n)^{18}}{n^7} = \frac{(\ln n)^{18}}{n}.$$\textcolor{white}{a}\QEDA

\

\begin{lema} A Probabilidade de existirem $p \in P$ e $\lbrace l_1,\ldots,l_5, l_1',\ldots,l_5'\rbrace \ $ tais que $ \ p \in l_1 \cap \ldots \cap l_5 \ $ e $ \ \lbrace l_i,l_i'\rbrace \ $ é perigoso para todo $i \in \lbrace 1,\ldots, 5\rbrace$ é menor que $\frac{(\ln n)^{45}}{\sqrt{n}}$.
\end{lema}

\textbf{Prova:} Sejam $p \in P \ $, $ \ l_1,\ldots,l_5,l_1',\ldots, l_5' \in L$, $A_{l_i,l_j}$ o evento ``$\lbrace l_i,l_j\rbrace$ é um par perigoso''. Pelo Corolário 2.5, temos que 

$$\textbf{P}( \ A_{l_1,l_1'} \ ) < \frac{(\ln n)^9}{n^{\frac{7}{2}}} \qquad \text{e} \qquad \textbf{P}( \ A_{l_{i + 1},l_{i + 1}'} \ \vert \ \bigcap_{k = 1}^i A_{l_i,l_i'} \ ) < \frac{(\ln n)^9}{n^{\frac{7}{2}}}. $$ 

Portanto $ \ \textbf{P}( \ \bigcap_{k = 1}^5 A_{l_i,l_i'} \ )< \frac{(\ln n)^{45}}{n^{35/2}}$. \begin{samepage}Como há menos que $ \ (n^2 + n + 1)\binom {n + 1}{5}(n^2 + n + 1)^5 < n^{17}$ formas de escolher $p,l_1,\ldots,l_5,l_1',\ldots,l_5'$ temos que (por um raciocínio análogo da prova do lema 2.7) a probabilidade de existirem tais $p,l_1,\ldots,l_5,l_1',\ldots,l_5'$ é menor que

$$n^{17}\frac{(\ln n)^{45}}{n^{\frac{35}{2}}} = \frac{(\ln n)^{45}}{\sqrt{n}}.$$\textcolor{white}{a}\QEDA

\

\end{samepage}

Com esses resultados, finalmente estamos aptos a provar o lema 2.3. Para melhor compreensão, o reenunciaremos de forma mais precisa:

\

\textbf{Lema 2.3.} A Probabilidade de nenhum ponto $p \in P$ estar em mais de quatro pares perigosos é pelo menos $1 - (\frac{(\ln n)^{18}}{n} + \frac{(\ln n)^{45}}{\sqrt{n}})$ (em particular, existe uma 8-coloração  $f$ de $F = P\setminus S$ na qual nenhum ponto $p  \in P$ pertence a mais do que quatro pares perigosos). 

\

\textbf{Prova:} Esse lema é uma consequência direta dos lemas 2.7 e 2.8. Pelo lema 2.8 a probabilidade de não existirem $l_1,\ldots,l_5 \in L$ e $p \in P$ como na figura 1  tais que existam $l_1',\ldots,\l_5' \in L$ com $\lbrace l_i,l_i'\rbrace$ par perigoso para todo $i \in \lbrace 1,\ldots,5\rbrace$ é pelo menos $1 - \frac{(\ln n)^{45}}{\sqrt{n}}$.

\

\begin{samepage}

\begin{center}
\begin{tikzpicture}[scale=0.4, every node/.style={scale=1.5}]

\draw [ line width=1.5] (-5,5) -- (5,-5) node[ pos=0.35, below=0.45] {$p$} node [pos=0.01, below] {$l_2$};

\draw [ line width=1.5] (-2,5) -- (2,-5) node[ pos=0.05, right] {$l_3$} ;

\draw [ line width=1.5] (2,5) -- (-2,-5)node[ pos=0.05, right] {$l_4$} ;

\draw [ line width=1.5] (5,5) -- (-5,-5) node[ pos=0.01, below=2] {$l_5$} ;

\draw [ line width=1.5] (-5,0) -- (5,0) node[ pos=0.04, above] {$l_1$} ;

\draw [fill = black] (0,0) circle [radius=0.2];

\end{tikzpicture} 
 
\end{center}

\begin{picture}(10,10)

\put(1,1){$\qquad\qquad\qquad\qquad\qquad\qquad\quad$\begin{footnotesize} 
 Figura 1

\end{footnotesize}
}
 
\end{picture}
\end{samepage}

\

 Porém, ainda poderiam existir 

\

\begin{samepage}

\begin{center}

\begin{tikzpicture}[scale=0.38, every node/.style={scale=1.3}]

\draw [ line width=1.5] (-5,5) -- (5,-5) node[ pos=0.35, below=0.45] {$p$} node [pos=0.01, below] {$l_2$};

\draw [ line width=1.5] (0,5) -- (0,-5) node[ pos=0.05, right] {$l_3$} ;

\draw [ line width=1.5] (5,5) -- (-5,-5) node[ pos=0.01, below=2] {$l_4$} ;

\draw [ line width=1.5] (-5,0) -- (5,0) node[ pos=0.04, above] {$l_1$} ;

\draw [fill = black] (0,0) circle [radius=0.2]; 

\end{tikzpicture}
\end{center}

\begin{picture}(10,10)

\put(1,1){$\qquad\qquad\qquad\qquad\qquad\qquad\quad$\begin{footnotesize} 
 Figura 2

\end{footnotesize}
}
 
\end{picture}

\end{samepage}

\

$l_1,\ldots,l_4 \in L$ e $p \in P$ como na figura 2 tais que existissem $l_1',l_1'',l_2',l_3', l_4' \in L$ com $\lbrace l_1,l_1'\rbrace$, $\lbrace l_1,l_1''\rbrace$, $\lbrace l_i,l_i'\rbrace$ pares perigosos para todo $i \in \lbrace 2,3,4\rbrace$. Mas, pelo lema 2.7, temos que a probabilidade de isso não acontecer é pelo menos $1 - \frac{(\ln n)^{18}}{n}$.\begin{samepage} Logo, a probabilidade de não acontecer nenhuma das duas situações acima (o que resulta no lema) é pelo menos $1 - (\frac{(\ln n)^{18}}{n} + \frac{(\ln n)^{45}}{\sqrt{n}})$.  \textcolor{white}{a}\textcolor{white}{a}\QEDA

\

\end{samepage}

\begin{ob} Em \cite{paper}, é conjecturado\footnote{Mais precisamente, é conjecturado que $6 \leq c(P_n) \leq 7$.} que $c(P_n) \leq 7$ . Não é difícil verificar que, com pequenas mudanças, quase todos os resultados dessa seção podem ser provados para $c=7$. A única exceção é o lema 2.8. Sem esse resultado não conseguimos provar aquilo que tínhamos como única finalidade: o lema 2.3, que efetivamente tira a dependência do grau máximo do gráfico de dependência com $n$. Portanto, para provar tal conjectura, se faz necessária outra abordagem. Nesse trabalho, apresentamos melhorias no resultado no sentido de achar valores para $c(P_n)$ para $n$'s pequenos e diminuir a exigência para $n$'s grandes. A conjectura $c(P_n) \leq 7$ segue um problema em aberto.  \textcolor{white}{a} \QEDB   
\end{ob} 

$$ \diamond \ \diamond \ \diamond $$

\chapter{O Método de Entropy Compression}

\

\begin{flushright}
\textit{``If something cannot go on  forever, it must stop.'' - \\ Herbert Stein}
\end{flushright} 

\

Como visto no capítulo 2, o Lema Local de Lovász prova a existência de um determinado objeto porém não o constrói. Para um caso especial do LLL chamado de \textit{variable version}\footnote{Esse nome é devido a Kolipaka e Szegedy por \cite{MT2}.}, tal problema foi resolvido por Moser e Tardos em \cite{MT}  pelo uso do chamado algoritmo de Moser-Tardos que busca no espaço de probabilidade por um ponto que evite todos os eventos ruins. O algoritmo foi posteriormente refinado e estendido para outras situações como em \cite{MT3,MT4,MT2,MT1}. O argumento principal em \cite{MT} é o chamado \textit{entropy compression method}\footnote{Esse nome é devido a Tao por \cite{Tao}.}. Grytczuk, Kozik, e Micek, em \cite{GKeM}, perceberam que aplicando o método de entropy compression diretamente pode-se, na maioria dos casos, conseguir melhores resultados do que aplicando o LLL. Muitos resultados foram obtidos dessa forma como em \cite{Entropy3,Entropy4, Entropy1, Entropy, Entropy2, Entropy5, Entropy6}. A próxima seção estuda como tal método funciona baseando-se no que foi exposto em \cite{Tao} por Terence Tao.

\

\section{O Argumento de Entropy Compression}

\

Uma estratégia comum para garantir que determinado algoritmo pára é explorar alguma propriedade de monotonicidade, isto é, mostrar que alguma quantidade limitada superiormente/inferiormente continua\\ crescendo/diminuindo a cada passo do algoritmo. O método de entropy compression se enquadra nesse tipo de estratégia. Nele, é mostrado que um algoritmo probabilístico pára fazendo-se uso da entropia de Shannon das variáveis do algoritmo. Mais precisamente, suponha que um algoritmo probabilístico faz uso de um vetor de bits aleatório $F$ e a cada passo do algoritmo toma um objeto $A$ (que pode ter sido inicializado aleatoriamente) e alguma porção do vetor $F$ para deterministicamente criar um novo objeto $A'$ e um vetor $F'$ formado pelo vetor $F$ menos as entradas utilizadas no passo em questão. Se, de posse $A'$, $F'$ e mais um pedaço de informação $R'$ formos capazes de resgatar $A$ e $F$ usando menos espaço, então estamos realizando uma compressão sem perdas de $A + F$ em $A' + F' + R'$. Porém, como exposto no apêndice C, a probabilidade de comprimirmos sem perdas $A + F$ em um vetor de tamanho menor que sua entropia de Shannon tende a zero. Logo se tivermos um limite inferior para a entropia de $A + F$ e se mostrarmos que a cada passo\footnote{Precisamos na verdade mostrar que existe um número $M \in \mathbb{N}$ para o qual após $M$ passos diminuímos significativamente em tamanho o vetor $A' + F' + R'$.} diminuímos significativamente em tamanho o vetor $A' + F' + R'$ então o algoritmo deve parar depois de um número finito de iterações. Ainda pode-se levantar a seguinte questão: digamos que $\vert F \vert = t$. É possível que, após $t$ iterações, o vetor $A + F$ ainda seja comprimível. Se aumentarmos o tamanho de $F$, o algoritmo terá mais iterações para continuar a compressão mas, em contrapartida, terá uma mensagem maior para comprimir, podendo assim ainda ficar acima do limite que é a entropia de $A + F$. O que garante que existe um $t$ grande o suficiente para o qual o algoritmo pára antes de $t$ iterações? Para responder essa questão precisamos ser um pouco menos simplistas: quando aumentamos o valor de $t$, aumentamos o tamanho do vetor comprimido (pois aumentamos o tamanho da mensagem original) mas aumentamos também a entropia da mensagem original. A questão é o que aumentamos mais. Se conseguirmos argumentar que a entropia tem um aumento maior que o tamanho do vetor comprimido, aí sim podemos garantir que existe um $t$ para o qual, se o algoritmo chegasse a $t$ iterações, a entropia da mensagem original seria maior que o tamanho da mensagem comprimida. Um absurdo, pois, como $t$ é grande o suficiente, isso contradiz o teorema C.26. Digamos que, para algum $c \in \mathbb{N}$, as entradas de $F$ são coletadas de variáveis aleatórias i.i.d. discretas Uniforme($c$) (como é o caso em aplicações do método de entropy compression). Logo, para cada unidade aumentada em $t$, temos um incremento de $\log_2{c}$ bits de entropia em $F$. Para armazenar essa informação extra, inicialmente são necessários $\lceil\log_2{c}\rceil$ bits. Mas estamos armazenando essa informação de forma comprimida, logo precisaremos de menos do que $\log_2{c}$ bits extras de armazenamento. Portanto o incremento de entropia supera o incremento no tamanho da mensagem comprimida e o algoritmo deve parar antes de $t$ iterações. Vejamos um exemplo exposto em \cite{LF}: 

\

\begin{ex} Considere o conjunto $\lbrace x_1,\ldots,x_n\rbrace$ de variáveis booleanas e sua negações $\lbrace \neg x_1,\ldots, \neg x_n\rbrace$, chamadas de literais. Fixado um natural $k > 1$, definimos uma cláusula (de tamanho $k$) como a disjunção de $k$ literais. Por exemplo 

$$ x_1 \vee x_3 \vee \neg x_2 $$

é uma cláusula de tamanho 3, falsa somente se $x_1$ e $x_3$ são falsos e $x_2$ é verdadeiro. Definimos o suporte de uma cláusula pelo conjunto de variáveis que a formam. Por exemplo a cláusula acima tem suporte $\lbrace x_1, x_2, x_3\rbrace$. Para evitar casos degenerados nós consideraremos apenas as cláusulas que tenham suporte de tamanho $k$. Note que se uma cláusula é falsa então sabemos os valores de todas as suas variáveis. Isso será um fato importante mais adiante. 
\

O problema enunciado a seguir é conhecido como problema de\\ $k$-satisfabilidade booleana: dado um conjunto $S$ de cláusulas de tamanho $k$ envolvendo $n$ variáveis booleanas $x_1,\ldots,x_n$, existe alguma configuração de valores de tais variáveis tal que todas as cláusulas de $S$ são satisfeitas?
\

Esse problema é simples para $k = 2$ e, segundo o teorema de Cook-Levin, NP-Completo para $k \geq 3$. Mas o problema se torna mais simples se impusermos restrições em $S$. Por exemplo, se todas as cláusulas de $S$ têm suportes disjuntos então o problema é resolvido facilmente (basta que uma variável de cada cláusula tenha o valor requerido pela mesma). Supondo agora que cada cláusula de $S$ não intersecta a maioria das outras cláusulas de $S$, nós estamos mais próximos de colocar o problema na linguagem do Lema Local de Lovász. 

\

\begin{teo} Seja $S$ um conjunto de cláusulas de tamanho $k$ tal que o suporte de cada cláusula $ s \in S$ intersecta no máximo $2^{k-C}$ suportes de cláusulas de $S$, onde $C$ é uma constante suficientemente grande. Então as cláusulas de $S$ são simultaneamente satisfazíveis.
\end{teo}

\textbf{Prova:} Atribuindo independente e uniformemente um valor a cada variável, temos que a probabilidade de cada cláusula não ser satisfeita é $2^{-k}$. Para cada $s \in S$, seja $A_s$ o evento ``a cláusula $s$ não é satisfeita''. Como cada cláusula não intersecta mais que $2^{k-C}$ cláusulas temos que o grafo de dependência dessa família de eventos tem grau máximo $2^{k-C}$. Logo, pelo caso simétrico do lema de Lovász, temos que se $C > \log_2 e$ então todos os eventos ruins podem ser evitados.\textcolor{white}{a}\QEDA

\

Como dito anteriormente, apesar de bastante simples, essa prova não nos dá valores para as variáveis para os quais todas as cláusulas estão satisfeitas. Usaremos agora um argumento de compressão de entropia para o mesmo problema:

\

\textbf{Prova:} De novo, começamos a prova atribuindo aleatoriamente valores $a_1,\ldots,a_n$ para as variáveis $x_1,\ldots,x_n$. Seja $ A = (a_1,\ldots,a_n)$. Se $A$ satisfizer todas as cláusulas $S$, não há nada mais a ser feito. Suponha que exista algum subconjunto não vazio $T \subset S$ cujas cláusulas não sejam satisfeitas por $A$.
\

Nossa tarefa agora é modificar $A$ de maneira a reduzir $\vert T \vert$. Se, por exemplo, sempre fosse possível achar uma modificação $A'$ de $A$ cujo conjunto $T'$ fosse menor que $T$ então nós poderíamos repetir esse processo até que todas as cláusulas em $S$ fossem satisfeitas. Uma maneira de tentar atingir esse objetivo é escolher uma cláusula $s$ não satisfeita por A e mudar os valores de seu suporte. Basta que uma variável de seu suporte mude de valor para que $s$ passe a ser satisfeita. Como ficará claro mais adiante, é desejável maximizar a quantidade de entropia do sistema no método de entropy compression. Portanto nós escolheremos aleatoriamente essa modificação de $A$. Nós escolheremos $k$ novos valores para o suporte de $s$ aleatoriamente para criar $A'$ (dessa forma há uma probabilidade de $2^{-k}$ de que essa cláusula continue a não ser satisfeita porém o argumento fica ligeiramente mais simples se não nos preocuparmos com isso).
\

Se todas as cláusulas tiverem suportes disjuntos, essa estratégia funciona sem dificuldade. Mas, quando não podemos garantir isso, temos um problema: toda vez que modificamos $A$ para satisfazer uma clausula $s$ podemos estar fazendo com que outras cláusulas que intersectam o suporte de $s$ deixem de ser satisfeitas podendo assim estar aumentando o tamanho de $T$. A observação chave de Moser é que cada falha de uma cláusula $s$ para um conjunto $A$ nos revela $k$  bits de informação - os valores das variáveis do suporte de $s$ em $A$. O plano é usar cada falha de uma cláusula para comprimir $A$ (e mais alguma informação) sem perdas. O ponto principal é que, a cada passo do processo, a cláusula que estamos consertando seja quase sempre uma que intersecta uma cláusula consertada no passo anterior. Dessa forma, o número total de possibilidades para cada cláusula a ser consertada, dada a cláusula consertada no passo anterior, é $2^{k-C}$, o que requer somente $k-C$ bits de armazenamento\footnote{Cada bite armazena 0 ou 1. Portanto $n$ bits podem representar $2^n$ objetos distintos. Logo, para armazenar $x$ objetos precisamos de $\log_2{x}$ bits.} contra $k$ bits de entropia eliminados. É justamente isso que nos garante que o processo pára em um número finito de passos.

\

Para construir essa prova de maneira mais detalhada, nós precisaremos do seguinte:

\

$\bullet$ Uma $n$-upla $A = \lbrace a_1,\ldots,a_n\rbrace$, que será iniciada aleatoriamente e será modificada ao longo do algoritmo (como mencionado acima);

\

$\bullet$ Um longo vetor de bits aleatório $F$, que será usado, em ordem, sempre que o algoritmo pedir uma entrada aleatória;

\

$\bullet$ Um algoritmo recursivo Fix($s$) que modifique a $n$-upla $A$ para satisfazer uma cláusula $s \in S$ ( e que pode ocasionalmente fazer com que $A$ obedeça outras cláusulas de $S$ anteriormente não satisfeitas) definido recursivamente da seguinte maneira:

\

\textbf{Passo 1.a:} Se $A$ satisfaz $s$ não faça nada.

\

\textbf{Passo 1.b:} Caso contrário, ordene o suporte de $s$ de maneira arbitrária. Leia os primeiros $k$ bits de $F$ ainda não usados (reduzindo $F$ de $k$ bits) e use a $j$-ésima entrada desse passo para redefinir o valor da $j$-ésima variável do suporte de $s$.

\

\textbf{Passo 2:} Encontre todas as cláusulas $s' \in S$ cujo suporte intersecta $s$ e que agora não estejam satisfeitas por $A$ - há no máximo $2^{k-C}$ cláusulas dessa forma. Ordene essas cláusulas arbitrariamente e execute Fix($s'$) para cada cláusula em ordem - assim, o algoritmo original Fix($s$) é colocado ``em espera'' em alguma pilha de uma CPU  enquanto todos os processos Fix($s'$) são executados. Fix($s$) termina assim que todos os Fix($s'$) terminem.

\

Por um argumento de indução é fácil mostrar que se Fix$(s)$ termina então a $n$-upla resultante $A$ satisfará $s$ e qualquer outra cláusula $s'$ já satisfeita por $A$ antes de Fix($s$) ser executado. Logo Fix$(s)$, quando termina, devolve (implicitamente) um conjunto de cláusulas não satisfeitas $T'$ menor do que o conjunto $T$ anterior. Podemos, então, consertar todas as cláusulas de $S$ chamando o algoritmo Fix para cada cláusula ainda não satisfeita no início de cada passo dado que o algoritmo Fix termina para cada uma dessas chamadas.
\

Cada vez que o passo 1.b do algoritmo Fix é chamado, a $n$-upla $A$ é transformada em uma nova $n$-upla $A'$ e o vetor aleatório  $F$ é reduzido a um vetor $F'$. Note que esse processo é reversível. Com efeito, dado que sabemos $\lbrace A', R':= s\}$ então $A$ pode ser reconstruido mudando os valores do suporte de $s$ para os únicos valores que violam $s$ enquanto $F$ pode ser reconstruido adicionando ao início de $F'$ os bits de $A$ no suporte de $s$. Esse tipo de reversibilidade não parece ser muito útil em um argumento de compressão de entropia porque enquanto $F'$ é menor que $F$ por $k$ bits, são necessários $\log_2 \vert S\vert$ bits para armazenar a cláusula $s$ (veja a nota de rodapé da página anterior). Então $F+A \mapsto F'+ A' + R'$ só é uma compressão se $\log_2 \vert S\vert < k$, o que não é o que foi assumido nesse exemplo. Porém, enquanto são necessários $\log_2 \vert S\vert$ bits para armazenar uma cláusula $s$, em uma aplicação recursiva do algoritmo Fix são necessário $k-C+O(1)$ bits para armazenar a cláusula $s'$ que está sendo consertada. Se $C$ for grande o suficiente temos que $k-C+O(1) < k$ o que torna possível nosso argumento de compressão de entropia. Com efeito, as cláusulas $s$ para as quais o algoritmo Fix($s$) é chamado podem ser dividas em duas categorias: 1) aquelas nas quais $s$ veio da lista $T$ original de cláusulas não satisfeitas. Cada uma delas requer $O(\log_2 \vert S\vert)$ bits de armazenamento. Como há no máximo $\vert T\vert$ chamadas desse tipo e $\vert T\vert \leq \vert S\vert$ temos que o total de bits necessários para armazena tais cláusulas é no máximo $O(\vert S \vert \log_2 \vert S\vert)$. 
Ou 2) aquelas para as quais Fix($s$) é chamado recursivamente a partir de uma chamada anterior Fix($s'$). Nesse caso $s$ é uma entre as, no máximo, $2^{k-C}$ cláusulas de $S$ cujo suporte intersecta o suporte de $s'$. Logo podemos nos referir a $s$ usando $s'$ e um número entre $1$ e $2^{k-C}$ representando a posição de $s$ (com respeito a uma ordenação arbitrária de $S$) na lista de todas as cláusulas de $S$ cujo suporte intersecta o suporte de $s'$. Chamaremos esse número de índice da chamada Fix$(s)$.
\

Como mencionado acima, a ideia é mantermos um registro $R$ com o qual possamos reconstruir os vetores originais. Seja $R$ o vetor que armazena 1) $s$ cada vez que uma das $\vert T\vert$ Fix($s$) originais for chamada 2) o índice das chamadas Fix($s$) subsequentes pelo processo recursivo do algoritmo Fix e 3) um símbolo de parada sempre que um algoritmo Fix terminar. Dessa forma, podemos deduzir qual cláusula $s$ está sendo consertada para cada chamada de Fix apenas analisando o registro $R$ até aquele ponto. Logo  em qualquer momento do algoritmo é possível reconstruir $A + F$ pelo estado atual $A' + F' + R'$ do processo.
\

Agora suponha, por contradição, que $S$ não é satisfazível. Logo a pilha de chamadas da função Fix deve ser infinita. Nós caminharemos nessa pilha dando $M$ passos, onde $M$ é algum número grande a ser escolhido. Nesse ponto do algoritmo, o vetor $F$ já foi reduzido de $Mk$ bits e o vetor $R'$ tem no máximo $O( \vert S\vert\log_2 \vert S\vert) + M(k - C + O(1))$ bits.\begin{samepage} Logo, nós temos uma compressão sem perdas $A+F \mapsto A'+ F' + R'$ de $n+\vert F\vert$ bits aleatórios para $n + \vert F \vert - Mk + O(\vert S\vert\log_2 \vert S\vert) + M(k-C+O(1))$ bits.  Mas, como $n + \vert F \vert$ bits aleatórios não podem ser comprimidos em um espaço menor que $n + \vert F \vert$, nós temos o limite\footnote{Por esse motivo, é desejável maximizar a entropia do sistema em um argumento de compressão de entropia.}

$$ n + \vert F \vert - Mk + O(\vert S\vert\log_2 \vert S\vert) + M(k-C+O(1)) \geq n+\vert F\vert \iff $$

$$ \iff O(\vert S\vert\log_2 \vert S\vert) + M(O(1) - C) \geq 0. $$

Logo, se M for grande o suficiente, chegamos à contradição desejada. O teorema 3.2 segue.\textcolor{white}{a}\QEDA \\ \textcolor{white}{a}\QEDB

\end{samepage}

\end{ex}

\section{O Método de Esperet e Parreau}

\

Em \cite{Entropy}, em 2013, Louis Esperet e Aline Parreau descreveram um método usando um argumento de compressão de entropia para resolver problemas de coloração onde algumas configurações são proibidas. Seguindo o proposto no artigo, estudaremos tal método pelo problema de

\

\subsection{Coloração Acíclica de Elos}

\

\begin{df} Uma $\lambda$-coloração dos elos de um grafo $G = (V,E)$ é uma função $f: E \rightarrow \lbrace 1,2,\ldots,\lambda\rbrace$. Uma $\lambda$-coloração é dita acíclica se dois elos adjacentes não possuem a mesma cor e se todo ciclo tem pelo menos três cores. 
\end{df}

\

\begin{df} a'($G$) é o menor número de cores que podem formar uma coloração acíclica de $G$.
\end{df}

\

 Fazendo uso do Lema Local de Lovász, Michael Molloy e Bruce Reed, em \cite{Molloy}, mostraram que se G tem grau máximo $\Delta$ então a'($G$) $\leq 64\Delta$. Usando uma versão mais recente do Lema Local devido a Bissacot \textit{el al.} (\cite{artigo Bissacot}), Ndreca \textit{el al.}, em \cite{Ndreca}, mostraram que a'($G$) $\leq \lceil 9,62(\Delta - 1)\rceil$\footnote{Conforme mostrado na proposição 1.20.}. Do método de entropy compression \cite{Entropy} resulta que a'($G$) $\leq 4(\Delta - 1)$.
 
\

Para tanto, faz-se uso de um simples algoritmo: seja $\gamma > 0$ e $ K = \lceil (2 + \gamma)(\Delta - 1)\rceil$. Ordene os elos de $G$ em $e_1, \ldots, e_m$ (onde $m = \vert E\vert$). Ao começo de cada passo, selecione o elo $e_i$ não colorido de menor índice e sorteie (uniformemente, por exemplo) uma cor em $\lbrace 1,\ldots, K \rbrace$ que não apareça em um elo adjacente a $e_i$. Se algum ciclo for colorido com apenas duas cores, descolora esse ciclo. Com o uso do método de entropy compression será mostrado que, dado que $K$ é grande o suficiente, um algoritmo muito similar a esse pára. 

\

\subsubsection{O Algoritmo}

\

Para a análise do algoritmo, construiremos aleatoriamente o vetor $F \in \lbrace1, \ldots, \lceil\gamma(\Delta - 1)\rceil\rbrace^t$, para um $t$ grande. No $i$-ésimo passo do algoritmo, a $i$-ésima entrada $F_i$ de $F$ será usada para atribuir uma cor ao elo não colorido de menor índice $e_j$  da seguinte forma: seja $e_j = uv$ e $S = \lbrace1,\ldots,K\rbrace \setminus S'$, onde $S'$ é o conjunto das cores atribuídas aos elos $xy \neq uv$ tais que:

\

(1) $x=u$ ou $x=v$

\

ou

\

(2) elos $ux$ e $vy$ existem e tem a mesma cor. 

\

Agora atribuímos o $F_i$-ésimo elemento de $S$ ao elo $e_j$. A condição (1) impede que, ao colorir $e_j$, formem-se dois elos adjacentes com a mesma cor. A condição (2) impede que se forme um ciclo com quatro elos com apenas duas cores, como mostra a figura a seguir. 

\

\begin{center}

\begin{tikzpicture}

\node (1) at ( 1, 0) {x};
\node (2) at ( 2, 2) {u};
\node (3) at ( 6, 2) {v};
\node (4) at ( 5, 0) {y};

\draw [color=brown, line width=1.5] (2) -- (1);

\draw [color=brown, line width=1.5] (3) -- (4); 
\draw [color=blue, line width=1.3] (1) -- (4);

\end{tikzpicture}

\end{center}

\begin{picture}(10,10)

\put(1,1){$\qquad\qquad\qquad$\begin{footnotesize} 
 Figura 1: Representação da condição (2). 

\end{footnotesize}
}
 
\end{picture}

\

Repare que essa instrução está bem definida pois $\vert S\vert \geq\lceil\gamma(\Delta -1)\rceil$. Com efeito, pela figura, é fácil notar que para cada cor contada em (2) uma cor é contada pelo menos duas vezes em (1). Logo $S'$ não contém mais cores do que o número de elos adjacentes a $e_j$. Logo $\vert S'\vert \leq 2(\Delta - 1)$. A condição (1) também impede que se formem ciclos de tamanho ímpar coloridos com apenas duas cores. Se, ao atribuirmos uma cor a $e_j$, forme-se um ciclo com duas cores (de tamanho 6 ou mais), digamos $e_{i_1},\ldots, e_{i_{2k}}, e_{i_1}$ com $e_{i_1} = e_j$  e $i_2 < i_{2k}$ então descolorimos todos os elos desse ciclo exceto $e_{i_2}$ e $e_{i_3}$.
\

A próxima seção se dedica à análise do algoritmo. Para tanto, definimos o vetor $R$ com $t$ entradas da seguinte forma: Suponha que, no $i$-ésimo passo do algoritmo, o elo $e_j$ seja colorido e um ciclo de tamanho $2k$ receba apenas duas cores. Como existe um número finito de possíveis ciclos de tamanho $2k$ que contém o elo $e_j$ (no máximo $(\Delta - 1)^{2k-2}$), podemos ordená-los em $C_1,\ldots, C_s$ com $s \leq (\Delta - 1)^{2k-2}$. Nesse caso, à $i$-ésima entrada $R_i$ de $R$ será atribuído o par $(k,l)$ onde $l \leq s$ é o índice do ciclo de tamanho $2k$ que contém $e_j$. Se, no $i$-ésimo passo do algoritmo, nenhum ciclo for formado então a $R_i$ nada será atribuído. Será de fundamental importância que, com o vetor $R$ e a coloração parcial ao passo $i$, sejamos capazes de reconstruir o vetor $F$ até sua $i$-ésima entrada. 
\

A ideia é mostrar que o conjunto de todos os vetores $F$ para os quais o algoritmo não pára antes do passo $t$ é menor que o conjunto de todas as possibilidades de pares de vetor $R$ com a coloração parcial ao passo $t$. Se mostrarmos que esse último é o($\lceil \gamma(\Delta - 1)\rceil^t$) então provaremos que o algoritmo pára com probabilidade maior que zero, pois $\lceil \gamma(\Delta - 1)\rceil^t$ é justamente o número de escolhas possíveis para o vetor $F$ implicando que o algoritmo pára para algum vetor $F$. 

\

\subsubsection{Análise do Algoritmo}

\

Seja $X_i$ o conjunto dos elos não coloridos após o passo $i$ e $\phi_i$ a coloração parcial de $G$ depois do passo $i$. Seja $F$ um vetor para o qual o algoritmo retorna o par ($R$,$\phi_t$). Os próximos dois lemas mostram que o par ($R$,$\phi_t$) determina o vetor $F$:

\

\begin{lema} A cada passo $i$, o conjunto $X_i$ é unicamente determinado pelo vetor $(R_j)_{j\leq i}$.
\end{lema}

\textbf{Prova:} Provaremos por indução em $i$. É imediato que $X_1$ é unicamente determinado por $(R_j)_{j\leq 1}$ pois $X_1 = E \setminus \lbrace e_1 \rbrace$. \

Agora, supomos que $X_{i-1}$ é unicamente determinado por $(R_j)_{j\leq i-1}$. Ao passo $i$, será colorido o elo $e_j$ onde $j = \min \{ j : e_j \in X_{i-1}\}$. Se $R_i$ é vazio então $X_i = X_{i-1} \setminus {e_j}$. Caso contrário, lendo $R_i$ sabemos exatamente quais vértices foram descoloridos nesse passo. Logo, também conhecemos o conjunto $X_i$.  \textcolor{white}{a}\QEDA

\

\begin{lema}  A cada passo $i$, a função que, a cada vetor $(F_j)_{j\leq i}$, associa o par $((R_j)_{j\leq i},\phi_i)$ é injetiva. 
\end{lema}

\textbf{Prova:} Novamente o argumento será por indução. Após o primeiro passo, $\phi_1$ é a configuração em que apenas o primeiro vértice está colorido e com a cor $F_1$. Supondo, para $i > 1$, que o lema é verdadeiro para $i -1$ provaremos que ele também vale para $i$. Pelo lema 3.5, nós conhecemos $X_i$ e $X_{i-1}$. Em particular nós conhecemos o elo $e_j$ que foi colorido no passo $i$.
\

Se $R_i$ é vazio então $\phi_{i-1}$ é obtido através de $\phi_i$ apenas descolorindo o elo $e_j$. Logo, por indução, nós conhecemos o vetor $(F_j)_{j < i}$ bastando determinar a entrada $F_i$. Seja $c$ a cor atribuída ao elo $e_j = uv$ e seja $\mathit{a}$ o número de diferentes cores  $\lbrace i : i < c\rbrace$ que aparecem em $\phi_{i-1}$ em (1) elos adjacentes a $e_j$ e (2) em elos $xy$ tais que $xu$ e $yv$ são elos de $G$ e possuem a mesma cor. Então $F_i = c - \mathit{a}$.
\

Agora, seja $R_i = (k,l)$. Então nós conhecemos o ciclo $e_{i_1},\ldots,e_{i_{2k}}, e_{i_1}$, com $e_{i_1} = e_j$, que foi colorido com duas cores no passo $i$. Nesse caso, $\phi_{i-1}$ é obtido a partir de $\phi_i$ apenas colorindo os elos $e_{i_5}, e_{i_7},\ldots,e_{i_{2k-1}}$ com a cor $\phi_i(e_{i_3})$ e os elos $e_{i_4}, e_{i_6},\ldots,e_{i_{2k}}$ com a cor $\phi_i(e_{i_2})$. Além disso, sabemos que, no passo $i$, $e_j$ recebeu a cor $\phi_i(e_{i_3})$ logo antes de ser descolorido. Pelo mesmo raciocínio acima, $(F_j)_{j\leq i-1}$ está determina pela hipótese de indução e $F_i$ está determinado pois conhecemos a cor que o elo $e_j$ recebeu no passo~ $i$. \textcolor{white}{a}\QEDA

\

Seja $\mathcal{F}_t$ o conjunto dos vetores $F$ tais que, ao passo $t$ do algoritmo, o grafo $G$ não está completamente colorido. Seja $\mathcal{R}_t$ o conjunto dos vetores $R$ que são gerados a partir de vetores de $\mathcal{F}_t$. Como há no máximo $(K + 1)^m$ possibilidades de colorações parciais $\phi_t$ de $G$, os dois lemas anteriores implicam que

\

\begin{lema} $\vert\mathcal{F}_t \vert \leq (K + 1)^m\vert\mathcal{R}_t \vert$. 
\end{lema}

\

O desafio agora é achar uma cota para $\vert\mathcal{R}_t \vert$ a fim de mostrar que, para $t$ grande o suficiente, $\vert\mathcal{F}_t \vert$ é menor que o conjunto de todos os vetores $F$ possíveis. Implicando que existe um vetor $F$ para o qual o algoritmo pára. Para tanto, transformaremos $\vert\mathcal{R}_t \vert$ em um conjunto mais familiar:

\

Considere uma palavra $w = w_1,\ldots,w_{2k-2}$ no alfabeto $\mathbf{A} = \lbrace 1,\ldots,\Delta - 1\rbrace$, e a função $\theta_k(w) = 1 + \sum_{i=1}^{2k-2}(w_i - 1)(\Delta - 1)^{i-1}$. A função $\theta_k$ tem imagem igual a $\lbrace 1, \ldots, (\Delta - 1)^{2k-2}\rbrace$ e é injetiva. Com efeito, seja $\gamma$ e $\psi$ tais que $\theta_k(\gamma) = \theta_k(\psi)$. Logo 

$$\theta_k(\gamma) - \theta_k(\psi) = \sum_{i=1}^{2k-2}(\gamma_i - 1)(\Delta - 1)^{i-1} - \sum_{j=1}^{2k-2}(\psi_i - 1)(\Delta - 1)^{j-1} =$$ 

$$= \sum_{i=1}^{2k-2}((\gamma_i - 1)(\Delta - 1)^{i-1} - (\psi_i - 1)(\Delta - 1)^{i-1}) = $$

$$ = \sum_{i=1}^{2k-2}((\gamma_i - 1) - (\psi_i - 1))(\Delta - 1)^{i-1} = \sum_{i=1}^{2k-2}(\gamma_i - \psi_i)(\Delta - 1)^{i-1}.$$ 

Mas $\vert (\gamma_i - \psi_i) \vert < (\Delta - 1)$, $\forall i \in \lbrace 1,\ldots, 2k - 2\rbrace$. Assim $\sum_{i=1}^{2k-2}(\gamma_i - \psi_i)(\Delta - 1)^{i-1}$ é uma combinação linear de potências de $(\Delta - 1)$ com coeficientes inteiros e menores, em modulo, que $(\Delta - 1)$. Portanto $\sum_{i=1}^{2k-2}(\gamma_i - \psi_i)(\Delta - 1)^{i-1} = 0 \iff \gamma_i  = \psi_i$, $\forall i \in \lbrace 1,\ldots, 2k - 2\rbrace$ pois esse somatório é uma soma de números na base $(\Delta - 1)$ onde as parcelas são da forma $a00\ldots000$ e existe apenas uma parcela não nula com $i$ zeros à direta para cada $i$. Logo $\theta_k(w)$ é injetiva. É imediato que $\theta_k(w)$ assume seu valor mínimo quando $w_i = 1$, $\forall i \in \lbrace 1,\ldots, 2k-2\rbrace$ - quando $\theta_k(w)$ é  igual a 1 - e assume seu valor máximo quando $w_i = (\Delta - 1)$, $\forall i \in \lbrace 1,\ldots, 2k-2\rbrace$ - quando $\theta_k(w)$ é igual a $(\Delta - 1)^ {2k - 2}$. Como existem exatamente $(\Delta - 1)^ {2k - 2}$ palavras $w$ diferentes e a função $\theta_k(w)$ é injetiva segue que a função tem imagem  $\lbrace 1,\ldots, (\Delta - 1)^{2k-2}\rbrace$.
\

Seja $R \in \mathcal{R}_t$ e $R^* = (R_i^*)_{i\leq t}$ uma sequência de $t$ palavras no alfabeto $\mathbf{A}^* = \mathbf{A} \cup \lbrace 0 \rbrace$ assim definida: para todo $1 \leq i \leq t$, se $R_i$ é vazio então $R_i^* = 0$. Caso contrário, existem $k$ e $l$ tais que $R_i = (k,l)$ e então $R_i^*$ será a concatenação de $0$ com $\theta_k^{-1}(l)$ ($R_i^*$ está bem definido pois, lembremos, $l \leq (\Delta - 1)^{2k-2}$). Agora consideramos a sequência de palavras $R^*$ como uma palavra $R^\bullet$ (a concatenação de todas as palavras de $R^*$) e definimos $R^\circ$ como a palavra em $\lbrace 0,1\rbrace$ onde $R_i^\circ = 0$ se $R_i^\bullet = 0$  e $R_i^\circ = 1$ caso contrário. Por exemplo, se $\Delta = 4$ e 

\

$R = (\emptyset, \emptyset,\emptyset,\emptyset,\emptyset,(3,4),\emptyset,\emptyset,\emptyset,(3,15))$, então teremos
\

$R^* = (0,0,0,0,0,01211,0,0,0,03221)$,
\

$R^\bullet = 000000121100003221$, e
\

$R^\circ = 00000111100001111.$

\

Observe que $R^* \rightarrow R^\bullet$ é uma injeção já que cada entrada de $R^*$ começa com um $0$ e não há outros zeros nas palavras de $R^*$. Logo $R \rightarrow R^\bullet$ também é uma injeção. Seguimos agora com algumas observações a respeito das palavras $R^\circ$ que provém de algum $R \in \mathcal{R}_t$.

\

\begin{df} Uma palavra Dyck parcial\footnote{Traduzido de \textit{partial Dyck word}.} é uma palavra $w$ no alfabeto $\lbrace 0,1\rbrace$ tal que qualquer prefixo de $w$ contém pelo menos o mesmo número de zeros do que de uns. Uma palavra Dyck de tamanho $2t$ é uma palavra Dyck parcial com $t$ zeros e $t$ uns. Uma descida\footnote{Traduzido de \textit{descent}.} em uma palavra Dyck (parcial) é uma sequência de uns consecutivos maximal.
\end{df}

\

\begin{lema} Para qualquer $R \in \mathcal{R}_t$, a palavra $R^\circ$ é uma palavra Dyck parcial com $t$ zeros e $t-r$ uns, onde $r$ é o número de elos coloridos após o passo $t$ do algoritmo. Ainda, todas as descidas de $R^\circ$ são pares e, se todos os ciclos de $G$ tem tamanho pelos menos $2l + 1$, para algum $l > 0$, então todas as descidas em $R^\circ$ têm tamanho pelo menos max(4,2$l$).
\end{lema}

\textbf{Prova:} Quando lemos $R^\circ$ da esquerda para a direita, cada $0$ corresponde a um elo ao qual foi atribuído uma cor e cada $1$ corresponde a um elo que foi descolorido. Como não se pode descolorir mais elos do que o número de elos coloridos, a primeira parte do lema segue. A segunda parte vem do fato que se todos os ciclos têm tamanho maior ou igual a $2l + 1$ então todos os ciclos pares (os que têm a possibilidade de ser coloridos com apenas duas cores pelo algoritmo) têm tamanho maior ou igual a $2l + 2$. Um ciclo colorido com apenas duas cores tem tamanho maior ou igual a 6. Logo cada descida é par e tem tamanho maior ou igual a max(4,$2l$). \textcolor{white}{a}\QEDA

\

Seja $R \in \mathcal{R}_t$. Se a palavra $R^\circ$ tem $t-r$ uns então a pré-imagem de $R^\circ$ pela função $R \rightarrow R^\circ$ tem, no máximo, cardinalidade $(\Delta - 1)^{t-r}$. Com efeito, como $R \rightarrow R^*$ e $R^* \rightarrow R^\bullet$ são injeções e cada 1 em $R^\circ$ corresponde a um elemento de $\lbrace 1,\ldots,\Delta - 1\rbrace$ em $R^\bullet$ o resultado segue.
\

Seja $\mathcal{R}_t^\circ = \lbrace R^\circ \vert R \in \mathcal{R}_t\rbrace$. A última afirmação junto com o fato que em uma palavra  $R^\circ$ o número de uns não é maior que o número de zeros (lema 3.9) implicam que $\vert \mathcal{R}_t \vert \leq (\Delta - 1)^t\vert \mathcal{R}_t^\circ\vert$. Então, pelo  lema 3.7:

\

\begin{lema} $\vert \mathcal{F}_t\vert \leq (K + 1)^m(\Delta - 1)^t\vert\mathcal{R}_t^\circ\vert$. 
\end{lema}

\

Nosso objetivo agora é contar palavras Dyck parciais que tenham as propriedades citadas no lema 3.9. Para facilitar, contaremos, na verdade, palavras Dyck com essas propriedades. O próximo lema mostra que contar esses dois objetos é quase equivalente, dado que $r$ (a diferença entre o número de zeros e o número de uns em uma palavra Dyck parcial) não é muito grande. 

\

\begin{lema} Sejam $t,r \in \mathbb{N}$, com $ r \leq t$, e $E \neq \lbrace 1\rbrace$ um conjunto não vazio de números naturais. Seja $C_{t,r,E}$ (resp. $C_{t,E}$) o número de palavras Dyck parciais com $t$ zeros, $t-r$ uns (resp. palavras Dyck com tamanho $2t$) e todas as descidas tendo tamanho em $E$. Então $C_{t,r,E} \leq C_{t+r(s-1),E}$, onde $s = min(E \setminus \lbrace 1\rbrace)$.
\end{lema}

\textbf{Prova:} Seja $\mathbf{D}_{t,r,E}$ (resp. $\mathbf{D}_{t,E})$ o conjunto das palavras Dyck parciais com $t$ zeros, $t-r$ uns (resp. palavras Dyck com tamanho $2t$) e todas as descidas tendo tamanho em $E$. Seja $\psi: \mathbf{D}_{t,r,E} \rightarrow \mathbf{D}_{t+r(s-1),E}$ a função que anexa ao final de uma palavra a palavra $(0^{s-1}1^s)^r$. É imediato que a função está bem definida e é injetiva o que termina a prova. \textcolor{white}{a}\QEDA

\

Existem várias formas de calcular cotas para $C_{t,E}$. Uma delas é construir uma bijeção com alguma estrutura bem conhecida. Nós usaremos uma bijeção com árvores planas\footnote{Uma árvore plana é uma árvore enraizada na qual uma ordem é dada para os filhos
de cada vértice. Tal nome vem do fato que  atribuir uma ordem aos filhos em uma árvore enraizada t é equivalente a fixar t ao plano.}. 

\

\begin{lema} O número $C_{t,E}$ de palavras Dyck com tamanho $2t$ e todas as descidas em $E$ é igual ao número de árvores planas com $t+1$ vértices tais que o grau de cada vértice está em $E \cup \lbrace 0\rbrace$.
\end{lema}

\textbf{Prova:} Existem bijeções entre os três seguintes objetos para qualquer $t \in \mathbb{N}$:

\

1. Árvores planas com $t + 1$ vértices tais que o grau de cada vértice está em $E \cup \lbrace 0\rbrace$;

\

2. Palavras Dyck com tamanho $2t$ nas quais o tamanho de qualquer sequência maximal de zeros consecutivos está em $E$;

\

3. Palavras Dyck de tamanho $2t$ tais que o tamanho de cada descida está em $E$.

\

A bijeção entre 1 e 2 pode ser construida da seguinte maneira: Em um passeio DFS na árvore, escreva, concatenadamente, cada vértice (exceto o último do passeio) tendo $i$ filhos como a palavra $0^i1$. \begin{samepage}É de fácil observação que a palavra obtida é uma palavra Dyck de tamanho $2t$ e com cada sequência maximal de zeros consecutivos em $E$. Também é fácil observar que essa é uma relação de bijeção. A bijeção entre 2 e 3 se dá pelo processo de pegar o espelho de uma palavra em 2 e trocar os zeros por uns e os uns por zeros. \textcolor{white}{a}\QEDA

\

\end{samepage}

Agora, usamos funções geradoras e o método descrito por Drmota em \cite{Drmota} para estimar $C_{t,E}$. Seja $X_E(z)$ a função geradora ordinária (ver apêndice B) associada ao número de árvores planas com $t+1$ vértices tais que o grau de cada vértice está em $E \cup \lbrace 0\rbrace$. Pelo lema anterior, $X_E(z) = z\sum_{t\in\mathbb{N}}C_{t,E}z^t$. Mas uma árvore plana como descrita acima é ou um único vértice ou, para algum $i \in E$, um vértice ligado a uma sequência de $i$ árvores planas tais que o grau de cada vértice está em $E \cup \lbrace 0\rbrace$. Logo, pelo método simbólico descrito no apêndice B, $X_E(z) = z(1+\sum_{i\in E}X_E(z)^i) = z\phi_E(X_E(z))$, com $\phi_E(x) = 1 + \sum_{i\in E}x^i$. O próximo lema é um corolário direto do teorema B.8. 

\

\begin{lema} Sejam $E \neq \lbrace 1\rbrace$ um conjunto não vazio de números inteiros não negativos e $\phi_E(x) := 1 + \sum_{i\in E}x^i$ tais que a equação $\phi_E(x) - x\phi_E^{'}(x) = 0$ tem uma solução $x = \tau$, com $0 < \tau < R$, onde $R$ é o raio de convergência de $\phi_E$. Então $\tau$ é a única solução da equação no intervalo aberto ($0,R$) e existe uma constante $c_E$ tal que $C_{t,E} \leq c_E\gamma^tt^{-\frac{3}{2}}$, onde $\gamma = \phi_E^{'}(\tau) = \phi_E(\tau) / \tau$.
\end{lema}

\

Agora estamos aptos a estabelecer cotas para a'($G$). Lembre que a cintura de um grafo $G$ é o tamanho de seu ciclo mais curto (se G não tem ciclos então $G$ tem cintura infinita).

\

\begin{teo} Sejam $l \in \mathbb{N}$, com $ l \geq 1$, $k = max(2,l)$. Então o polinômio $P(x) = (2k -3)x^{2k+2} + (1 - 2k)x^{2k} + x^4 - 2x^2 + 1$ tem uma única raiz $\tau$ no intervalo aberto $(0,1)$ e todo grafo com grau máximo $\Delta$ e cintura maior ou igual a $2l + 1$ tem uma coloração acíclica de elos com no máximo $\lceil (2 + \gamma)(\Delta - 1)\rceil$ cores, onde $\gamma = (\tau^{2k} - \tau^2 + 1)/(\tau - \tau^3)$.
\end{teo}

\textbf{Prova:} Seja $E = 2\mathbb{N} + 2k$. Então $\phi_E(x) = 1 + \sum_{i \in E}x^i = 1 + \frac{x^{2k}}{1 - x^2}$. Logo $\phi_E^{'}(x) = (2kx^{2k-1} - (2k-2)x^{2k+1})/(1-x^2)^2$, e a equação $\phi_E(x) - x\phi_E^{'}(x) = 0$ é equivalente a $P(x) = 0$. O raio de convergência de $\phi_E$ é 1 e como $P(0) = 1$ e $P(1) = -2$ o polinômio $P$ tem uma raiz $\tau$ no intervalo aberto $(0,1)$. Pelo  lema 3.13, essa é a única raiz em $(0,1)$. Também pelo lema 3.13, sabemos que existe uma constante $c_E$ tal que $C_{t,E} \leq c_E\gamma^tt^{-\frac{3}{2}}$, onde $\gamma = \phi_E^{'}(\tau) = \phi_E(\tau)/\tau = (\tau^{2k} - \tau^2 + 1)/(\tau - \tau^3)$.
\

Para provar o teorema, precisamos mostrar a existência de um vetor $F \in \lbrace 1,\ldots,\lceil \gamma(\Delta - 1)\rceil\rbrace^t$ tal que o algoritmo, recebendo $G$ e $F$ como entradas, devolve uma coloração acíclica dos elos de $G$. Pelo lema 3.10, $\vert \mathcal{F}_t\vert \leq (\lceil (2 + \gamma)(\Delta - 1)\rceil + 1)^m(\Delta - 1)^t\vert \mathcal{R}_t^{\circ}\vert$. Observe que, para cada $R \in \mathcal{R}_t$, o número de zeros e uns em cada prefixo de $R^\circ$ difere de no máximo $m - 1$ (onde $m$ é o número de elos de $G$) já que no máximo $m - 1$ elos estão coloridos a cada passo do algoritmo. Pelos lemas 3.9 e 3.11, isso implica que $\vert \mathcal{R}_t^{\circ}\vert \leq \sum_{r=0}^{m-1}C_{t+r(2k-1),E} \leq c_E^{'}\gamma^{t + m(2k-1)}t^{-\frac{3}{2}}$, onde $c_E^{'} = c_E/(\gamma^{2k-1} - 1)$. Isso implica que $\vert \mathcal{F}_t\vert \leq c_E^{'}(\lceil(2+\gamma)(\Delta - 1)\rceil + 1)^m(\Delta - 1)^t\gamma^{t+m(2k-1)}t^{-\frac{3}{2}}$, e $\vert \mathcal{F}_t\vert/ \lceil \gamma(\Delta - 1)\rceil^t$ tende a zero quando $t$ tende a infinito. Em particular, para $t$ grande o suficiente $\vert \mathcal{F}_t\vert <  \lceil \gamma(\Delta - 1)\rceil^t$, o que significa que existe um vetor $F$ para o qual o algoritmo para em menos de $t$ passos e devolve uma coloração acíclica dos elos de $G$ com no máximo $\lceil(2+\gamma)(\Delta - 1)\rceil$ cores.\textcolor{white}{a}\QEDA

\

Muthu \textit{et al.} em \cite{Muthu} provaram, em 2007, que grafos com grau máximo não maior que $\Delta$ e cintura no mínimo 9 têm uma coloração acíclica de elos com no máximo 6$\Delta$ cores e grafos com cintura no mínimo 220 essa cota foi melhorada para 4.52$\Delta$. Ndreca \textit{et al.} em \cite{Ndreca} mostraram as seguintes cotas para grafos com grau máximo $\Delta$ e cintura $g$: a'($G$) $\leq \lceil 9.62(\Delta - 1)\rceil$, a'($G$) $\leq \lceil 6.42(\Delta - 1)\rceil$, se $g \geq 5$, a'($G$) $\leq \lceil 5.77(\Delta - 1)\rceil$, se $g \geq 7$ e a'($G$) $\leq \lceil 4.52(\Delta - 1)\rceil$, se $g \geq 53$. O próximo corolário segue diretamente do teorema 3.12 e melhora significativamente todas essas cotas.

\

\begin{cor} Seja $G$ um grafo de grau máximo $\Delta$ e cintura $g$. Então
\end{cor}

\

1. a'($G$) $\leq 4\Delta - 4$; 

\

2. se $g \geq 7$, a'($G$) $\leq \lceil 3.74(\Delta - 1)\rceil$;

\

3. se $g \geq 53$, a'($G$) $\leq \lceil 3.14(\Delta - 1)\rceil$;

\

4. se $g \geq 220$, a'($G$) $\leq \lceil 3.05(\Delta - 1)\rceil$.

\

A tabela 3.1 mostra os valores de $\tau$ e $\gamma$ bem como $E$ e $P(x)$ referentes aos casos mencionados no corolário 3.13.

\

\begin{table}[h]
\centering

\begin{tabular}{|c|c|c|c|c|}

\hline 
$g$ & $E$ & $P(x)$ & $\tau$ & $\gamma$ \\ 
\hline                               
3   & 2$\mathbb{N} + 4$   & $x^6 - 2x^4 - 2x^2 + 1$ & $\frac{1}{2}(\sqrt{5} - 1)$ & 2\\
\hline
7   & 2$\mathbb{N} + 6$   & $3x^8 - 5x^6 + x^4 - 2x^2 + 1$ & 0.66336 & 1.73688\\
\hline
53  & 2$\mathbb{N} + 52$  & $49x^{54} - 51x^{52} + x^4 - 2x^2 + 1$ & 0.89610 & 1.13481\\
\hline
220 & 2$\mathbb{N} + 218$ & $215x^{220} - 217x^{218} + x^4 - 2x^2 + 1$ & 0.96341 & 1.04225 \\ 
\hline 

\end{tabular}

\vspace{0.2cm}
\caption{Valores referentes ao corolário 3.15.}

\end{table}

\subsection{O Caso Geral de \cite{Entropy}}

\

O método apresentado na seção anterior pode ser aplicado a qualquer problema de coloração de vértices (ou elos) que:

\renewcommand{\labelenumi}{\roman{enumi}}

\begin{enumerate}

\item possa ser definido como colorações que evitam um conjunto de configurações proibidas de cores -- uma configuração é um par ($H,c$) onde $H$ é um grafo  com uma coloração $c$.

\item para cada configuração ruim ($H,c$), para cada vértice $v$ de $H$ exista um número $k$ de vértices diferentes de $v$ em $H$ para os quais se soubermos suas cores existe uma única extensão dessa coloração parcial de $H$ para uma coloração congruente a $c$.

\end{enumerate}

Seja o conjunto das configurações proibidas indexado por um conjunto $I \in \mathbb{N}$. O objetivo é encontrar uma coloração $C$ do grafo $G$ tal que, para qualquer $i \in I$, e qualquer cópia $H$ de um grafo $H_i$ em $G$, a restrição de $C$ a $H$ não seja congruente a $c_i$ (duas colorações do mesmo grafo são ditas congruentes se uma pode ser obtida da outra por uma simples permutação de nomes das cores). 
\

Ainda, para cada $i$, sejam $l_i = \vert V(H_i)\vert - k_i$ e $E = \lbrace l \in \mathbb{N} \vert \exists \textit{ } i, l_i = l\rbrace$. Para $l \in E$, seja $d_l$ o máximo, entre os vértices $v$ de $G$, do número de subgrafos que contém o vértice $v$ e são isomorfos a algum $H_i$ com $l_i = l$. Seja $\gamma$ definido como no lema 3.13 da seção anterior usando esse conjunto $E$. Analisando da mesma forma que anteriormente, podemos provar que existe uma coloração de $G$ com $\gamma\cdot \sup_{l\in E}d_l^{1/l}$ cores, tal que nenhuma cópia de $H_i$ tem uma coloração congruente a $c_i$, para qualquer $i$.

\subsubsection{Prova do Caso Geral}

\

Considere um problema de coloração de vértices de $V$ de um grafo $G$ com conjunto $\mathcal{H} = \{ H_1,\ldots,H_n\}$ de configurações proibidas. Para cada $H_i$, sejam $k_i, l_i, E, d_l$ e $\gamma$ como definidos acima.

\subsubsection{O Algoritmo}

\

Considere o seguinte algoritmo: colete $t$ amostras aleatórias da distribuição discreta Uniforme$( \lceil \gamma \sup_{l \in E} d_l^{1/l}\rceil)$ e armazene-as no vetor $F$. Como anteriormente, esse vetor fornecerá as entradas do algoritmo. Estabeleça uma ordem para o conjunto de vértices $V$. Ao $i$-ésimo passo do algoritmo, atribua a $i$-ésima entrada do vetor $F$ ao vértice $v_j$ descolorido de menor índice. Se, ao fazê-lo, todos os vértices de $V$ ficam coloridos e nenhuma configuração proibida é criada, o algoritmo pára. Caso contrário, se nenhuma configuração proibida é criada mas ainda restam vértices a serem coloridos, o algoritmo segue para o próximo passo. Finalmente, se ao colorir o vértice no $i$-ésimo passo uma configuração ruim $H_s$ é criada, descolore-se os $l_s$ vértices de $H_s$. Ainda, define-se uma ordem aos no máximo $d_{l_s}$ subgrafos de $G$ que contém $v_j$ e são isomorfos a algum $H_r$ com $l_r = l_s$ e atribui-se a $i$-ésima entrada de um vetor $R$ um par contendo $l_s$ e a ordem do subgrafo observado em tal passo. O algoritmo então segue para o próximo passo.

\

\subsubsection{Análise do Algoritmo}

\

Seja $X_i$ o conjunto dos vértices não coloridos após o passo $i$ e $\phi_i$ a coloração parcial de $G$ depois do passo $i$. Seja $F$ um vetor para o qual o algoritmo retorna o par ($R$,$\phi_t$). Os próximos dois lemas mostram que o par ($R$,$\phi_t$) determina o vetor $F$.

\

\begin{lema} A cada passo $i$, o conjunto $X_i$ é unicamente determinado pelo vetor $(R_j)_{j\leq i}$.
\end{lema}

\textbf{Prova:} Provaremos por indução em $i$. É imediato que $X_1$ é unicamente determinado por $(R_j)_{j\leq 1}$ pois $X_1 = V \setminus \lbrace v_1 \rbrace$. \

Agora, supomos que $X_{i-1}$ é unicamente determinado por $(R_j)_{j\leq i-1}$. Ao passo $i$, será colorido o vértice $v_j$ onde $j = min( j : v_j \in X_{i-1})$. Se $R_i$ é vazio então $X_i = X_{i-1} \setminus \{ v_j \}$. Caso contrário, com as entradas de $R_i$ nós sabemos exatamente qual configuração proibida foi gerada e em qual subgrafo de $G$. Portanto sabemos quais vértices foram descoloridos. Logo $X_i$ também está unicamente determinado. \textcolor{white}{a}\QEDA

\

\begin{lema}  A cada passo $i$, a função que, a cada vetor $(F_j)_{j\leq i}$, associa o par $((R_j)_{j\leq i},\phi_i)$ é injetiva. 
\end{lema}

\textbf{Prova:} Novamente o argumento será por indução. Após o primeiro passo, $\phi_1$ é a configuração em que apenas o primeiro vértice está colorido e com a cor $F_1$. Supondo, para $i > 1$, que o lema é verdadeiro para $i - 1$ provaremos que ele também vale para $i$. Pelo lema 3.16 nós conhecemos $X_i$ e $X_{i-1}$. Em particular nós conhecemos o vértice $v_j$ que foi colorido no passo $i$.
\

Se $R_i$ é vazio então $\phi_{i-1}$ é obtido através de $\phi_i$ apenas descolorindo o vértice $v_j$. Logo, por indução, nós conhecemos o vetor $(F_j)_{j < i}$ bastando determinar a entrada $F_i$. Claramente, $F_i$ é a cor atribuída a $v_j$ em $\phi_i$.
\

\begin{samepage}

Agora, seja $R_i = (l,q)$. Logo nós sabemos qual configuração proibida $H_s$ envolvendo $v_j$ foi criada no passo $i$. Nesse caso, $\phi_{i-1}$ é obtido a partir de $\phi_i$ apenas colorindo os $l_s$ vértices descoloridos no passo $i$ da única forma que gera $H_s$. Pelo mesmo raciocínio acima, $(F_j)_{j\leq i-1}$ está determinado pela hipótese de indução e $F_i$ está determinado pois conhecemos a cor que o vértice $v_j$ recebeu no passo $i$ (pois só há uma maneira de colorir os $l_s$ vértices de forma a gerar $H_s$). \textcolor{white}{a}\QEDA

\

\end{samepage}

Seja $\mathcal{F}_t$ o conjunto dos vetores $F$ tais que, ao passo $t$ do algoritmo, o grafo $G$ não está completamente colorido. Seja $\mathcal{R}_t$ o conjunto dos vetores $R$ que são gerados pelo vetores de $\mathcal{F}_t$. Como há no máximo $( \lceil \gamma \sup_{l \in E} d_l^{1/l}\rceil + 1)^{\vert V \vert}$ possibilidades de colorações parciais $\phi_t$ de $G$, os dois lemas anteriores implicam que

\

\begin{lema} $\vert\mathcal{F}_t \vert \leq ( \lceil \gamma \sup_{l \in E} d_l^{1/l}\rceil + 1)^{\vert V \vert}\vert\mathcal{R}_t \vert$. 
\end{lema}

\

O desafio agora é achar uma cota para $\vert\mathcal{R}_t \vert$ a fim de mostrar que, para $t$ grande o suficiente, $\vert\mathcal{F}_t \vert$ é menor que o conjunto de todos os vetores $F$ possíveis implicando que existe um vetor $F$ para o qual o algoritmo pára. Novamente, transformaremos $\vert\mathcal{R}_t \vert$ em um conjunto mais familiar:

\

Considere uma palavra $w = w_1,\ldots,w_l$ no alfabeto $\mathbf{A} = \lbrace 1,\ldots,d_l^{1/l}\}$, e a função $\theta_l(w) = 1 + \sum_{i=1}^{l}(w_i - 1)(d_l^{1/l})^{i-1}$. A função $\theta_l$ tem imagem igual a $\lbrace 1, \ldots, d_l\}$ e é injetiva (com prova muito similar a realizada anteriormente). 
\

Seja $R \in \mathcal{R}_t$ e $R^* = (R_i^*)_{i\leq t}$ uma sequência de $t$ palavras no alfabeto $\mathbf{A}^* = \mathbf{A} \cup \lbrace 0 \rbrace$ assim definida: para todo $1 \leq i \leq t$, se $R_i$ é vazio então $R_i^* = 0$. Caso contrário, seja $R = (l,q)$ com $1 \leq q \leq d_l$. Então $R_i^*$ será a concatenação de $0$ com $\theta_l^{-1}(q)$. Agora consideramos a sequência de palavras $R^*$ como uma palavra $R^\bullet$ (a concatenação de todas as palavras de $R^*$) e definimos $R^\circ$ como a palavra em $\lbrace 0,1\rbrace$ onde $R_i^\circ = 0$ se $R_i^\bullet = 0$  e $R_i^\circ = 1$ caso contrário. Observe que $R^* \rightarrow R^\bullet$ é uma injeção já que cada entrada de $R^*$ começa com um $0$ e não há outros zeros nas palavras de $R^*$. Logo $R \rightarrow R^\bullet$ também é uma injeção. Usaremos novamente o conceito de palavras Dyck para analisar $\vert R^\circ\vert$.

\

\begin{lema} Para qualquer $R \in \mathcal{R}_t$, a palavra $R^\circ$ é uma palavra Dyck parcial com $t$ zeros e $t-r$ uns, onde $r$ é o número de vértices coloridos após o passo $t$ do algoritmo. Ainda, $E$ é o conjunto dos tamanhos de descidas de $R^\circ$.
\end{lema}

\textbf{Prova:} Quando lemos $R^\circ$ da esquerda para a direita, cada $0$ corresponde a um vértice ao qual foi atribuído uma cor e cada $1$ corresponde a um vértice que foi descolorido. Como não se pode descolorir mais vértices do que o número de vértices coloridos, a primeira parte do lema segue. A segunda parte segue da definição do conjunto $E$. \textcolor{white}{a}\QEDA

\

Seja $R \in \mathcal{R}_t$. Se a palavra $R^\circ$ tem $t-r$ uns então a pré-imagem de $R^\circ$ pela função $R \rightarrow R^\circ$ tem, no máximo, cardinalidade $(d_l^{1/l})^{t-r}$. Com efeito, como $R \rightarrow R^*$ e $R^* \rightarrow R^\bullet$ são injeções e cada 1 em $R^\circ$ corresponde a um elemento de $\lbrace 1,\ldots,d_l^{1/l}\rbrace$ em $R^\bullet$ o resultado segue.
\

Seja $\mathcal{R}_t^\circ = \lbrace R^\circ \vert R \in \mathcal{R}_t\rbrace$. A última afirmação junto com o fato que em uma palavra  $R^\circ$ o número de uns não é maior que o número de zeros (lema 3.19) implicam que $\vert \mathcal{R}_t \vert \leq (\sup_{l \in E} d_l^{1/l})^t\vert \mathcal{R}_t^\circ\vert$. Então, pelo  lema 3.18:

\

\begin{lema} $\vert \mathcal{F}_t\vert \leq ( \lceil \gamma \sup_{l \in E} d_l^{1/l}\rceil + 1)^{\vert V\vert}(\sup_{l \in E} d_l^{1/l})^t\vert\mathcal{R}_t^\circ\vert$. 
\end{lema}

\

Nosso objetivo agora é contar palavras Dyck parciais que tenham as propriedades citadas no lema 3.19. Nós seguiremos a mesma metodologia usada anteriormente. Em particular, usaremos o lema 3.13 supondo que existe $\tau$ solução da equação $\phi_E(x) - x\phi_E^{'}(x) = 0$ para finalmente provar o 

\

\begin{teo} Existe uma coloração dos vértices de $G$ com $\lceil \gamma \sup_{l \in E} d_l^{1/l}\rceil$ cores que não contém subgrafos congruentes a alguma configuração proibida.  
\end{teo}

\textbf{Prova:} Para provar o teorema, precisamos mostrar a existência de um vetor $F \in \lbrace 1,\ldots, \gamma \sup_{l \in E} d_l^{1/l}\rbrace^t$ tal que o algoritmo, recebendo $G$ e $F$ como entradas, devolve uma coloração de $G$ que evita configurações proibidas. Pelo lema 3.20, $\vert \mathcal{F}_t\vert \leq ( \lceil \gamma \sup_{l \in E} d_l^{1/l}\rceil + 1)^{\vert V\vert}(\sup_{l \in E} d_l^{1/l})^t\vert\mathcal{R}_t^\circ\vert$. Observe que, para cada $R \in \mathcal{R}_t$, o número de zeros e uns em cada prefixo de $R^\circ$ difere de no máximo $\vert V\vert - 1$ já que no máximo $\vert V\vert - 1$ vértices estão coloridos a cada passo do algoritmo. Pelos lemas 3.19 e 3.11, isso implica que $\vert \mathcal{R}_t^{\circ}\vert \leq \sum_{r=0}^{\vert V\vert -1}C_{t+r(s-1),E} \leq c_E^{'}\gamma^{t + \vert V\vert(s-1)}t^{-\frac{3}{2}}$, onde $s = min(E \setminus \lbrace 1\rbrace)$ e $c_E^{'} = c_E/(\gamma^{s-1} - 1)$, com $c_E$ e $\gamma$ como definidos na seção anterior. Isso implica que $\vert \mathcal{F}_t\vert \leq c_E^{'}(\lceil \gamma \sup_{l \in E} d_l^{1/l}\rceil + 1)^{\vert V\vert} (\sup_{l \in E} d_l^{1/l})^t\gamma^{t+\vert V\vert(s-1)}t^{-\frac{3}{2}}$ logo $\vert \mathcal{F}_t\vert/ \lceil \gamma\sup_{l \in E} d_l^{1/l}\rceil^t$ tende a zero quando $t$ tende a infinito. Em particular, para $t$ grande o suficiente, $\vert \mathcal{F}_t\vert <  \lceil \gamma\sup_{l \in E} d_l^{1/l}\rceil^t$, o que significa que existe um vetor $F$ para o qual o algoritmo para em menos de $t$ passos e devolve uma coloração dos vértices de $G$ que evita todas as configurações proibidas com no máximo $\lceil \gamma\sup_{l \in E} d_l^{1/l}\rceil$ cores.\textcolor{white}{a}\QEDA

\

\begin{ex}

\ Coloquemos o problema de coloração acíclica de elos na linguagem que acabamos de apresentar: as configurações proibidas são $H_1 = ``\textit{um caminho com dois elos com uma única cor}"$ e, para cada $i \geq 2$, $H_i = \textit{um ciclo com 2i elos colorido com apenas duas cores}$ $\textit{(sem que ocorra}$ $H_1)$. Logo $k_1 = l_1 = 1$ e, para cada $i \geq 2$, $k_i = 2$ e $l_i = 2i-2$. Isso implica que $E = \lbrace 1\rbrace \cup 2\mathbb{N} + 2$ e $\phi_E(x) = 1 + x + \frac{x^2}{1-x^2}$. Assim obtemos $\gamma = 3.6$. Nós temos $d_1 \leq 2\Delta$ e, para cada $i \geq 2$, $d_{2i-2} \leq \Delta^{2i - 2}$. Então todo grafo com grau máximo $\Delta$ tem uma coloração acíclica de elos com $3.6\cdot 2\Delta = 7.2\Delta$ cores. Tal cota é pior do que a alcançada pelo teorema 3.12 pois nele configurações menores foram tratadas de outra maneira de forma a diminuir sua influência no resultado final. \QEDB

\end{ex}

\

\subsection{Variable Version via Compressão de Entropia} 

\

Nessa seção, desenvolvemos uma variação do método descrito em \cite{Entropy} e provamos que essa variação se aplica a todo problema na forma variable version do LLL. Na sequência, relaxamos as exigências da seção anterior não requisitando a restrição \rnum{2}. Isto é, provaremos que o método de entropy compression pode ser  aplicado a todo problema de coloração em um grafo $G = (V,E)$ tal que

\begin{enumerate}

\item existem \textit{colorações ruins} a serem evitadas. Por coloração ruim, nos referimos a um par $(Z,c)$ onde $Z \subset V$ e $c$ é uma dada coloração de $Z$.

\end{enumerate}

Provaremos que todo problema na linguagem variable version do Lema Local de Lovász se enquadra nessa categoria. Para tanto precisamos do seguinte vocabulário:

\

Seja  $G$ um grafo (ou a hipergrafo) com conjunto de vértices finito $V$ para o qual estamos buscando uma coloração que evite um conjunto de \textit{colorações ruins}. Uma coloração ruim é um par $(Z, c)$ onde $Z \subset V$ e $c$ é uma dada coloração de $Z$. Sejam:

\begin{itemize}

\item $\mathcal{Z} := \{ Z : \text{ Existe c para o qual (Z,c) é uma coloração ruim}\} \text{ indexado}$ $\text{por um conjunto finito}  I \subset \mathbb{N}$.

\item Para cada $i \in I$: 

\begin{itemize}

\item $C_i := \{c : (Z_i,c) \text{ é uma coloração ruim}\}$. 

\medskip

\item $\mathscr{Z}_i := \{ (Z_i, c) : c \in C_i\}$.

\medskip

\item Para cada $v$ de $Z_i$, escolheremos $k_i$ vértices fixos diferentes de $v$ para os quais, se soubermos suas cores, existem no máximo $m_i$ maneiras de estender essa coloração parcial a uma coloração de $Z_i$ a alguma coloração $c \in C_i$. O que sempre é possível pois trabalharemos com um número finito de cores.

\medskip

\item $l_i := \vert Z_i\vert - k_i.$

\end{itemize} 

\item $E := \{ l \in \mathbb{N} : \exists \; i \in I \text{ tal que } l_i = l\}$.

\item Para cada $l \in E$ e $v \in G$, seja $Q_l(v)$ o conjunto de subconjuntos de $V$ que contém $v$ e são iguais a algum $Z_j \in \mathcal{Z}$ com $l_j = l$. Seja ainda 

$$d_l := \max_{v \in V} \vert Q_l(v) \vert.$$ 
 
\end{itemize}

\subsubsection{O Algoritmo}

Sejam $\Psi = \{\psi_1, \ldots, \psi_n\}$ um conjunto de variáveis aleatórias discretas independentes com imagem finita e  $\mathcal{A} = \{ A_x : x \in X \}$ uma família de eventos (ruins) determinados por tais variáveis\footnote{Para aplicar o LLL, a partir desse ponto, basta construir o grafo de dependência $G$ para a família de eventos $\mathcal{A}$ e encontrar os números $(r_x)_{x \in X}$  em [0,1) tais que, para cada $x \in X$, $$ \mathbf{P}(A_x) = p_x \leq r_x \prod_{y \in \Gamma(x)}(1 - r_y).$$}. Para cada $x \in X$, seja $c(A_x)$ o conjunto de  configurações das variáveis de vbl($A_x$)\footnote{Ver \cite{MT}.} que fazem $A_x$ ocorrer. Para aplicar o entropy compression basta definirmos as classes de colorações ruins. Seja $G$ um grafo com $V(G) = \Psi$. Nós olharemos para o problema de achar uma configuração de $\Psi$ que evite todos os evento de $\mathcal{A}$ como o problema de achar uma coloração de $V(G)$ que evite todas as colorações ruins, onde dizemos que se para algum $i \in [n]$ uma variável $\phi_i$ é amostrada então o valor amostrado é a cor de $\phi_i$. Logo, $\mathcal{Z} = \mathcal{A}$, $I = X$ e $\mathscr{Z}_x = \{ (vbl(A_x),c) : c \in c(A_x) \}$. Agora, sejam $\Psi_1, \ldots, \Psi_n$  vetores que armazenam $t$ amostragens das variáveis $\psi_1, \ldots, \psi_n$, respectivamente e $x_{\psi_1},\ldots, x_{\psi_n}$ variáveis que podem assumir valores em $Im(\psi_1) \cup \{ \circ \}, \ldots, Im(\psi_n) \cup \{ \circ  \} $, respectivamente. O elemento $\circ$ é um caractere com o qual todas as variáveis $x_{\psi_1},\ldots, x_{\psi_n}$ são inicializadas. Para cada $j \in [n]$ e cada $l \in E$, defini-se uma ordem no conjunto $\{ vbl(A_x) : x \in X, \ \psi_j \in \text{vbl}(A_x),  \  l_x = l \}$ -- cuja cardinalidade é no máximo $d_l$. Para cada $x \in X$, definimos também uma ordem para as $m_x$ possíveis formas de reamostrar as $l_x$ variáveis de $vbl(A_x)$ recriando uma coloração ruim de $\mathcal{H}_x$. No $i$-ésimo passo do algoritmo, atribua a primeira entrada ainda não utilizada do vetor $\Psi_j$ à variável $x_{\psi_j}$ com valor $\circ$ de menor índice. Se, ao fazê-lo, todas as variáveis ficam diferentes de $\circ$ e nenhuma coloração ruim é criada, o algoritmo pára. Caso contrário, se nenhuma coloração ruim é criada mas ainda restam variáveis com o valor $\circ$ atribuído, o algoritmo segue para o próximo passo. Finalmente, se uma coloração ruim $(vbl(A_y),c)$ é criada, atribuí-se o valor $\circ$ às $l_y$ variáveis de $vbl(A_y)$ (note que, desta forma, a cor de $\psi_s$ é o atual valor da variável $x_s$, para cada $s \in [n]$). Ainda, considere um vetor $R$ para o qual atribui-se a $i$-ésima entrada uma terna  contendo $\alpha = l_y$, $\beta \in [d_{l_y}]$ -- que nos fala a ordem de $vbl(A_y)$ em $\{ vbl(A_y) : x \in X, \ \psi_j \in \text{vbl}(A_x), \ l_x = l_y \}$ e $\gamma \in [m_y]$ -- que nos fala quais eram os valores das $l_y$ variáveis de $vbl(A_y)$. O algoritmo então segue para o próximo passo. O vetor $R$ servirá de registro dado que será de fundamental importância para a análise do algoritmo guardar a informação de o que ocorreu durante sua execução.

\begin{ob} A diferença em relação à abordagem original de Esperet e Parreau \cite{Entropy} é a adição de uma nova informação no registro $R$: a etiqueta $\gamma$ de configurações parciais. Dado que em geral trabalharemos com um número finito de parâmetros (vértices, cores, etc.), sempre podemos inicialmente fixar uma ordem para as configurações ruins parciais de forma a evitar a restrição \rnum{2} da seção anterior evitando o número $k_{(H,c)}$ da abordagem original. \textcolor{white}{aa} \QEDB    

\end{ob}

\subsubsection{Análise do Algoritmo}

\

Seja $X_i$ o conjunto das variáveis com valor $\circ$ após o passo $i$ e $\phi_i$ a configuração das variáveis no passo $i$. Sejam $\Psi_1, \ldots, \Psi_n$ vetores para os quais o algoritmo retorna o par ($R$,$\phi_t$). Essa seção se limitará a provar os próximos dois lemas que mostram que o par ($R$,$\phi_t$) determina os vetores $\Psi_1, \ldots, \Psi_n$. 

\

\begin{lema} A cada passo $i$, o conjunto $X_i$ é unicamente determinado pelo vetor $(R_j)_{j\leq i}$.
\end{lema}

\textbf{Prova:} Provaremos por indução em $i$. É imediato que $X_1$ é unicamente determinado por $(R_j)_{j\leq 1}$ pois $X_1 = \{ x_{\psi_2},\ldots, x_{\psi_n} \}$.

Agora, supomos que $X_{i-1}$ é unicamente determinado por $(R_j)_{j\leq i-1}$. Ao passo $i$, será atribuído um valor à variável $x_{\psi_j}$, onde $j = min( j : x_{\psi_j} \in X_{i-1})$. Se $R_i$ é vazio então $X_i = X_{i-1} \setminus \{ x_{\psi_j} \}$. Caso contrário, com as entradas de $R_i$ nós sabemos exatamente qual coloração ruim foi gerada e em qual conjunto de v.a.. Portanto sabemos quais variáveis receberam o valor $\circ$. Logo $X_i$ também está unicamente determinado. \textcolor{white}{a}\QEDA 

\

\begin{lema}  $((R_j)_{j\leq i},\phi_i)$ armazena informação suficiente para concluir quais cores foram usadas em cada passo do algoritmo até o passo $i$.  
\end{lema}

\textbf{Prova:} Novamente o argumento será por indução. Após o primeiro passo, $\phi_1$ é a configuração em que apenas a primeira variável tem o valor diferente de $\circ$ e vale $(\Psi_1)_1$. Supondo, para $i > 1$, que o lema é verdadeiro para $i - 1$ provaremos que ele também vale para $i$. Pelo lema 3.23, nós conhecemos $X_{i-1}$. Em particular, nós conhecemos a variável $x_{\psi_j}$ que recebeu um valor no passo $i$. Se $R_i$ é vazio então $\phi_{i-1}$ é obtido através de $\phi_i$ apenas atribuindo o valor $\circ$ à variável $x_{\psi_j}$. Logo, pela hipótese de indução, nós sabemos quais cores foram utilizadas em cada passo até o passo $i-1$. Claramente, a cor utilizada no passo $i$ é o valor de $x_{\psi_j}$ em $\phi_i$. Agora, seja $R_i = (\alpha, \beta, \gamma)$. Então nós sabemos qual coloração ruim $(vbl(A_y),c)$ envolvendo $\psi_j$ foi criada no passo $i$ e quais eram os valores de suas $l_y$ variáveis. Nesse caso, $\phi_{i-1}$ é obtido a partir de $\phi_i$ apenas atribuindo esses valores às $l_y$ variáveis que receberam o valor $\circ$ no passo $i$. Pelo mesmo raciocínio acima, Portanto novamente conhecemos as cores utilizadas até o passo $i-1$ pela hipótese de indução. Por $\gamma$ nós conhecemos a cor atribuída à $\psi_j$ no passo $i$. \textcolor{white}{a} \QEDA

\

A partir desse ponto a análise segue de perto a análise feita na seção anterior e resulta no seguinte teorema:

\

\begin{teo} 
Seja  $G$ um grafo (ou a hipergrafo) com conjunto de vértices finito $V$ para o qual estamos buscando uma coloração que evite um conjunto de \textit{colorações ruins}. Uma coloração ruim é um par $(Z, c)$ onde $Z \subset V$ e $c$ é uma dada coloração de $Z$. Sejam:

\begin{itemize}

\item $\mathcal{Z} := \{ Z : \text{ Existe c para o qual (Z,c) é uma coloração ruim}\} \text{ indexado}$ $\text{por um conjunto finito}  I \subset \mathbb{N}$.

\item Para cada $i \in I$: 

\begin{itemize}

\item $C_i := \{c : (Z_i,c) \text{ é uma coloração ruim}\}$. 

\medskip

\item Para cada $v$ de $Z_i$, escolheremos $k_i$ vértices fixos diferentes de $v$ para os quais, se soubermos suas cores, existem no máximo $m_i$ maneiras de estender essa coloração parcial a uma coloração de $Z_i$ a alguma coloração $c \in C_i$. O que sempre é possível pois trabalharemos com um número finito de cores.

\medskip

\item $l_i := \vert Z_i\vert - k_i.$

\end{itemize} 

\item $E := \{ l \in \mathbb{N} : \exists \; i \in I \text{ tal que } l_i = l\}$.

\item Para cada $l \in E$ e $v \in G$, seja $Q_l(v)$ o conjunto de subconjuntos de $V$ que contém $v$ e são iguais a algum $Z_j \in \mathcal{Z}$ com $l_j = l$. Seja ainda $d_l := \max_{v \in V} \vert Q_l(v) \vert.$ 

\item $\phi_E(x) := 1 + \sum_{i\in E}x^i$ e $\gamma := \phi_E^{'}(\tau)$.
 
\end{itemize}

Então existe uma coloração dos vértices de $G$ usando $\lceil \gamma \sup_{i \in I} (d_{l_i} m_i)^{1/l_i} \rceil$ cores que evita todos os eventos ruins.  
\end{teo}

\subsection{O Código do Entropy Compression}

\

Como visto, o método descrito em \cite{Entropy} faz uso da teoria da informação -- em particular, da teoria de compressão de dados -- mas não explicita o código\footnote{O leitor interessado é convidado a visitar o apêndice C.} utilizado. Isso será feito nessa seção. 

\

A mensagem a ser codificada é o vetor $F$ das entradas do algoritmo\footnote{Em geral, a mensagem é composta também pela configuração inicial do sistema mas esta, em \cite{Entropy}, não é aleatória (o sistema começa com todos os vértices/elos descoloridos).}. $F$ é formado pela coleção de $t$ amostragens independentes da variável aleatória discreta $X$ com distribuição Uniforme($c$), onde $c$ é o número de cores disponíveis. Aqui fica clara a principal dificuldade na hora de definir o código: seu domínio é a imagem de $X$. O problema é que, se atribuirmos uma palavra para cada cor, cada vez que uma cor for usada pelo algoritmo ela vai ser codificada pela mesma palavra. Porém são duas situações distintas atribuir uma cor para um elo e formar ou não uma configuração proibida, logo esse evento deve ser codificado de maneira diferente (lembrando que nossa missão é definir um código univocamente decodificável). Precisaremos de uma estrutura um pouco mais complexa para definir tal código. 

\

Sem perda de generalidade, suponha que estamos colorindo os vértices de um grafo $G = \{V,E\}$. Sejam $C$ o conjunto das $c$ cores disponíveis , $\Phi$ o conjunto de todas as configurações possíveis, $\mathbf{R}$ o conjunto de todas as entradas possíveis do registro $R$ (incluindo o vazio), $D$ o alfabeto binário e $D^*$ o conjunto de todas sequências finitas de símbolos do alfabeto $D$. 
\

Considere a variável aleatória $Y: \Phi \times C \rightarrow \mathbb{R}$ que associa: o número 1 a todos os pares que não formam uma configuração proibida se a cor for associada ao vértice não colorido de menor índice da configuração deixando ainda vértices a serem coloridos e aos pares em que não há vértices não coloridos na configuração (ou seja, aos pares em que a configuração já atingiu uma coloração que evita todos os pares ruins\footnote{Incluímos os pares em que não há vértices não coloridos na configuração por uma questão de formalidade.}); a todos os pares que geram um registro distinto, associa um elemento distinto do conjunto $\mathbb{N} \setminus \{1\}$ -- chamemos o conjunto dos números utilizados para os registros de $\mathcal{R}$; a todos os pares que colorem por completo de forma distinta o conjunto $V$ evitando configurações ruins associe um elemento distinto do conjunto $\mathbb{N} \setminus \ ( \{1\} \cup \mathcal{R})$. Digamos, por simplicidade, que $Im(Y) = [K]$, para algum $K \in \mathbb{N}$. 
\

Agora, usando as informações do problema de coloração em questão, devemos calcular as probabilidades da v.a. $Y$ assumir cada valor de $[K]$. Chamemos esses valores de $p_i$, para cada $i \in [K]$. Conforme exposto na seção C.3, para definir o melhor código univocamente decodificável para $Y$ devemos achar a distribuição q($Y$) 2-ádica que mais se aproxime, no sentido da entropia relativa, da distribuição de p($Y$), isto é, devemos achar os valores das probabilidades $q_i$ que minimizem

$$ D(p \vert\vert q) = \sum_{i = 1}^K p_i\log \frac{p_i}{q_i} = E_p[\log\frac{p(Y)}{q(Y)}]. $$

Agora podemos definir o código $C_Y: [K] \rightarrow D^*$, onde, para cada $i \in [K]$, os comprimentos das palavras são dados por 

$$ l_i = -\log_D q_i.$$

Com o código $C_Y$ definido, basta um pouco de formalidade para chegarmos ao nosso código: para cada $i \in [t]$, sejam $X_i$ a variável aleatória que amostrará a $i$-ésima entrada do vetor $F$ e $\phi_i$ a configuração no $i$-ésimo passo do algoritmo. Sejam $C_{X_i}$ os códigos definidos por

$$ C_{X_i}: [c] \rightarrow D^*$$
$$   \qquad\qquad\qquad x  \mapsto  C_Y(\phi_i,x) \  $$   

O código desejado é dado por

\begin{samepage} 

$$ C_F:  [c] \times \stackrel{\text{t vezes}}{\overbrace{\ldots}} \times [c] \to D^* \times \stackrel{\text{t vezes}}{\overbrace{\ldots}} \times D^*$$
$$    \qquad\qquad\qquad  (x_1,\ldots,x_t)  \mapsto  C_{X_1}(x_1) \ldots C_{X_t}(x_t) \ \    $$

\end{samepage}

$$ \diamond \ \diamond \ \diamond $$

\chapter{Novos Resultados para $c(P_n)$}

\
\begin{minipage}[t]{0.298\textwidth}
$ \ $
\end{minipage}\begin{minipage}[t]{0.702\textwidth}

\textit{``Pois não é o esforço da memória que constitui o verdadeiro conhecimento na matemática, ele sobretudo restringe suas competências, não as aumenta. ...  fazem bem os que abdicam desse trabalho tedioso que enfraquece o espírito de invenção e de pesquisa.'' - \begin{flushright} Sylvestre François Lacroix 
\end{flushright}}

\end{minipage}

\

\

Como explicado no capítulo 2, sendo $(P,L,R)$ um plano projetivo finito de ordem $n$, o objeto de estudo desse capítulo é o número mínimo de cores necessárias -- denotado por $c(P_n)$ -- para que exista uma $c$-coloração legítima de $P_n$. Em \cite{paper}, Noga Alon e Zoltan F$\ddot{\text{u}}$redi impõe que se a ordem de $P_n$ for suficientemente grande (maior do que $3\cdot 10^{250}$)\footnote{Esse número é obtido do lema 2.8.}, então $5 \leqslant c(P_n) \leqslant 8$. Em outras palavras, para um $P_n$ de ordem suficientemente grande, se $c < 5$ então em toda $c$-coloração de $P_n$ existem pelo menos duas linhas com o mesmo tipo e se $c = 8$ então existe uma $c$-coloração legítima de $P_n$. 
\

No capítulo 2 mostramos que existe um conjunto $S \subset P$ e uma 8-coloração $f$: $P\setminus S \rightarrow \lbrace 1,\ldots,8 \rbrace$ na qual nenhum ponto $p$ $\in$ $P$ pertence a mais do que quatro linhas envolvidas em pares perigosos e nenhuma linha forma um par perigoso com outras duas linhas, para um $n$ grande. Fazendo pequenas alterações nos lemas de \cite{paper}, não é difícil mostrar que, a partir de um determinado $n$, existe uma $c$-coloração de $P \setminus S$ na qual nenhum ponto $p$ $\in$ $P$ pertence a mais do que $b$ linhas envolvidas em pares perigosos e nenhuma linha forma um par perigoso com outras $a$ linhas. Como ilustrado na Figura 1 isso implica que nenhum ponto pertence a mais do que $a \cdot b$ pares perigosos.

\begin{center}

\begin{tikzpicture}

 \node (2) at ( 0, 5) {}     ;
 \node (3) at ( 2, 0) {} ;
 \node (4) at ( 2.88, 5) {} ;

 \draw [ color = brown, line width=1.5] (5,0) -- (0,5)   node[ near start ,above = -1.88pt] {$
 \bullet \quad$};

 \draw [ color = brown, line width=1.5] (5,0) -- (0.7,4.3)   node[ near start ,above = 5pt] {\color{black!35!red}{$
 p_j $}};

 \draw [ color = brown, line width=1.5] (5,0) -- (4,1)   node[ near start ,above = 0pt] {$\qquad \qquad \qquad  \qquad$ \color{black!15!brown}{$\longrightarrow $ b possibilidades}};


 \draw [ , line width=1.5] (1,4) -- (0, 5)   node[ near end , below = 100pt] {\color{black!25!blue}{a possibilidades $\longleftarrow $ }};

 \draw [ color = brown, line width=1.5] (3,2) -- (0, 5)   node[ near end , above = 8pt] {$l_r $ };

 \draw [ color = black!25!blue, line width=1.5] (2,0) -- (2.88, 5)   node[ near end , above = 10pt] {$\qquad \quad l_s $ };


 \end{tikzpicture}

 \end{center}

\begin{picture}(10,10)

\put(-5,2.3){\begin{footnotesize}
$\qquad\qquad$ Figura 1: Representação dos pares perigosos que contém $p_j$
\end{footnotesize}
}

\end{picture}

\begin{picture}(10,10)

\put(-5,1.3){\begin{footnotesize}
$\qquad\qquad\qquad\quad$ nos quais $p_j$ não está na intersecção das linhas. 
\end{footnotesize}
} 

\end{picture}

\

Usaremos tais resultados combinados com o Teorema 3.26 para mostrar que, fixado um $P_n$, existe uma coloração legítima dos pontos de $P_n$ com um número limitado de cores. Seja $P_n = (P,L)$ um plano projetivo finito com $P = \{p_1, \ldots, p_{n^2+n+1}\}$ e $L = \{\ell_1, \ldots, \ell_{n^2+n+1}\}$. Fixemos um conjunto $S \subset P$ e uma coloração parcial $f \colon P\setminus S \rightarrow [\mathsf{d}]$ para a qual nenhum ponto pertence a mais do que $a\cdot b$ pares perigosos. Definiremos nossas colorações ruins da seguinte maneira:

\begin{itemize}

\item $\mathcal{Z} = \{ Z \subset P : \exists  \ \ell_{\kappa},\ell_j, \ell_{\kappa} \neq \ell_j \text{ tais que } Z = (\ell_\kappa \cup \ell_j)\backslash (\ell_\kappa \cap \ell_j)  \} \text{ indexado}$ $\text{por um conjunto finito }  I \subset \mathbb{N}$.

\item Para cada $i \in I$: 

\begin{itemize}

\item  $C_i = \{ c : Z_i = (\ell_\kappa \cup \ell_j)\backslash (\ell_\kappa \cap \ell_j) \ \text{e} \  t_{\ell_\kappa,c} = t_{\ell_j,c} \}.$
\medskip

\item Para cada ponto $p$ de $Z_i$, considere o conjunto $Y_i$ contendo $p$ e outros  $\underline{m}-1$ pontos de menor índice da linha contendo $p$, onde $\underline{m} = \argmin_m \frac{m}{m-1}(m!a \cdot b (m-1))^{(1/m)}$ (essa escolha se tornará clara mais adiante). Fixando os pontos de $Z_i \setminus Y_i$ é fácil mostrar que $m_i = \underline{m}!$.

\medskip

\item $l_i = \vert Z_i \vert - \vert Z_i \setminus Y_i\vert =  \underline{m}$.

\end{itemize} 

\item $E = \lbrace \underline{m} \rbrace$.

\item $d_{\underline{m}} = a\cdot b$.

\end{itemize}

Então o Teorema 3.26 nos permite provar o seguinte resultado:

\begin{prop} Seja $P_n$ um plano projetivo finito de ordem $n$. Existe uma coloração legítima de $P_n$ com $\max(\mathsf{d}(n, a, b), \lceil \min_m\frac{m}{m-1}(m!a\cdot b (m-1))^{1/m}\rceil)$ cores, onde $\mathsf{d}(n,a,b)$ é o número de cores necessárias para que exista uma coloração parcial $f: P\setminus S \rightarrow [\mathsf{d}(n,a,b)]$  na qual nenhum ponto pertença a mais do que $a\cdot b$ pares perigosos.
\end{prop}

\textbf{Prova:} Nosso primeiro passo é fixar uma coloração parcial $f$ de $P_n$ na qual nenhum ponto pertença a mais do que $a\cdot b$ pares perigosos usando $\mathsf{d}(n,a,b)$ cores, nós estimares esse número em um segundo momento. Como $E = \lbrace \underline{m} \rbrace$ nós temos que $\ \phi_E(x) = 1 + x^{\underline{m}} \ $ e $\phi_E(\tau) - \tau\phi_E^{'}(\tau) = 0$ sse  $\tau = (\underline{m} - 1)^{-1 / \underline{m}}$. Então $\phi_E^{'}(\tau) = \underline{m}(\underline{m}-1)^{(1-\underline{m})/\underline{m}}.$ Portanto, usando $\lceil \min_m \frac{m}{m-1}(m!a\cdot b (m-1))^{1/m}\rceil$ cores é possível estender $f$ para uma coloração $C$ na qual nenhum par ruim é criado.  \QEDA

\

A seguir cotaremos a cardinalidade do conjunto $S$ para deduzir uma relação entre $\mathsf{d}$, $n$, $a$ e $b$.

\begin{lema} A cardinalidade de $S$ é cotada superiormente por $\frac{(n^2+n+1) 11\ln n}{n+1}$.

\end{lema}

\textbf{Prova:} Como existem $n+1$ linhas passando por cada ponto de $S$ e $\vert L \vert = n^2 + n + 1$, nós temos que, em média, cada linha $\ell \in L$ terá $\frac{\vert S \vert (n+1)}{n^2+n+1}$ pontos em $S$. Portanto  $\max_{\ell \in L} \vert \ell \cap S \vert   \geqslant \frac{\vert S \vert (n+1)}{n^2+n+1}$. Como, para cada linha $\ell \in L$ nós temos que $\vert \ell \cap S \vert \leqslant 11\ln n$, então $\vert S \vert \leqslant \frac{(n^2+n+1)11\ln n}{n+1}$.\QEDA

\

Agora seja
$$K := 2^\mathsf{d} \binom{22 \log n + \mathsf{d} }{\mathsf{d}}\frac{\mathsf{d}^{\frac{\mathsf{d}}{2}}}{[2\pi(n + 1 - 11 \log n - n^{1/2})]^{\frac{\mathsf{d}-1}{2}}}.$$ 

Generalizando o lema 2.7, descobre-se que a probabilidade $P_a$ de que exista uma linha em $L$ que forme pares perigosos com outras $a + 1$ linhas (ao invés de 2) obedece
$$P_a \leq (K)^{a+1} (n^2 + n + 1)\binom{n^2 + n}{a+1}.$$

De forma similar, generalizando o lema 2.8, descobre-se que a probabilidade $P_b$ de que exista um ponto de $P$ que pertença a mais do que $b+1$ linhas envolvidas em pares perigosos com outras linhas (ao invés de 5) obedece 
$$ P_b \leq \frac{(K)^{b+1}11\log n (n^2 + n + 1)\binom{n + 1}{b+1}(n^2 + n - (b + 1))^{b+1}}{n+1}. $$

A condição que precisa ser satisfeita para que o conjunto $f$ tenha as propriedades mencionadas é
$$ 1 - ( P_a + P_b ) > 0. $$

Como visto no capítulo 2, $K$ é uma melhor estimativa do que a feita no corolário 2.5 de \cite{paper} e está provada no corolário 2.5 do capítulo 2. Derivar $P_a$ e $P_b$ de $K$ e a condição $1 - ( P_a + P_b ) > 0$ é trivial. Temos então o seguinte problema de otimização:
$$ min \quad f(a,b,m,n) = n $$
$$ \qquad\qquad\qquad s.a. \quad \begin{cases}
\frac{m}{m-1}(m!a b (m-1))^{1/m} \leq c \\
 P_a + P_b < 1

\end{cases} $$

que busca a menor ordem de um plano projetivo que pode ser colorido legitimamente com $c$ cores. 

\begin{ob} Nos últimos anos o Método de Entropy Compression tem vencido do Método Probabilístico em inúmeros problemas. O problema de otimização acima serve como altar para o inesperado matrimônio entre os dois métodos ministrado pela aplicação aqui exposta. De fato a primeira restrição do problema nasce do Método de Entropy Compression e a segunda do Método Probabilístico. \QEDB

\end{ob}

A figura a seguir mostra alguns resultados obtidos manualmente do problema acima:

 \

\begin{center} 
\begin{tikzpicture}[scale = 0.23, xscale = 1.5, every node/.style={scale=0.45},decoration={
    markings,
   mark=at position 0.5 with {\arrow{>}}}]
\draw [help lines, thick, black, ->]  (0,0) -- (23,0);
\draw [help lines, thick, black, ->] (0,0) -- (0,43);
\draw [dashed] (0,0.5) -- (2,0.5);
\node [left] at (0,0.5) {$2$};
\draw [dashed] (0,1.5) -- (3,1.5);
\node [left] at (0,1.5) {$3$};
\draw [dashed] (0,2.5) -- (4,2.5);
\node [left] at (0,2.5) {$4$};
\draw [dashed] (0,3.5) -- (5,3.5);
\node [left] at (0,3.5) {$6$};
\draw [dashed] (0,4.5) -- (6,4.5);
\node [left] at (0,4.5) {$8$};
\draw [dashed] (0,5.5) -- (7,5.5);
\node [left] at (0,5.5) {$11$};
\draw [dashed] (0,6.5) -- (8,6.5);
\node [left] at (0,6.5) {$15$};
\draw [dashed] (0,7.5) -- (9,7.5);
\node [left] at (0,7.5) {$20$};
\draw [dashed] (0,9.745) -- (18,9.745);
\node [left] at (0,9.75) {$10^5$};

\draw [dashed] (0,42) -- (2,42);
\node [left] at (0, 42) {$10^{250}$};
\draw [dashed] (0,40.4) -- (3,40.4);
\node [left] at (0, 40.4) {$ 10^{64}$};
\draw [dashed] (0,37.2) -- (2,37.2);
\node [left] at (0, 37.2) {$ 10^{54}$};
\draw [dashed] (0,34) -- (4,34);
\node [left] at (0,34) {$ 10^{37}$};
\draw [dashed] (0,31) -- (5,31);
\node [left] at (0, 31) {$ 10^{27}$};
\draw [dashed] (0,28.2) -- (3,28.2);
\node [left] at (0, 28.2) {$ 10^{25}$};
\draw [dashed] (0,25.6) -- (6,25.6);
\node [left] at (0, 25.6) {$ 10^{22}$};
\draw [dashed] (0,23.2) -- (7,23.2);
\node [left] at (0, 23.2) {$ 10^{19}$};
\draw [dashed] (0,21) -- (8,21);
\node [left] at (0, 21) {$ 10^{17}$};
\draw [dashed] (0,19) -- (9,19);
\node [left] at (0, 19) {$ 10^{15}$};
\draw [dashed] (0,17.2) -- (5,17.2);
\node [left] at (0, 17.2) {$ 10^{13}$};
\draw [dashed] (0,15.6) -- (6,15.6);
\node [left] at (0, 15.6) {$ 10^{11}$};
\node [left] at (0, 14.2) {$ 10^{9}$};
\draw [dashed] (0,14.2) -- (7,14.1);
\node [left] at (0, 13) {$ 10^{8}$};
\draw [dashed] (0,13) -- (8,13);
\node [left] at (0, 12) {$ 10^{7}$};
\draw [dashed] (0,12) -- (9,12);
\node [below] at (2,0) {$8$};
\node [below] at (3,0) {$9$};
\node [below] at (4,0) {$10$};
\node [below] at (5,0) {$11$};
\node [below] at (6,0) {$12$};
\node [below] at (7,0) {$13$};
\node [below] at (8,0) {$14$};
\node [below] at (9,0) {$15$};
\node [below] at (18,0) {$42$};
\path (2,45) coordinate (A) (2,42) coordinate (B) (3,42) coordinate (C)(3,40.4) coordinate (D) (4,40.4) coordinate (E) (4,34) coordinate (F)(5,34) coordinate (G) (5,31) coordinate (H) (6,31) coordinate (I) (6,25.6) coordinate (J) (7,25.6) coordinate (K)(7,23.2) coordinate (L) (8,23.2) coordinate (M) (8,21) coordinate (N) (9,21) coordinate (O) (9,19) coordinate (P) (12.25,19) coordinate (Q) (12.25,17.7) coordinate (R) (15.5,17.7) coordinate (S) (15.5,16.7) coordinate (T) (17.75,16.7) coordinate (U)(17.75,16) coordinate (V) (24, 16) coordinate (X)(24,45) coordinate (Y) ;
  \draw[pattern=north east lines, dashed] (A) -- (B) -- (C) -- (D) -- (E) -- (F) -- (G) -- (H) -- (I) -- (J)-- (K)-- (L)-- (M)-- (N)-- (O)-- (P)-- (Q)-- (R) -- (S) -- (T) -- (U) -- (V) -- (X) -- (Y) -- cycle;
  
  \path  (2,0) coordinate (A)(2,0.5) coordinate (B) (3,0.5) coordinate (C)(3,1.5) coordinate (D) (4,1.5) coordinate (E) (4,2.5) coordinate (F) (5,2.5) coordinate (G) (5,3.5) coordinate (H) (6,3.5) coordinate (I) (6,4.5) coordinate (J) (7,4.5) coordinate (K) (7,5.5) coordinate (L)(8,5.5) coordinate (M)(8,6.5) coordinate (N)(9,6.5) coordinate (O) (9,7.5) coordinate (P) (11.25, 7.5) coordinate (Q) (11.25,8.06) coordinate (R) (13.5,8.06) coordinate (S) (13.5, 8.62) coordinate (T) (15.75,8.62) coordinate (U) (15.75,9.18) coordinate (V) (18,9.18) coordinate (X)(18,9.75) coordinate (Y) (18, 10.31) coordinate (Z) (15.75,10.31) coordinate (a) (15.75,10.87) coordinate (b) (13.5,10.87)  coordinate (c) (13.5,11.43) coordinate (d) (11.25,11.43) coordinate (e)  (11.25,12) coordinate (f) (9,12) coordinate (g) (9,13) coordinate  (h) (8,13) coordinate (i) (8,14.2) coordinate (j) (7,14.2) coordinate (k) (7,15.6) coordinate (l) (6,15.6) coordinate (m) (6,17.2) coordinate (n) (5,17.2) coordinate (o) (5,21) coordinate (p) (4,21) coordinate (q) (4,28.2) coordinate (r) (3,28.2) coordinate (s) (3,37.2) coordinate (t) (2,37.2) coordinate (u)(2,42) coordinate (v) (24,42) coordinate (x) (24,0) coordinate (y);
  \draw[pattern=vertical lines] (A) -- (B) -- (C) -- (D) -- (E) -- (F) -- (G) -- (H) -- (I) -- (J) -- (K) -- (L) -- (M) -- (N) -- (O) -- (P) -- (Q) -- (R) -- (S) -- (T) -- (U) -- (V) -- (X) -- (Y) -- (Z) -- (a) -- (b) -- (c) -- (d) -- (e) -- (f) -- (g) -- (h) -- (i) -- (j) -- (k) -- (l) -- (m) -- (n) -- (o) -- (p) -- (q) -- (r) -- (s) -- (t) -- (u) -- (v) -- (x) -- (y) -- cycle;

\node at (13.5,43.3) {\small\textbf{Região Anterior}};
  
\node at (15.3,31) {\normalsize\textbf{Região do LLL}};
 \node at (18,5) {\small\textbf{Região do}};
 \node at (18,4.2) {\small\textbf{Entropy}}; 
 \node at (18,3.4) {\small\textbf{Compression}};  
 \node[align=center, below] at (11.5,-1)%
 {Figura 1: Novos resultados para $c(P_n)$.};

 \end{tikzpicture}
\end{center}

\newpage

\begin{ob} Exemplos onde a região obtida pelo entropy compression é contida pela região obtida pelo lema local de lovász fogem ao meu conhecimento - possivelmente por inexistência dos mesmos. Logo, o fato de isso ter ocorrido em nosso problema não é surpreendente. O que surpreende é a região essencialmente diferente encontrada: não apenas conseguimos resultados melhores via compressão de entropia como também obtivemos resultados para valores pequenos de $n$, o que era inacessível para Lovász.  \QEDB 
\end{ob}

\

\begin{ob} A figura 1 levanta a seguinte questão:

\

\textcolor{black!35!red}{\textbf{Questão:}} Pode-se colorir de forma legítima qualquer plano projetivo finito usando 8 cores?

\

\ A pergunta, a priori sem cabimento, agora se faz natural pois, fixado um número de cores, conseguimos colorir de forma legítima planos projetivos até uma determinada ordem, depois deixamos de conseguir e voltamos a conseguir após uma ordem mais elevada. Isso levanta a questão se essa região é intrínseca ao problema ou se é possível preenche-la provando que com 8 cores consegue-se colorir de forma legítima um plano projetivo de qualquer ordem. Uma possível explicação para a região pertencer ao problema é a falta de complexidade do plano projetivo quando sua ordem é baixa. Após um ganho considerável de complexidade é preciso produzir mais tipos com o mesmo número de cores conforme a ordem cresce, como elucidado no capítulo 2. Porém tal região pode ser um ponto cego do método que usamos. Na tentativa de responder essa questão, desenvolvi um algoritmo que toma como entradas a ordem $n$ de um plano projetivo e o numero $c$ de cores a serem usadas para formar uma coloração legítima. O programa retorna, após um processo de exaustão no conjunto das colorações do plano de ordem $n$ usando $c$ cores, se existe uma coloração legítima para as entradas dadas. O código em python é o seguinte:

\

\begin{lstlisting}[language=Python]
import numpy as np
import itertools as itt

def Compara_tipos(c, t_L1, t_L2):

    k = 0

    for i in range(c):

        if( t_L1[i] == t_L2[i] ):

            k += 1

    if( k == c ):
        return 1
    else:
        return 0



def Busca_Coloracao_Legitima(q,c):

    
    MI = MatrizIncidencia(q) #Essa funcao esta definida no apendice A. 

    k = 0

    for co in itt.product( range(1,c+1), repeat = q*q + q + 1 ):

        if( k == 1 ): #Se os tipos de duas linhas sao iguais, na proxima iteracao comparamos primeiramente essas duas linhas.
        
            t_L1 = np.zeros(c)

            t_L2 = np.zeros(c)

            for i in range(q*q + q + 1):

                if( MI[i][parr[0]] == 1 ):

                    t_L1[ co[i]-1 ] += 1
            
            for i in range(q*q + q + 1):

                if( MI[i][parr[1]] == 1 ):

                    t_L2[ co[i]-1 ] += 1

            if( Compara_tipos(c, t_L1, t_L2) ):
                continue


        for par in itt.combinations( range(q*q + q + 1), 2 ):

            t_L1 = np.zeros(c)

            t_L2 = np.zeros(c)

            for i in range(q*q + q + 1):

                if( MI[i][par[0]] == 1 ):

                    t_L1[ co[i]-1 ] += 1
            
            for i in range(q*q + q + 1):

                if( MI[i][par[1]] == 1 ):

                    t_L2[ co[i]-1 ] += 1

            if( Compara_tipos(c, t_L1, t_L2) ):
                k = 1
                parr = par
                break
        else:
        
        	print("Existe uma coloracao legitima!\n\n ", co,"\n\n para a seguite matriz de incidencia\n\n",MI)

            break

    else:
        
        print("Nao existe uma coloracao legitima para o plano projetivo de ordem",q,"usando",c,"cores!")

\end{lstlisting}

\

Infelizmente esse código não tem sentido prático pois seu tempo de execução é $O(c^n)$. Um caminho natural para melhorar o algoritmo é desconsiderar previamente um número considerável de colorações trivialmente não-legítimas. Dado a complexidade do plano projetivo tal tarefa se mostrou nada trivial. \textcolor{white}{a}\QEDB

\end{ob}

$$ \diamond \ \diamond \ \diamond $$

\chapter*{Conclusão}

Neste trabalho utilizamos diferentes versões do Lema Local de Lovász (LLL) \cite{livrinho, artigo Bissacot} -- contido no arsenal do \textit{Método Probabilístico} \cite{AMP} -- bem como o recente método chamado \textit{entropy compression} \cite{Entropy, Tao}, para atacar o problema de existência de colorações legítimas em planos projetivos finitos proposto por  Noga Alon e Zoltan F\"uredi em \cite{paper}. Obtemos em \cite{BL} uma variação do método sistematizado por Louis Esperet e Aline Parreau em \cite{Entropy} e provamos que esta variação pode ser aplicada em todos os problemas formulados dentro da linguagem \textit{variable version} do LLL \cite{Tao, MT},  onde se aplica o algoritmo de Moser-Tardos. Como aplicação, mostramos que todo plano projetivo com ordem maior que $10^{54}$ admite uma coloração legítima com 8 cores -- o resultado original dos autores em \cite{paper} exige uma ordem maior que $10^{250}$. Além disso, obtemos colorações legítimas com 8 cores para planos projetivos de ordem pequena. Permitindo o uso de mais cores, obtemos progressivamente melhores resultados para ordens grandes e pequenas. De fato, provamos que com 42 cores pode-se colorir legitimamente planos projetivos finitos de qualquer ordem. A prova da existência de colorações legítimas para planos de ordens grandes e pequenas e com diferentes números de cores levanta a questão se todo plano projetivo pode ser colorido de forma legítima usando-se poucas cores independente de sua ordem. Com base nos resultados obtidos, conjecturamos que pode-se colorir legitimamente qualquer plano projetivo usando-se somente 8 cores (e não 42). Os resultados obtidos estão contidos na próxima figura:

\newpage

\begin{center} 
\begin{tikzpicture}[scale = 0.45,xscale = 1.08, every node/.style={scale=0.8},decoration={
    markings,
   mark=at position 0.5 with {\arrow{>}}}]
\draw [help lines, thick, black, ->]  (0,0) -- (23,0);
\draw [help lines, thick, black, ->] (0,0) -- (0,43);
\draw [dashed] (0,0.5) -- (2,0.5);
\node [left] at (0,0.5) {$2$};
\draw [dashed] (0,1.5) -- (3,1.5);
\node [left] at (0,1.5) {$3$};
\draw [dashed] (0,2.5) -- (4,2.5);
\node [left] at (0,2.5) {$4$};
\draw [dashed] (0,3.5) -- (5,3.5);
\node [left] at (0,3.5) {$6$};
\draw [dashed] (0,4.5) -- (6,4.5);
\node [left] at (0,4.5) {$8$};
\draw [dashed] (0,5.5) -- (7,5.5);
\node [left] at (0,5.5) {$11$};
\draw [dashed] (0,6.5) -- (8,6.5);
\node [left] at (0,6.5) {$15$};
\draw [dashed] (0,7.5) -- (9,7.5);
\node [left] at (0,7.5) {$20$};
\draw [dashed] (0,9.745) -- (18,9.745);
\node [left] at (0,9.75) {$10^5$};

\draw [dashed] (0,42) -- (2,42);
\node [left] at (0, 42) {$10^{250}$};
\draw [dashed] (0,40.4) -- (3,40.4);
\node [left] at (0, 40.4) {$ 10^{64}$};
\draw [dashed] (0,37.2) -- (2,37.2);
\node [left] at (0, 37.2) {$ 10^{54}$};
\draw [dashed] (0,34) -- (4,34);
\node [left] at (0,34) {$ 10^{37}$};
\draw [dashed] (0,31) -- (5,31);
\node [left] at (0, 31) {$ 10^{27}$};
\draw [dashed] (0,28.2) -- (3,28.2);
\node [left] at (0, 28.2) {$ 10^{25}$};
\draw [dashed] (0,25.6) -- (6,25.6);
\node [left] at (0, 25.6) {$ 10^{22}$};
\draw [dashed] (0,23.2) -- (7,23.2);
\node [left] at (0, 23.2) {$ 10^{19}$};
\draw [dashed] (0,21) -- (8,21);
\node [left] at (0, 21) {$ 10^{17}$};
\draw [dashed] (0,19) -- (9,19);
\node [left] at (0, 19) {$ 10^{15}$};
\draw [dashed] (0,17.2) -- (5,17.2);
\node [left] at (0, 17.2) {$ 10^{13}$};
\draw [dashed] (0,15.6) -- (6,15.6);
\node [left] at (0, 15.6) {$ 10^{11}$};
\node [left] at (0, 14.2) {$ 10^{9}$};
\draw [dashed] (0,14.2) -- (7,14.1);
\node [left] at (0, 13) {$ 10^{8}$};
\draw [dashed] (0,13) -- (8,13);
\node [left] at (0, 12) {$ 10^{7}$};
\draw [dashed] (0,12) -- (9,12);
\node [below] at (2,0) {$8$};
\node [below] at (3,0) {$9$};
\node [below] at (4,0) {$10$};
\node [below] at (5,0) {$11$};
\node [below] at (6,0) {$12$};
\node [below] at (7,0) {$13$};
\node [below] at (8,0) {$14$};
\node [below] at (9,0) {$15$};
\node [below] at (18,0) {$42$};
\path (2,45) coordinate (A) (2,42) coordinate (B) (3,42) coordinate (C)(3,40.4) coordinate (D) (4,40.4) coordinate (E) (4,34) coordinate (F)(5,34) coordinate (G) (5,31) coordinate (H) (6,31) coordinate (I) (6,25.6) coordinate (J) (7,25.6) coordinate (K)(7,23.2) coordinate (L) (8,23.2) coordinate (M) (8,21) coordinate (N) (9,21) coordinate (O) (9,19) coordinate (P) (12.25,19) coordinate (Q) (12.25,17.7) coordinate (R) (15.5,17.7) coordinate (S) (15.5,16.7) coordinate (T) (17.75,16.7) coordinate (U)(17.75,16) coordinate (V) (24, 16) coordinate (X)(24,45) coordinate (Y) ;
  \draw[pattern=north east lines, dashed] (A) -- (B) -- (C) -- (D) -- (E) -- (F) -- (G) -- (H) -- (I) -- (J)-- (K)-- (L)-- (M)-- (N)-- (O)-- (P)-- (Q)-- (R) -- (S) -- (T) -- (U) -- (V) -- (X) -- (Y) -- cycle;
  
  \path  (2,0) coordinate (A)(2,0.5) coordinate (B) (3,0.5) coordinate (C)(3,1.5) coordinate (D) (4,1.5) coordinate (E) (4,2.5) coordinate (F) (5,2.5) coordinate (G) (5,3.5) coordinate (H) (6,3.5) coordinate (I) (6,4.5) coordinate (J) (7,4.5) coordinate (K) (7,5.5) coordinate (L)(8,5.5) coordinate (M)(8,6.5) coordinate (N)(9,6.5) coordinate (O) (9,7.5) coordinate (P) (11.25, 7.5) coordinate (Q) (11.25,8.06) coordinate (R) (13.5,8.06) coordinate (S) (13.5, 8.62) coordinate (T) (15.75,8.62) coordinate (U) (15.75,9.18) coordinate (V) (18,9.18) coordinate (X)(18,9.75) coordinate (Y) (18, 10.31) coordinate (Z) (15.75,10.31) coordinate (a) (15.75,10.87) coordinate (b) (13.5,10.87)  coordinate (c) (13.5,11.43) coordinate (d) (11.25,11.43) coordinate (e)  (11.25,12) coordinate (f) (9,12) coordinate (g) (9,13) coordinate  (h) (8,13) coordinate (i) (8,14.2) coordinate (j) (7,14.2) coordinate (k) (7,15.6) coordinate (l) (6,15.6) coordinate (m) (6,17.2) coordinate (n) (5,17.2) coordinate (o) (5,21) coordinate (p) (4,21) coordinate (q) (4,28.2) coordinate (r) (3,28.2) coordinate (s) (3,37.2) coordinate (t) (2,37.2) coordinate (u)(2,42) coordinate (v) (24,42) coordinate (x) (24,0) coordinate (y);
  \draw[pattern=vertical lines] (A) -- (B) -- (C) -- (D) -- (E) -- (F) -- (G) -- (H) -- (I) -- (J) -- (K) -- (L) -- (M) -- (N) -- (O) -- (P) -- (Q) -- (R) -- (S) -- (T) -- (U) -- (V) -- (X) -- (Y) -- (Z) -- (a) -- (b) -- (c) -- (d) -- (e) -- (f) -- (g) -- (h) -- (i) -- (j) -- (k) -- (l) -- (m) -- (n) -- (o) -- (p) -- (q) -- (r) -- (s) -- (t) -- (u) -- (v) -- (x) -- (y) -- cycle;

\node at (13.5,43.3) {\small\textbf{Região Anterior}};
  
\node at (15.3,31) {\normalsize\textbf{Região do LLL}};
 \node at (18,5) {\small\textbf{Região do}};
 \node at (18,4.2) {\small\textbf{Entropy}}; 
 \node at (18,3.4) {\small\textbf{Compression}};  
 \node[align=center, below] at (11.5,-1)%
 {Figura 1: Novos resultados para $c(P_n)$.};

 \end{tikzpicture}
\end{center}

\appendix

\chapter{Planos Projetivos Finitos}

Esse capítulo se propõe a expor de forma sucinta as definições e resultados básicos relacionados a planos projetivos finitos. Espera-se com isso municiar o leitor para a leitura e entendimento dos capítulos 2 e 4. Para um estudo mais aprofundado no tema, indico ao leitor interessado os artigos \cite{PP1} e \cite{PP2}. 

\  

Um plano projetivo pode ser obtido pela expansão do plano euclidiano por meio de dois acréscimos:

\

1. Um ponto para cada família de retas paralelas. Sendo esse a interseção entre elas. 

\

2. Uma reta que contém somente os novos pontos adicionados.

\

Note que, com o acréscimo desses pontos, não existe mais o conceito de retas paralelas. Mais rigorosamente, podemos definir um plano projetivo de forma axiomática: sejam $P$ um conjunto cujos elementos são chamados \textit{pontos}, $L$ um conjunto cujos elementos são chamados \textit{retas }ou \textit{linhas} e $R \colon P$ $\rightarrow$ $L$ uma relação entre esses dois conjuntos. A terna $(P,L,R)$ é um plano projetivo se:

\

Axioma 1 - Dados $b$, $c$ $\in$ $P$, $b$ $\neq$ $c$, existe e é único $a \in$ $L$ tal que $(b, a)$, $(c, a)$ $\in$ $R$.

\

Axioma 2 - Dados $b$, $c$ $\in$ $L$, $b$ $\neq$ $c$, existe e é único $a$ $\in$ $P$ tal que $(a, b)$, $(a, c)$ $\in$ $R$. 

\

Axioma 3 - Existem $a_1$, $a_2$, $a_3$, $a_4$ $\in$ $P$ tais que não existem $b$ $\in$ $L$ e $i$, $j$, $k$ $\in$ $ \{1, 2, 3, 4\}$ tais que $(a_i, b), (a_j, b) (a_k, b)$ $\in$ $R$, com $k$ $\neq$ $i$ $\neq$ $j$ $\neq$ $k$.

\

Duas observações são necessárias:

\

\noindent 1 - O Axioma 1 se assemelha ao primeiro axioma de Euclides que diz ``\textit{pode-se traçar uma única reta ligando quaisquer dois pontos}''. Porém os dois diferem em algo essencial: o conceito de reta. Para nós, uma reta é apenas um elemento do conjunto $L$. A reta de Euclides, ou o segmento de reta entre dois pontos, é um comprimento de dimensão 1. Tais conceitos nem fazem sentido em nossa construção. Essa diferença ficará mais clara, mais adiante, com a análise do Plano de Fano (figura 1).

\

\noindent 2 - O Axioma 2 exclui o conceito de retas paralelas.

\

Feitas as observações, uma convenção: se $(p, l)$ $\in$ $R$ diremos que $l$ passa por $p$ ou é incidente a $p$. Outras vezes, diremos que $p$ é incidente a $l$ ou $p$ pertence a $l$ (muitas vezes $l$ será visto como um conjunto). Mais uma vez, cuidado para não confundir com a noção euclidiana de ponto pertencente a uma reta.  
A relevância do Axioma 3 será comentada mais adiante. 
\

Se $P$ e $L$ são conjuntos com cardinalidade finita então $(P,L,R)$ é dito ser um plano projetivo finito. Pode-se provar  que, para um plano projetivo finito:

\

Resultado 1: Existe $n$ $\in$ $\mathbb{N}$ tal que $\vert P \vert = \vert L \vert = n^2 + n + 1$ (a esse $n$, dá-se o nome de ordem do plano projetivo finito).

\

Resultado 2: Se $l$ $\in$ $L$ então existem somente $p_1,\ldots,p_{n+1}$ $\in$ $P$ tais que $(p_i,l)$ $\in$ $R$, para todo $i$ $\in$ $\{1,...,n+1\}$ ($n$ fixado no Resultado 1).

\

Resultado 3: Seja $p$ $\in$ $P$ então existem somente $l_1,\ldots,l_{n+1}$ $\in$ $L$ tais que $(p,l_i)$ $\in$ $R$, para todo $i$ $\in$ $\{1,...,n+1\}$ ($n$ fixado no Resultado 1).

\

Daqui em diante, nos referiremos a um plano projetivo finito pela terna $(P,L,R)$.

\

\begin{ex} Um bom exemplo para nos afastar de nossa intuição geométrica euclidiana e assimilar melhor essa construção axiomática é o plano projetivo finito de ordem 2 ou Plano de Fano:

\

\setlength{\unitlength}{0.5cm}
\begin{center}

\includegraphics[scale=0.18]{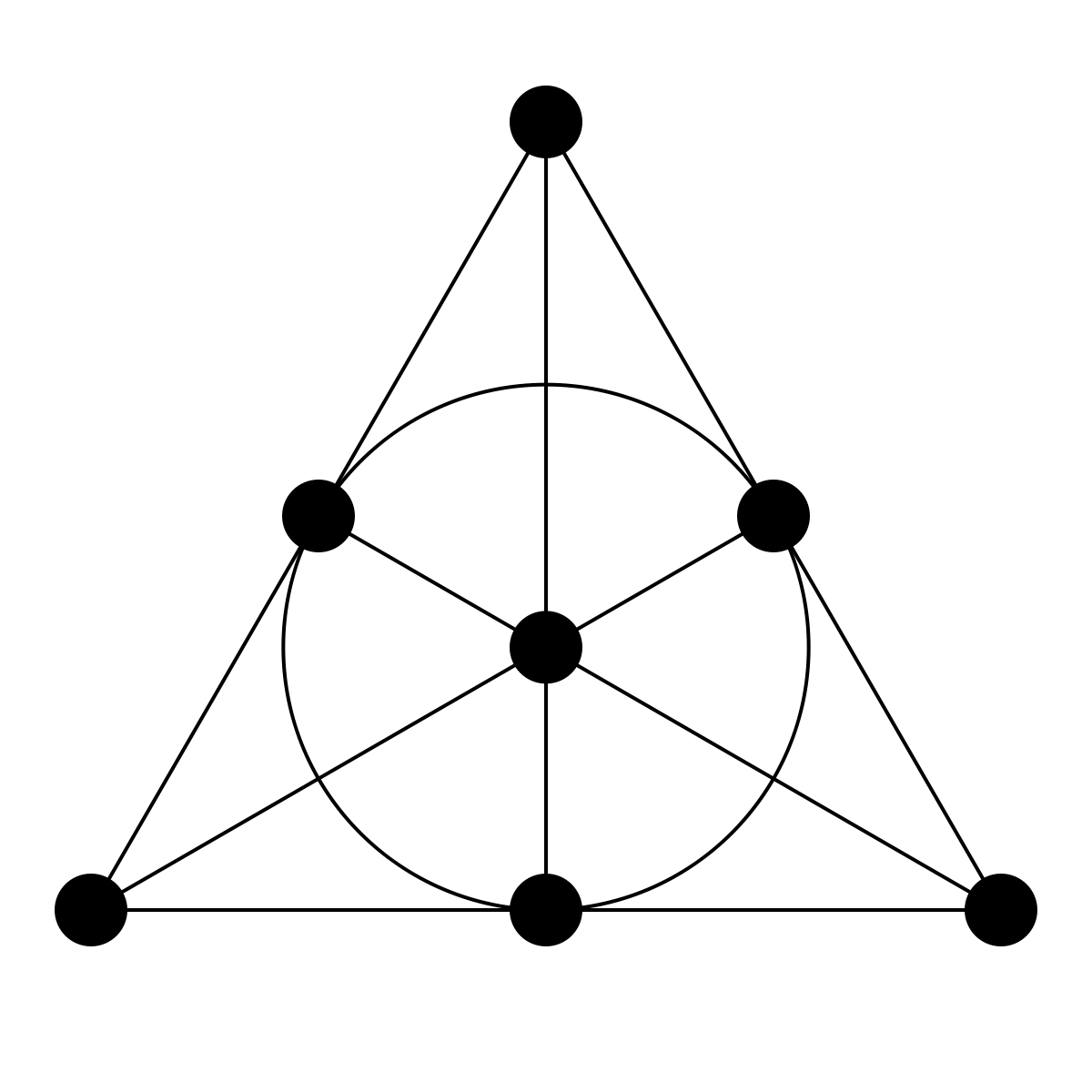} 

\end{center}

\begin{picture}(1,1)

\put(7.5,2.7){\begin{footnotesize}
Figura 1 - Plano de Fano.

\end{footnotesize}
}
\end{picture}

Repare que:

\

$\bullet$ (Axioma 1) Cada par de pontos tem exatamente um reta de incidência em comum. 

\

$\bullet$ (Axioma 2) Cada par de retas tem exatamente um ponto de incidência em comum (em nossa construção, uma reta não é um conjunto infinito de pontos. Não é porque duas retas se cruzam no desenho que deve ter um ponto incidindo com essas duas retas).

\
 
$\bullet$ (Axioma 3) Existem quatro pontos tais que nenhuma reta é incidente a mais de dois desses pontos (os três vértices mais o ponto central, por exemplo).\

$\bullet$ (Resultado 1) Como o $(P,L,R)$ é de ordem 2, há exatamente sete pontos e sete retas (você contou seis?).

\

$\bullet$ (Resultado 2) Cada ponto é incidente a três ($= n+1$) retas.\\

\

$\bullet$ (Resultado 3) Cada reta é incidente a três ($= n+1$) pontos.\textcolor{white}{a}\QEDB

\end{ex}

\

Analisar um $(P,L,R)$ por um desenho que o represente pode gerar confusões. Isso se dá pela dificuldade de desvencilhar-nos dos axiomas da geometria euclidiana. Uma abordagem, na minha opinião, mais adequada é a construção de uma matriz de incidência. 

\

\begin{ex} A figura a seguir é a representação de um $(P,L,R)$ de ordem três por sua matriz de incidência:

\

\


\setlength{\unitlength}{0.8cm}
\begin{picture}(20,15)
\linethickness{0.075mm}
\multiput(0,1.5)(1,0){15}%
{\line(0,1){14}}
\multiput(0,1.5)(0,1){15}%
{\line(1,0){14}}
\put(1.29,14.8){$l_1$}
\put(2.29,14.8){$l_2$}
\put(3.29,14.8){$l_3$}
\put(4.29,14.8){$l_4$}
\put(5.29,14.8){$l_5$}
\put(6.29,14.8){$l_6$}
\put(7.29,14.8){$l_7$}
\put(8.29,14.8){$l_8$}
\put(9.29,14.8){$l_9$}
\put(10.29,14.8){$l_{10}$}
\put(11.29,14.8){$l_{11}$}
\put(12.29,14.8){$l_{12}$}
\put(13.29,14.8){$l_{13}$}
\put(0.29,13.8){$p_1$}
\put(0.29,12.8){$p_2$}
\put(0.29,11.8){$p_3$}
\put(0.29,10.8){$p_4$}
\put(0.29,9.8){$p_5$}
\put(0.29,8.8){$p_6$}
\put(0.29,7.8){$p_7$}
\put(0.29,6.8){$p_8$}
\put(0.29,5.8){$p_9$}
\put(0.27,4.8){$p_{10}$}
\put(0.27,3.8){$p_{11}$}
\put(0.27,2.8){$p_{12}$}
\put(0.27,1.8){$p_{13}$}
\put(4.29,1.8){1}
\put(7.29,1.8){1}
\put(9.29,1.8){1}
\put(11.29,1.8){1}
\put(4.29,2.8){1}
\put(6.29,2.8){1}
\put(8.29,2.8){1}
\put(13.29,2.8){1}
\put(4.29,3.8){1}
\put(5.29,3.8){1}
\put(10.29,3.8){1}
\put(12.29,3.8){1}
\put(3.29,4.8){1}
\put(7.29,4.8){1}
\put(8.29,4.8){1}
\put(12.29,4.8){1}
\put(3.29,5.8){1}
\put(6.29,5.8){1}
\put(10.29,5.8){1}
\put(11.29,5.8){1}
\put(3.29,6.8){1}
\put(5.29,6.8){1}
\put(9.29,6.8){1}
\put(13.29,6.8){1}
\put(2.29,7.8){1}
\put(7.29,7.8){1}
\put(10.29,7.8){1}
\put(13.29,7.8){1}
\put(2.29,8.8){1}
\put(6.29,8.8){1}
\put(9.29,8.8){1}
\put(12.29,8.8){1}
\put(2.29,9.8){1}
\put(5.29,9.8){1}
\put(8.29,9.8){1}
\put(11.29,9.8){1}
\put(1.29,10.8){1}
\put(11.29,10.8){1}
\put(12.29,10.8){1}
\put(13.29,10.8){1}
\put(1.29,11.8){1}
\put(8.29,11.8){1}
\put(9.29,11.8){1}
\put(10.29,11.8){1}
\put(1.29,12.8){1}
\put(5.29,12.8){1}
\put(6.29,12.8){1}
\put(7.29,12.8){1}
\put(1.29,13.8){1}
\put(2.29,13.8){1}
\put(3.29,13.8){1}
\put(4.29,13.8){1}
\put(2.05,0.83){\begin{footnotesize}
Figura 2 - Matriz de incidência de $(P,L,R)$ de ordem 3.
\end{footnotesize}
}

\end{picture}

Onde $M_{i,j} = 1$ se $(p_i, l_j)$ $\in$ $R$ e $M_{i,j} = 0$ (aparece em branco na matriz) se $(p_i, l_j)$ $\notin$ $R$. Repare que, considerando que $p_i$ é incidente a $l_j$ se $M_{i,j} = 1$:

\

$\bullet$ (Axioma 1) Cada par de linhas tem exatamente uma coluna de incidência em comum.

\

$\bullet$ (Axioma 2) Cada par de colunas tem exatamente uma linha de incidência em comum.

\

$\bullet$ (Axioma 3) Existem quatro linhas tais que não existe coluna incidente a mais do que duas dessas linhas ($p_1$,$p_4$,$p_5$ e $p_8$, por exemplo).

\

$\bullet$ (Resultado 1) Como o $(P,L,R)$ é de ordem 3, a matriz tem exatamente 13 linhas e 13 colunas.

\

$\bullet$ (Resultado 2) Cada linha é incidente a 4 ($= n+1$) colunas.

\

$\bullet$ (Resultado 3) Cada coluna é incidente a 4 ($= n+1$) linhas.\textcolor{white}{a}\QEDB

\end{ex}

\

No capítulo 4, a observação 4.8 apresenta um código que faz uso da matriz de incidência de um plano projetivo de ordem $n$ primo. Com base em resultados obtidos em \cite{MI}, o código a seguir monta tal matriz:

\

\begin{lstlisting}[language=Python]
import numpy as np


def Posicao (a,A): #Acha a matriz posicao de A.

    B = np.zeros((len(np.transpose(A)[0]),len(A[0])*len(a)))
    for i in range(len(np.transpose(A)[0])):

        for j in range(len(A[0])*len(a)):

            if A[i][j % len(A[0])] == a[j//len(A[0])]:

                B[i][j] = 1

    return B


def MatrizA (u,q): #Cria os blocos para montar-se a matrix A mencionada no teorema 3.1 de [7].


    B = np.zeros((q,q))

    for i in range(q):

        for j in range(q):

            if -j%q == i:

                B[i][j] = q
            else:

                B[i][j] = u*(i+j)%q

    return B

def Oq (q): #Monta a matriz O_q mencionada em [7].

    B = np.zeros((q,q))

    for i in range(q):

        for j in range(q):

            B[i][j] = i + 1

    return B

def OQ (q): #Monta a matriz O_{q+1,q} mencionada em [7].

    B = np.zeros((q + 1,q))

    for i in range(q + 1):

        for j in range(q):

            B[i][j] = i + 1

    return B

def Ot (q): #Monta as ultimas colunas da matriz de incidencia.

    a = np.zeros((q,1))
    b = np.ones((q,1))

    
def MatrizIncidencia(q):

    

    B = np.zeros((q-1, q, q)) #Sera usado para armazenar os blocos da matriz A.

    P = np.zeros((q-1, q, q**2)) #Sera usado para armazenar as matrizes posicoes.

    a = list(range(1,q+1)) #Armazena os simbolos do campo de Galois.

    o = Oq(q)

    O = OQ(q)

    I = np.identity(q)

    AI = np.identity(q)

    j = np.ones(q + 1)

    h = np.zeros(q**2)

    for u in range(1,q):

        B[u-1] = MatrizA(u,q)

    for u in range(q-1):

        P[u] = Posicao(a,B[u])

    MI = P[0]

    for i in range(1, q-1):
        
        MI = np.vstack((MI,P[i]))

    for i in range(q-1):

        AI = np.hstack((AI,I))

    PO = Posicao(a,o)

    MI = np.vstack((MI,AI))

    MI = np.vstack((MI,PO))

    MI = np.vstack((MI,h))

    Ot = np.transpose(Posicao(list(range(1,q+2)), O))

    Ot = np.vstack((Ot, j))

    MI = np.hstack((MI,Ot)) #E isso completa nossa matriz de incidencia.

    return MI

\end{lstlisting}

\

O código acima, apesar de bastante didático, não é muito eficiente quanto ao uso da memória. O próximo código inverte os valores sendo em certa medida confuso mas usando apenas a memória necessária: 

\

\begin{lstlisting}[language=Python]
import numpy as np

def Posicao (a,A): #Acha a matriz posicao de A.

    B = np.zeros((len(np.transpose(A)[0]),len(A[0])*len(a)))
    for i in range(len(np.transpose(A)[0])):

        for j in range(len(A[0])*len(a)):

            if A[i][j % len(A[0])] == a[j//len(A[0])]:

                B[i][j] = 1

    return B

def PSE_MI(a, q, MI): #Monta a parte superior esquerda da matriz de incidencia.

    for i in range(q * (q - 1)):
    

        for  j in reversed(range(q * q)):
        
            if( MI[i][j % q] == a[j//q] ):
                MI[i][j] = 1 ;
            else:
                MI[i][j] = 0 ;

    return MI

def MatrizA (MI, q): #Cria os blocos para montar-se a matrix A mencionada em [7].


    for i  in range( q * (q - 1) ):

        for j in  range( q ):

            if -j % q ==  i % q :
                MI[i][j] = q 
            else:
                MI[i][j] =  ( (i//q + 1) * (i + j) ) % q

    return MI

def Oq (q): #Monta a matriz O_q mencionada em [7].

    B = np.zeros((q,q))

    for i in range(q):

        for j in range(q):

            B[i][j] = i + 1

    return B

def OQ (q): #Monta a matriz O_{q+1,q} mencionada em [7].

    B = np.zeros((q + 1,q))

    for i in range(q + 1):

        for j in range(q):

            B[i][j] = i + 1

    return B

def MatrizIncidencia(q):

    a = list(range(1,q+1)) #Armazena os simbolos do campo de Galois.

    MI = np.zeros((q*(q-1), q*q))

    o = Oq(q)

    O = OQ(q)

    I = np.identity(q)

    AI = np.identity(q)

    j = np.ones(q + 1)

    h = np.zeros(q**2)

    MI = MatrizA(MI, q)

    MI = PSE_MI(a, q, MI)
        
    for i in range(q-1):

        AI = np.hstack((AI,I))

    PO = Posicao(a,o)

    MI = np.vstack((MI,AI))

    MI = np.vstack((MI,PO))

    MI = np.vstack((MI,h))

    Ot = np.transpose(Posicao(list(range(1,q+2)), O))

    Ot = np.vstack((Ot, j))

    MI = np.hstack((MI,Ot)) #E isso completa nossa matriz de incidencia.

    return MI

              
\end{lstlisting}

\

Agora nos resta entender a relevância do axioma 3. Observe a figura a seguir:

\newpage

\begin{center}

\begin{tikzpicture}[scale=0.7]

\draw [fill = black] (-1,3) circle [radius=0.2];

\draw [ line width=1.5] (3,3) -- (7,3);   

\draw [ line width=2] (-1,0.5) -- (-1,-3.5);

\draw [ line width=1.5] (-3,0.5) -- (1,-3.5);

\draw [ line width=1.5] (1,0.5) -- (-3,-3.5);

\draw [ line width=1.5] (-3,-1.5) -- (1,-1.5);

\draw [ line width=1.5] (3,-1.5) -- (7,-1.5);

\draw [fill = black] (-1,-1.5) circle [radius=0.2];

\draw [fill = black] (3.5,-1.5) circle [radius=0.2];

\draw [fill = black] (4.5,-1.5) circle [radius=0.2];

\draw [fill = black] (5.5,-1.5) circle [radius=0.2];

\draw [fill = black] (6.5,-1.5) circle [radius=0.2];

\draw [fill = black] (-1,-7.5) circle [radius=0.2];

\draw [fill = black] (3.5,-7.5) circle [radius=0.2];

\draw [fill = black] (4.5,-7.5) circle [radius=0.2];

\draw [fill = black] (5.5,-7.5) circle [radius=0.2];

\draw [fill = black] (6.5,-7.5) circle [radius=0.2];


\draw [ line width=2] (-1,-5.5) -- (-1,-9.5);

\draw [ line width=1.5] (-3,-5.5) -- (1,-9.5);

\draw [ line width=1.5] (1,-5.5) -- (-3,-9.5);

\draw [ line width=1.5] (-3,-7.5) -- (7,-7.5);

\draw [ line width=1.5] (-1.5,-12) -- (5.5,-12);

\draw [ line width=1.5] (-0.8,-11.5) -- (2.5, -14.5);

\draw [ line width=1.5] (0.3,-11.5) -- (2.5, -14.5);

\draw [ line width=1.5] (1.4,-11.5) -- (2.5, -14.5);

\draw [ line width=1.5] (2.6,-11.5) -- (2.5, -14.5);

\draw [ line width=1.5] (3.7,-11.5) -- (2.5, -14.5) node[pos = 1, below = 20]{Figura 3: Planos projetivos degenerados.};

\draw [ line width=1.5] (4.8,-11.5) -- (2.5, -14.5);

\draw [fill = black] (-0.32,-12) circle [radius=0.2];

\draw [fill = black] (0.65,-12) circle [radius=0.2];

\draw [fill = black] (1.56,-12) circle [radius=0.2];

\draw [fill = black] (2.55,-12) circle [radius=0.2];

\draw [fill = black] (3.49,-12) circle [radius=0.2];

\draw [fill = black] (4.42,-12) circle [radius=0.2];

\draw [fill = black] (2.5,-14.4) circle [radius=0.2];









\end{tikzpicture}
 
\end{center}

\

Os seis ``planos projetivos'' acima mais o espaço vazio formam os sete tipos diferentes de planos projetivos degenerados. Eles obedecem aos dois primeiros axiomas (verifique!) mas não ao terceiro. A inclusão do axioma 3 exclui os casos degenerados. Desse forma sobram apenas os planos projetivos com estrutura mais complexa. 

\

Finalmente, repare que o resultado 1 garante que para cada $(P,L,R)$ existe um $n$ $\in$ $\mathbb{N}$ tal que $\vert P\vert = \vert L \vert = n^2 + n + 1$. Porém, não garante que para cada $n$ $\in$ $\mathbb{N}$ existe um $(P,L,R)$ de ordem $n$. Quais $n$ $\in$ $\mathbb{N}$ são ordem de um $(P,L,R)$ ainda é uma questão em aberto. Em \cite{GF}, é provado que se $n$ é uma potência de um primo então existe um $(P,L,R)$ de ordem $n$. Outro resultado importante é o obtido por Bruck e Ryser em \cite{BR} que impõe que se $n$ $=$ 1(mod4) ou $n$ = 2(mod4) então $n$ deve ser a soma de dois quadrados para ser ordem de algum $(P,L,R)$ (isso elimina  $n = 6$ como ordem possível, por exemplo). Os únicos $(P,L,R)$ conhecidos são aqueles cuja ordem é uma potência de um primo. Porém é sabido que não existe um $(P,L,R)$ de ordem 10 conforme provado em \cite{paper10}, mostrando-se, computacionalmente, que não é possível construir uma matriz de incidência de ordem 10.

\chapter{Funções Geradoras}

Como visto no capítulo 3, o método descrito em \cite{Entropy} faz uso de funções geradoras para contar o número de árvores planas com $t+1$ vértices tais que o grau de cada vértice está contigo em um conjunto $E$. Esse capítulo se baseia, em sua maior parte, em \cite{FG} e se propõe a mostrar as definições e propriedades básicas das funções geradoras ordinárias e elucidar como tal ferramenta pode ser usada em problemas de contagem. É importante ressaltar, porém, que esse método é muito mais abrangente do que esse capítulo foi capaz de mostrar. Apenas citando algumas aplicações possíveis:  na próxima seção, resolveremos relações de recorrência; em \cite{FG1} o leitor poderá ver essa técnica aplicada a problemas computacionais (como ordenação, procura, algoritmos dinâmicos, etc.);  em \cite{FG2} funções geradoras são usadas para contar grafos; e em \cite{FGA} e \cite{GJ} o leitor encontrará muitas outras aplicações não citadas nesse capítulo. O leitor interessado em uma abordagem mais completa desde a base da teoria é convidado a consultar \cite{Co}. Por sua vez, o leitor interessado em se aprofundar na questão poderá consultar \cite{GJ}.

\

Uma função geradora é um varal onde pendura-se uma sequência de números. A próxima seção procura entender o sentido dessa frase. Enquanto isso, suponha que tenhamos um problema cuja resposta é uma sequência $a_0,a_1,a_2,\ldots$ . Encontrar uma fórmula simples para $a_n$ seria a melhor maneira de resolver o problema. Se descobrirmos que $a_n = n^2 + 3$ para cada $n = 0 , 1, 2, \ldots ,$ então não há dúvidas que respondemos a questão. Mas e se não houver uma fórmula simples para o membros da sequência desconhecida? Afinal, algumas sequências são complicadas. Suponha, por exemplo, que a sequência desconhecida seja 2, 3, 5, 7, 11, 13, 17, 19, $\ldots$ , onde $a_n$ é o n-ésimo número primo. Nesse caso não podemos esperar por uma fórmula simples. Funções geradoras oferecem uma outra abordagem para o problema: ao invés de procurar por uma fórmula simples para os membros da sequência podemos encontrar uma fórmula simples para a soma da série de potência que tem tal sequência como coeficientes. 

\

\section{Relações de Recorrência}

\

\begin{ex} Considere a sequência $a_0, a_1, \ldots ,$ satisfazendo 

$$a_{n+1} = 2a_{n + 1} \ , \ n \geq 0 \ ; \ a_0 = 0.$$

Nosso desafio é descobrir quais são os elementos da sequência. Para ganharmos alguma intuição, calculemos alguns elementos primeiro. A sequência começa com 1, 3, 7, 15, 31, $\ldots$ , indicando que a sequência deve obedecer $2^n - 1, (n \geq 0)$. De fato não é difícil provar tal fato por um argumento de indução. Vejamos, porém, como é a abordagem usando funções geradoras: seja 

$$ A(x) = \sum_{n \geq 0} a_n x^n. $$ 

Para achar $A(x)$, basta multiplicar ambos os lados da relação de recorrência por $x_n$ e somar nos valores para os quais a recorrência é válida, i.e., $n \geq 0$ e tentar relacionar tais somas à função geradora $A(x)$. Do lado esquerdo ficamos com 

$$ \sum_{n \geq 0} a_{n+1} x^n. $$

Já quase encontramos $A(x)$. O índice de '$a$' é uma unidade maior do que deveria. Porém

$$ \sum_{n \geq 0} a_{n+1} x^n = a_1 + a_2 x + a_3 x^2  + \ldots  = $$

$$ = \frac{(a_0 + a_1 x + a_2 x^2 + a_3 x^3 + \ldots) - a_0}{x} = \frac{A(x)}{x}, $$

já que, nesse problema, $a_0 = 0$. O lado direito da ralação de recorrência nos fornece 

$$ \sum_{ n \geq 0} (2a_n + 1) x^n = 2A(x) + \sum_{ n \geq 0} x^n = 2A(x) + \frac{1}{1 - x}.$$

Logo, nós temos que 

$$ \frac{A(x)}{x} = 2A(x) + \frac{1}{1 - x}.$$
 
E

$$ A(x) = \frac{x}{(1 - x)(1 - 2x)}.$$

Essa é a função geradora para o problema. os elementos desconhecidos da sequências estão dispostos no varal: $a_n$ é o coeficiente de $x^n$ na expansão em série de $A(x)$. Nós temos que 

$$ \frac{x}{(1 - x)(1 - 2x)} = x (\frac{2}{1 - 2x} - \frac{1}{1-x} = (2x + 2^2 x^2 + 2^3 x^3 + \ldots ) -$$

$$- (x + x^2 + x^3 + \ldots ) = (2-1) x + (2^2 - 1) x^2 + (2^3 - 1) x^3 + \ldots .$$

E de fato temos que $a_n = 2^n - 1$, para $n \geq 0$.\textcolor{white}{a}\QEDB

\end{ex}

\

Esse exemplo gera algum desconforto. De fato, usamos um canhão para matar uma formiga. O interessante desse método é que ele é de fato um canhão! A mesma abordagem pode ser aplicada para problemas bem mais complicados. Para tanto, dada uma relação de recorrência para uma sequência $(a_n)_{n \in \mathbb{N}}$, basta seguir o seguinte método:

\

$\bullet$ Assegure-se de que o conjunto dos índices para os quais a relação é válida está bem definido.

\

$\bullet$ Dê um nome para a função geradora a ser descoberta e a escreva em função dos termos da sequência ($A(x)$, por exemplo, e  $A(x) = \sum a_n x^n$). 

\

$\bullet$ Multiplique ambos os lados da relação de recorrência por $x^n$ e some em todos os valores de $n$ para os quais a relação é válida. 

\

$\bullet$ Expresse ambos os lados da equações oriunda do passo anterior em função da função geradora $A(x)$.

\

$\bullet$ Resolva a equação na incógnita $A(x)$.

\

$\bullet$ Escolha algum método conveniente para expandir $A(x)$ em sua série de potência.

\

\begin{ex} Agora, apliquemos tal método à sequência de Fibonacci:

$$ F_{n+1} = F_n + F{_n-1} \quad , \quad (n \geq 1 ; F_0 = 0; F_1 = 1).$$ 

Chamemos nossa função geradora de $F(x) = \sum F_n x^n.$. Agora, devemos multiplicar ambos os lados da relação de recorrência por $x^n$ e somar para $n \geq 1$. O lado esquerdo nos dá

$$ F_2 x + F_3 x^2 + F_4 x^3 + \ldots = \frac{F(x) - x}{x},$$
 
e o lado direito 

$$ ( F_1 x + F_2 x^2 + F_3 x^3 + \ldots) + ( F_0 x + F_1 x^2 + F_2 x^3 + \ldots) = F(x)+ x F(x).$$

Logo 

$$ (F(x) - x)/x = F(x) + xF(x) \quad \therefore \quad F(x) = \frac{x}{1 - x - x^2}. $$

Sejam $r\pm = (1 \pm \sqrt{5})/2$. Logo

$$ 1 - x - x^2 = (1 - x r_{+})(1 - x r_{-})$$

e

$$ \frac{x}{1 - x - x^2} = \frac{x}{(1 - x r_{+})(1 - x r_{-})} = \frac{1}{r_{+} - r_{-} } ( \frac{1}{1 - x r_{+}} - \frac{1}{1 - x r_-}) = $$

$$ = \frac{1}{\sqrt{5}}( \sum_{j \geq 0} r_+^j x^j - \sum_{j \geq 0} r_-^j x^j).$$

Logo 

$$ F_n = \frac{1}{\sqrt{5}} (r_+^n - r_-^n) \ , \ n \geq 0.$$

Esse exemplo nos mostra, inclusive, um pouco mais o tipo de informação que podemos adquirir usando funções geradoras. Nós não obtivemos somente a resposta exata mas também uma aproximação, tão mais próxima da verdade quanto maior for o valor de $n$. Para valores grandes de $n$, uma ótima aproximação para $F_n$ é

$$ F_n = \frac{1}{\sqrt{5}} ( \frac{1 + \sqrt{5}}{2} )^n.$$

Ainda cabe a pergunta: de que serve uma aproximação quando se tem o valor exato? Um resposta possível é que, as vezes, a resposta exata é muito complicada e a aproximada é mais informativa. Mesmo nesse caso, em que a resposta exata não é muito complicada, ainda podemos aprender algo com a aproximação. De fato, desconsiderando $r_-$ em $F_n$ comete-se um erro inferior a 0.5. Consequentemente, $F_n$ é exatamente igual a 

$$ \lceil \frac{1}{\sqrt{5}} ( \frac{1 + \sqrt{5}}{2} )^n \rceil.$$ 

E assim obtivemos uma fórmula mais simples para $F_n$.\textcolor{white}{a}\QEDB

\end{ex}

\

No capítulo 3, aplicamos funções geradoras para contar árvores com determinada propriedade. Árvores são definidas como grafos conexos sem ciclos. Suas propriedades são básicas na teoria de grafos. Por exemplo, um grafo conexo é uma árvore se, e somente se, o número de elos é igual ao número de vértices menos 1. Ou cada par de vértices é conectado por um único caminho, etc..
\

Se marcarmos um vértice $r$ de uma árvore $T$, a transformamos em uma \textit{árvore enraizada}. Árvores enraizadas são organizadas em \textit{gerações}. A raiz ($r$) é a geração zero. Seus vizinhos formam a primeira geração e em geral vértices que distam $k$ da raiz formam a $k$-ésima geração. Se um vértice $v$ da $k$-ésima geração tem vizinhos na ($k+1$)-ésima geração, então esses vizinhos são chamados de \textit{sucessores} de $v$. Logo, se $v$ é um vértice da uma árvore enraizada $T$, então $v$ pode ser considerado como a raiz da subárvore $T_v$ composta por todos os sucessores iterado de $v$. Isso significa que uma árvore enraizada pode ser definida de forma recursiva. Essa propriedade de árvores enraizadas as tornam objetos mais fáceis de se trabalhar em problemas de contagem em comparação com árvores não-enraizadas. 
\
Existem várias aplicações de árvores enraizadas em problemas de computação. Elas aparecem naturalmente em estrutura de dados - a estruturas recursiva de pastas em um computador é basicamente uma árvore enraizada, por exemplo. Algoritmos como quicksort e o algoritmo de compressão de dados de Lampel-Zig têm uma relação direta com árvores enraizadas. Como citado, árvores enraizadas também aparecem em teoria da informação, abordada no apêndice C. Por exemplo, códigos instantâneos em um alfabeto $m$-ário são facilmente codificados como as folhas\footnote{Vértices que não possuem sucessores} de uma árvore $m$-área\footnote{Uma árvore $m$-área é  uma generalização de árvore binária onde cada vértice está ligado a no máximo $m$ sub-árvores.}. 
\

Funções geradoras são especialmente apropriadas para problemas de contagem de árvores enraizadas pois sua estrutura recursiva se traduz facilmente em uma relação de recorrência, como mostra o seguinte exemplo:

\

\begin{ex} Uma árvore binária é uma estrutura definida recursivamente como sendo um único vértice ou um vértice conectado com duas outras árvores binárias - uma sub-árvore esquerda e uma sub-árvore direita. Nosso problema é descobrir quantas árvores binárias com $n$ folhas existem. Seja $T_n$ o número de árvores binárias com $n+1$ folhas. A próxima figura representa as primeiras árvores dessa sequência:


\includegraphics[scale=0.75]{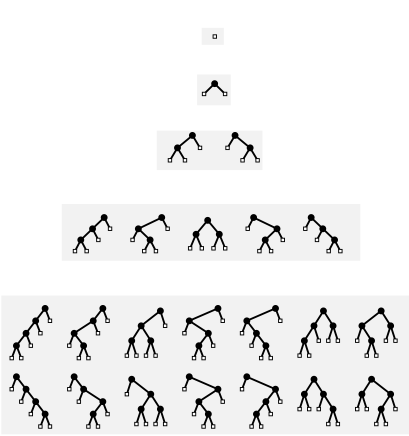}

\begin{picture}(2,2)

\put(2,-1){\begin{footnotesize}
$\qquad\qquad\quad$ Figura 1: Árvores binárias com $n \in \lbrace 1,\ldots, 5\rbrace$ folhas.

\end{footnotesize}
}
\end{picture}


\

Temos, então, $T_0 = 1$, $T_1 = 1$, $T_2 = 2$, $T_3 = 5$ e $T_4 = 14$. Calculemos $T_n$: se, em uma árvore binária com $n+1$ folhas, a sub-árvore esquerda tem $k$ folhas (existem $T_{k-1}$ sub-árvores possíveis dessa forma) então a sub-árvore direita deverá ter $n - k + 1$ (existem $T_{n - k}$ sub-árvores possíveis dessa forma). Assim, $T_n$ deve satisfazer

$$ T_n = \sum_{1\leq k\leq n}T_{k-1}T_{n-k} \ , \ \textit{para } n > 0, \ \textit{com } T_0 = 1.$$

\

Multiplicando por $z^n$ e somando em $n$ temos

$$ \sum_{n \geq 1}T_nz^n = \sum_{n \geq 1}(\sum_{1 \leq k \leq n}T_{k-1}T_{n-k})z^n  \therefore$$ 

$$  \sum_{n \geq 1}T_nz^n + T_0 - T_0 = \sum_{n \geq 1}(\sum_{0 \leq k \leq n - 1}T_{k}T_{n - k - 1})z^{n - 1}z \therefore  $$

$$  \sum_{n \geq 0}T_nz^n - T_0 = z\sum_{n - 1 \geq 0}(\sum_{0 \leq k \leq n - 1}T_{k}T_{n-k - 1})z^{n -1}.$$

Pondo  $m = n - 1$ e considerando a regra do produto de Cauchy temos

$$  \mathbf{T}(z) - T_0 = z\mathbf{T}^2(z) \therefore  \mathbf{T}(z) = z\mathbf{T}^2(z) + 1.$$

O que se reduz a uma equação do segundo grau onde 

$$\mathbf{T}(z) = \frac{1 \pm \sqrt{1 - 4z}}{2z}.$$

Se escolhermos o sinal positivo, o limite do lado direito para $z \rightarrow 0$ será $\infty$. Mas $T$(0) = 1. Se escolhermos o sinal negativo e fizermos uso da regra de L'Hospital, veremos, aliviados, que o lado direito atinge 1 no limite para $z \rightarrow 0$. Para extrair os outros coeficientes usaremos o teorema binomial generalizado para o expoente $\frac{1}{2}$:

$$ z\mathbf{T}(z) = -\frac{1}{2}\sum_{n \geq 1}\binom{\frac{1}{2}}{n}(-4z)^n.$$

E, ajustando os coeficientes, temos

$$ T_n = -\frac{1}{2}\binom{\frac{1}{2}}{n + 1}(-4)^{n+1} = -\frac{1}{2}\frac{\frac{1}{2}(\frac{1}{2} - 1)\ldots(\frac{1}{2} - n)(-4)^{n+1}}{(n + 1)!} = $$

$$ = \frac{1 \cdot 3 \cdot 5\ldots(2n - 1)\cdot2^n }{(n + 1)!} = \frac{1}{n+1}\frac{1\cdot 3\cdot 5 \ldots (2n - 1)}{n!}\frac{2\cdot 4\cdot 6\ldots 2n}{1\cdot 2\cdot 3\ldots n} = $$

$$ = \frac{1}{n + 1}\binom{2n}{n}.$$

Esses números são conhecidos como os números de Catalan.\textcolor{white}{a}\QEDB

\end{ex}

 \

Para mais exemplos (mais complicados), o leitor pode consultar \cite{FG}.

\

\section{Métodos de Obtenção de Funções Geradoras}

A seção anterior nos mostrou um método para se obter funções geradoras. Dependendo da natureza do problema, outros métodos podem ser mais adequados. Essa seção apresenta ao leitor dois outros métodos.

\subsection{O Método Simbólico}

\

O método simbólico, utilizado no capítulo 4, é uma ferramenta poderosa para traduzir definições formais de objetos combinatórios em equações de funções geradoras. Essa seção apenas dá a primeira pincelada em uma grande tela. Para uma abordagem mais completa e profunda, consulte \cite{FGA}.  

\

Quando contamos com funções geradoras, nós consideramos classes de objetos combinatórios com uma noção de tamanho definida para cada objeto. Para uma classe $\mathcal{A}$, nós denotamos o número de membros da classe de tamanho $n$  por $a_n$. E a função geradora toma a forma

$$ \mathbf{A}(z) = \sum_{n \geq 0}a_nz^n = \sum_{a \in \mathcal{A}}z^{\vert a\vert}.$$

Dadas duas classes $\mathcal{A}$ e $\mathcal{B}$ de objetos combinatórios, podemos criar outras classes como $\mathcal{A} + \mathcal{B}$, que consiste na classe composta pela união disjunta dos elementos das duas classes ou $\mathcal{A} \times \mathcal{B}$, que consiste na classe composta pelos pares ordenados de elementos de $\mathcal{A}$ e $\mathcal{B}$. Não é difícil mostrar que essas novas classes terão funções geradoras dadas por $\mathbf{A}(z) +  \mathbf{B}(z)$ e $\mathbf{A}(z) \mathbf{B}(z)$, respectivamente. 
\

O método simbólico consiste em combinar classes elementares para construir uma classe desejada. As classes elementares são: a classe vazia $\emptyset$, com função geradora $0$, o objeto nulo $\epsilon$, com função geradora $1$ e objetos atômicos como $0$ ou $1$ em problemas binários, $\circ$ representando vértices em problemas de grafos, com função geradora $z$. 

\

\begin{ex} Seja $G$ a classe de sequências binárias que não possuem dois zeros consecutivos. Tais sequências podem ser $\epsilon$, um único $0$, ou um 1 ou um 01 seguido de uma sequência sem dois zeros consecutivos. O que se traduz simbolicamente em
	
$$ G = \epsilon + \lbrace 0\rbrace + \lbrace 1,01\rbrace \times G.$$

Então

$$ \mathbf{G}(z) = 1 + z + (z + z^2)\mathbf{G}(z) \ \therefore \ \mathbf{G}(z) = \frac{1 + z}{(1 - z - z^2)}.$$

Logo, pelo resultado do exemplo B2, o número de sequências binárias de tamanho $n$ que não possuem dois zeros consecutivos obedece a sequência de Fibonacci $G_n + G_{n+1} = G_{n + 2}$.  \textcolor{white}{a}\QEDB

\end{ex}

\

\begin{ex} Mostraremos que o número $p_n$ de árvores enraizadas planares\footnote{Árvores enraizadas onde cada vértice possui um número arbitrário de sucessores com uma ordem natural da esquerda para a direita.} com $n \geq 1$ vértices é dado por 

$$ p_n = \frac{1}{n} \binom{2n - 2}{n - 1}.$$ 

Seja $\mathcal{P}$ o conjunto de árvores enraizadas planares. Logo, pela definição de árvores enraizadas planares,

$$ \mathcal{P} = \circ + \circ \times \mathcal{P} + \circ \times \mathcal{P}^2 + \circ \times \mathcal{P}^3 + \ldots .$$

Sendo $P(x) = \sum p_n x^n$ temos então

$$ P(x) = x + xP(x) + xP(x)^2 + xP(x)^3 + \ldots = \frac{x}{1 - P(x)}. $$

Logo

$$ P(x) = \frac{1 - \sqrt{1 - 4x}}{2} = xT(x).$$

Portanto

$$ p_n = T_{n-1} = \frac{1}{n}\binom{2n - 2}{n - 1}.$$\textcolor{white}{a}\QEDB

\end{ex}

\begin{ob}

A relação $p_n = T_{n-1}$ tem um significado mais profundo. Existe uma bijeção entre árvores enraizadas planares com $n$ vértices e árvores binárias com $n$ folhas. Comecemos com uma árvore enraizada planar com $n$ vértices e apliquemos o seguinte procedimento:

\

1. Delete a raiz e todos os elos envolvendo a raiz.

\

2. Se um vértice possui sucessores, delete todos os elos envolvendo o elo e os sucessores exceto o elo mais a esquerda.

\

3. Ligue todos os sucessores do passo anterior com um caminho (horizontal).

\

4. Rotacione esses novos caminhos 90º. 

\

5. Os $n-1$ vértices restantes são agora considerados vértices internos de uma árvore binária bastando anexar as $n+1$ folhas que faltam. 

\

Não é difícil verificar que esse procedimento é uma bijeção.\textcolor{white}{a}\QEDB

\end{ob}

\

\begin{ex} O problema de contar árvores binárias com $n$ folhas também pode ser atacado por esse procedimento:

\

Uma árvore binária é composta por ou apenas um vértice ou um vértice ligado a duas sub-árvores que são, por sua vez, árvores binárias. Logo, sendo $\tau$ a classe das árvores binárias, temos

$$ \tau = \epsilon + \circ \times \tau \times \tau$$

Como a condição inicial é $T_0 = 1$ temos, como esperado

$$ \mathbf{T}(z) = 1 + z\mathbf{T}^2(z).$$\textcolor{white}{a}\QEDB

\end{ex}

\subsection{Árvores Simplesmente Geradas}

\

Árvores simplesmente geradas\footnote{Traduzido de simply generated trees.} foram introduzidas por Meir e Moon em \cite{Meir} e são generalizações de vários tipos de árvores enraizadas. Nessa seção, definimos esse objeto e provamos o teorema B.8 usado no capítulo 4. Seja 

$$\phi(x) = \phi_0 + \phi_1x + \phi_2x^2 + \ldots $$

uma série de potências com nenhum coeficiente negativo. Em particular, com $\phi_0 > 0$ e $\phi_j > 0$ para algum $j \geq 2$. Seja $T$ uma árvore enraizada finita. Definimos o peso $w(T)$ dessa árvore por

$$ w(T) = \prod_{j \geq 0}\phi_j^{D_j(T)},$$

onde $D_j(T)$ é o número de vértices de $T$ com $j$ sucessores. Definindo

$$ y_n = \sum_{\vert T\vert = n}w(T), $$

então $y_n$ denota um número ponderado de árvores de tamanho $n$. Por exemplo, se $\phi_j = 1$ para todo $j \geq 0$ (ou seja, se $\phi(x) = 1/(1-x)$) então todas as árvores enraizadas tem o mesmo peso $w(T) = 1$ e $y_n = p_n$ é o número de árvores planas. Se $\phi(x) = 1 + x + x^2$ então apenas árvores enraizadas cujos vértices tem menos do que três sucessores recebem um peso diferente de zero dado por $w(T) = 1$. Tais árvores são conhecidas como árvores de Motzkin e, nesse caso, $y_n$ é o número de árvores de Motzkin de tamanho $n$. As árvores binárias também podem ser expressas dessa maneira: se $\phi(x) = 1 + 2x + x^2 = (1 + x)^2$ então os vértices com apenas um sucessor recebem peso 2. Isso se dá pois, em uma árvore binária, existem dois tipos de vértices com apenas um sucessor,  aqueles com uma sub-árvore esquerda mas sem sub-árvore direita e aqueles com uma sub-árvore direita mas sem sub-árvore esquerda. De forma similar, uma árvore $m$-ária, uma generalização da árvore binária onde cada vértice está ligado a no máximo $m$ sub-árvores, são contadas com o uso da função $\phi(x) = (1 + x)^m$.
\

Logo a função geradora

$$ y(x) = \sum_{n \geq 1}y_nx^n $$

satisfaz a equação

$$ y(x) = x\phi(y(x)).$$

Com efeito, dada estrutura recursiva das árvores simplesmente geradas, usaremos o método simbólico para atingir a equação desejada. A dificuldade consiste no fato de que não estamos trabalhando com uma estrutura fixa e sim com uma classe de objetos. Para driblar esse problema usaremos os coeficientes $\phi_i$. Pelas definições e exemplos dados acima, é possível entender a informação ``$\phi_i = a$'' como ``\textit{os vértices com $i$ filhos possuem $a$ formas de dispor seus descendentes}''. Portanto, para uma árvore simplesmente gerada $T$ temos a relação recursiva ponderada

$$ T = \phi_0\cdot\circ+\phi_1\cdot\circ\times T + \phi_2\cdot\circ\times T^2 + \cdots.$$ 

Portanto $y(x) = x\phi(y(x))$

Provaremos agora o resultado que dá sentido a esse apêndice pois sustenta o lema 3.13:

\begin{teo} Seja $R$ o raio de convergência de $\phi(t)$ e suponha que existe $\tau$ com $0 < \tau < R$ tal que $\tau\phi'(\tau) = \phi(\tau)$. Seja $d = mdc\{ j>0 : \phi_j > 0\}$. Então

$$ y_n = d \sqrt{\frac{\phi(\tau)}{2\pi \phi''(\tau)}} \frac{\phi'(\tau)^n}{n^{3/2}}(1 + O(n^{-1}) \ , \ n \ \equiv 1 ( mod \ d),$$

e $y_n = 0$ se $n \not\equiv 1 ( mod \ d)$.

\end{teo}

\textbf{Prova} Aplicaremos o teorema B.38 para $F(x,y) = x\phi(y)$. Assumindo $d=1$ as condições do teorema B.38 são satisfeitas e o resultado segue. Para $d > 1$ temos que $y_n = 0$ se $n \not\equiv 1$(mod $d$). Logo $y(x) = \widehat{y}(x^d)/x^{d-1}$ e $ \phi(x) = \widehat{\phi}(x^d)$ com $\widehat{y}$ e $\widehat{\phi}$ funções analíticas satisfazendo $\widehat{y}(x) = x\widehat{\phi}(\widehat{y}(x))$ com mdc $\widehat{d} = 1$. Logo estamos novamente nas condições do teorema B.38 e o resultado segue.  \textcolor{white}{a}\QEDA

\

\begin{ob} Em tempo, cabe ressaltar que existe uma surpreendente relação entre árvores simplesmente geradas e o processo de ramificação de Galton-Watson. O leito interessado é convidado a consultar \cite{Drmota}.\textcolor{white}{a}\QEDB

\end{ob}

\section{Séries de Potência Formais}

\

Discutir a teoria formal de série de potências, e não a teoria analítica, significa observar tais séries somente como objetos algébricos em sua função de varal, i.é., com as potências servindo apenas como índices para seus coeficientes sem se importar com a função para a qual essa série possa estar convergindo. O objetivo de se estudar séries sob esse ponto de vista é conseguir fazer manipulações como as feitas na seção B.1 sem se preocupar com questões de convergência podendo, por exemplo, calcular a derivada de uma função geradora sem saber se ela converge para uma função. Essa seção se propõe a mostrar que não há falta de rigor ao se trabalhar dessa forma. Tais manipulações podem ser feitas no anel de séries de potências formais onde não existe o conceito de convergência. Nós podemos executar o método por completo, achar a função geradora em sua forma de série de potência, e só então descobrir se tal série converge ou não. Caso não convirja, ainda haverá muita informação a ser coletada da série formal mas não teremos acesso a nenhuma informação analítica, como fórmulas assintóticas para os coeficientes por exemplo. A série

$$ f = 1 + x + 2 x^2 + 6 x^3 + 24 x^4 + 120 x^5 + \ldots + n! x^n + \ldots \ ,$$ 

por exemplo, pertence ao anel das séries de potências formais mas não converge para nenhum valor de $x$. Apesar disso, essa série desempenha um papel importante em alguns problemas de contagem.

\

\begin{df} Uma série de potência formal é uma expressão da forma 

$$ a_0 + a_1x + a_2x^2 + \ldots,$$

onde a sequência $(a_n)_{n \in \mathbb{N}}$ é chamada de sequência dos coeficientes. Ainda, dizer que duas séries de potência são iguais significa dizer que os coeficientes são iguais.

\end{df}

\

Algumas operações podem ser feitas entre séries de potências formais como adição e subtração:

$$ \sum a_n x^n \pm \sum b_n x^n = \sum (a_n \pm b_n) x^n.$$

E multiplicação, usando a regra do produto de Cauchy:
 
$$ \sum a_n x^n \sum b_n x^n = \sum c_n x^n \ , \ \textit{onde } c_n = \sum_k a_k b_{n-k}.$$

Essa regra para o produto confere uma aplicabilidade de séries de potências a muitos problemas de combinatória. De fato, é frequente conseguirmos construir todos os objetos $a_n$ de objetos do tipo $n$ de uma família escolhendo um objeto do tipo $k$ e outro do tipo $n-k$ e colocando-os juntos para forma um objeto do tipo $n$ (isso é feito no exemplo B.3). O número de maneiras de se fazer isso será $a_k a_{n-k}$, e se somarmos em $k$ nós veremos que o produto de Cauchy é de fundamental importância para o problema em questão. 
\

Usando a propriedade distributiva, nós obtemos 

$$ (1 - x)(1 + x + x^2 + x^3 + \ldots ) = 1.$$

Logo, podemos dizer que a série $(1-x)$ tem um inverso multiplicativo dado por $1 + x + x^2 + \ldots $ (e vice-versa).

\

\begin{prop} Uma série de potência formal $f = \sum a_n x^n$ tem uma inverso multiplicativo se, e somente se, $ a_0 \neq 0$. Nesse caso o inverso multiplicativo é único.
\end{prop}

\textbf{Prova:} Seja $f$ uma série de potência com inverso multiplicativo $1/f = \sum b_n x^n$. Então $f(1/f) =
1$ e $a_0 b_0 = 1.$ Logo $a_0 \neq 0$. Além disso, sabemos que, para $n \geq 1$, $\sum_k a_k b_{n-k} = 0$, de onde tiramos

$$ b_n = \frac{-1}{a_0}\sum_{k \geq 1} a_k b_{n-k}.$$

Logo $b_1, b_2,\ldots $ estão unicamente determinados. Reciprocamente, suponha que $a_0 \neq 0$. Então podemos determinar $b_0, b_1, b_2, \ldots $ como acima e a série resultante é o inverso multiplicativo de $F$. \textcolor{white}{a}\QEDA

\

As séries de potencia formais com as regras aritméticas descritas formam um anel cujos elementos invertíveis são as séries com termo independente diferente de zero. É importante ressaltar que o conceito de inverso multiplicativo de uma série difere do conceito de inverso de uma série. O inverso de uma série $f$, se existir, é uma série $g$ que obedece $f(g(x)) = g(f(x)) = x$, para todo $x$. Vejamos quando uma série possui um inverso.

\

\begin{df} se $ f = \sum a_n x^n$ então $f(g(x)) = \sum a_n g(x)^n$.
\end{df} 

\

Se a série $g(x)$ tiver elemento independente igual a zero, então nós sabemos calcular os coeficientes de $f(g(x))$. De fato, se quisermos, para algum $k$, calcular o coeficiente de $x^k$, podemos observar que em

$$ a_n g(x)^n = a_n(g_1 x + g_2 x^2 + \ldots)^n = a_n x^n(g_1 + g_2x + \ldots)^n $$

se $n > k$ só existirão potências de $x$ maiores que $k$ e não precisamos nos preocupar com eles. Logo teremos um processo bem definido (com número finito de passos) para calcular o coeficiente de $x^k$. Se $g_0 \neq 0$, o mesmo não pode ser feito e, se não for polinomial, o processo para calcular cada coeficiente é infinito e dependente da convergência da série para estar bem definido. Como no anel de séries formais não existe o conceito de convergência, $f(g(x))$ só está definido se $f$ for polinomial ou se $g_0 = 0$. 
\

No anel das séries de potência formais existem outras operações análogas às operações definidas no cálculo (porém sem o conceito de limite). Por exemplo, a derivada de uma série de potência formal $f = \sum a_n x^n$ é a série $f' = \sum na_n x^{n-1}$. A derivação segue as regras usuais do calculo tais como soma, produto e cociente. Muitas delas, inclusive, são mais fáceis de provar nesse contexto. Por exemplo:

\

\begin{prop} Se $f' = 0$ então $f = a_0$ é constante.
\end{prop}

\textbf{Prova:} $ f' = 0$ significa que a série de potência formal $f'$ é idêntica à série de potência formal 0, o que significa que todos os coeficientes de $f'$ são iguais á zero. Mas os coeficientes de $f'$ são $a_1, 2a_2, 3a_3, \ldots$ . Logo $a_j = 0$ para todo $j \geq 1$ e $f = a_0$ é constante. \textcolor{white}{a}\QEDA

\

\begin{prop} Se $f' = f$ então $ f = ce^x$.
\end{prop}

\textbf{Prova:} Como $f' = f$, os coeficientes de $x^n$ devem ser os mesmo em $f$ e em $f'$, para todo $n \geq 0$. Logo $(n + 1)a_{n+1} = a_n$, para todo $n \geq 0$ e $a_{n+1} = a_n/(n + 1), (n \geq 0)$. Por indução em $n$, $a_n = a_0/n!$, para todo $n \geq 0$. Portanto $f = a_0e^x$.\textcolor{white}{a}\QEDA

\

\section{Funções Geradoras Ordinárias}

\

Operações com séries formais implicam em operações correspondentes em seus coeficientes. Se a série convergir para uma função, então operações com tal função correspondem a operações nos coeficientes da série de potência obtida ao expandir a função. Nessa seção, nós exploraremos algumas dessas relações cujo entendimento facilita o manuseio e entendimento de como se usar as funções geradoras. Usaremos a seguinte definição:

\

\begin{df} O símbolo $f \stackrel{\text{fgo}}{\longleftrightarrow} (a_n)_{n \in \mathbb{N}}$ significa que $f$ é a função geradora ordinária para a sequência $(a_n)_{n \in \mathbb{N}}$, i.e., que

$$f = \sum a_n x^n.$$

\end{df}

\subsection{Propriedades das Funções Geradoras Ordinárias}

\

Sejam $f \stackrel{\text{fgo}}{\longleftrightarrow} (a_n)_{n \in \mathbb{N}}$. Procuraremos agora a função que gera $(a_{n + 1})_{n \in \mathbb{N}}$. Para tanto, precisamos fazer alguns cálculos:

$$ \sum a_{n+1} x^n = \frac{1}{x} \sum_{m \geq 1} a_m x^m = \frac{f(x) - f(0)}{x}.$$

Então 

$$f \stackrel{\text{fgo}}{\longleftrightarrow} (a_n)_{n \in \mathbb{N}} \Rightarrow (\frac{f - a_0}{x}) \stackrel{\text{fgo}}{\longleftrightarrow} (a_{n+1})_{n \in \mathbb{N}}.$$

Se transladarmos os índices em duas unidades apenas precisaremos repetir o cálculo feito acima para concluir que

$$ (a_{n + 2})_{n \in \mathbb{N}} \stackrel{\text{fgo}}{\longleftrightarrow} \frac{((f - a_0)/x) - a_1}{x} = \frac{f - a_0 - a_1x}{x^2}.$$

Note que esse fato facilita entender a relação de recorrência de Fibonacci $F_{n+2} = F_{n+1} + F_n$, para $n \geq 0$, $F_0 = 0$ e  $F_1 = 1$. De fato, a partir da relação, sendo $f(x) = \sum F_n x^n$, temos que

$$ \frac{f - x}{x^2} = \frac{f}{x} + f$$.

A próxima proposição -- cuja prova é deixada para o leitor -- generaliza o que foi feito acima. 

\

\begin{prop} sejam $f \stackrel{\text{fgo}}{\longleftrightarrow} (a_n)_{n \in \mathbb{N}}$ e $h > 0$. Então

$$ (a_{n+h})_{n \in \mathbb{N}} \stackrel{\text{fgo}}{\longleftrightarrow} \frac{f - a_0 - \ldots - a_{h-1} x^{h-1}}{x^h}.$$
\end{prop}

Agora, vejamos o que acontece quando multiplicamos uma série por uma potência de $n$. De novo, sejam $f \stackrel{\text{fgo}}{\longleftrightarrow} (a_n)_{n \in \mathbb{N}}$. Para achar a função geradora de $(n a_n)_{n \in \mathbb{N}}$ basta observar que $\sum (n a_n) x^n = xf^{'}$. Logo, multiplicar o $n$-ésimo elemento da sequência por $n$ implica em uma 'multiplicação' por $x(d/dx)$ em sua função geradora cuja notação será $xD$. i.e.:

$$ f \stackrel{\text{fgo}}{\longleftrightarrow} (a_n)_{n \in \mathbb{N}} \Rightarrow (xDf) \stackrel{\text{fgo}}{\longleftrightarrow} (n a_n)_{n \in \mathbb{N}}. $$

\

\begin{ex} Considere a seguinte relação de recorrência

$$(n + 1)a_{n+1} = 3a_n +1 \ , \ n \geq  0 , \ a_0 = 1.$$

Se $f$ é a função geradora ordinária  da sequência $(a_n)_{n \in \mathbb{N}}$, então, pelo o que foi feito acima e pela proposição B.16, temos que

$$ f' = 3f + \frac{1}{1 - x},$$

uma equação diferencial de primeira ordem na incógnita $f$.\textcolor{white}{a}\QEDB
\end{ex}

\

Para achar a função que gera $(n^2 a_n)_{n \in \mathbb{N}}$ basta aplicarmos o operador $xD$ novamente obtendo $(xD)^2 f$ como resposta. Em geral

$$ (xD)^k f \stackrel{\text{fgo}}{\longleftrightarrow} (n^k a_n)_{n \in \mathbb{N}}.$$

O que deve gerar $((3 - 7n^2)a_n)_{n \in \mathbb{N}}?$ De novo, basta aplicarmos ao operador $xD$ o mesmo que é aplicado a $n$. i.e., a resposta é $(3 - 7(xD)^2)f$. O caso geral -- também deixado como exercício para o leitor -- é dado por

\

\begin{prop} Sejam  $f \stackrel{\text{fgo}}{\longleftrightarrow} (a_n)_{n \in \mathbb{N}}$ e $P$ um polinômio. Então

$$ P(xD)f \stackrel{\text{fgo}}{\longleftrightarrow} (P(n)a_n)_{n \in \mathbb{N}}.$$

\end{prop}

\

\begin{ex} Vamos usar o que foi discutido para achar uma fórmula fechada para a soma dos quadrados dos primeiros $N$ inteiros positivos. Para tanto, comecemos com a identidade

$$ \sum_{n=0}^N x^n = \frac{x^{N+1} - 1}{x - 1}$$

Se aplicarmos $(xD)^2$ em ambos os lados e fixarmos $x=1$, o lado esquerdo será a soma dos quadrados enquanto o lado direito nos dará a resposta:

$$ \sum_{n=1}^N n^2 = (xD)^2 (\frac{x^{N+1} - 1}{x - 1})\vert_{x=1}.$$

Após duas diferenciações e alguma álgebra, chegamos na resposta 

$$ \sum_{n=1}^N n^2 = \frac{N(N + 1)(2N + 1)}{6} \ , \ \text{para } N \in \{1, 2, \ldots \}.$$

De fato essa é uma fórmula já muito conhecida. O que é importante ressaltar aqui é que as funções geradoras resolvem de forma mecânica muitos problemas complexos envolvendo somas\footnote{Se algo não parece estar certo no exemplo B19, veja o exemplo B27.}.  \textcolor{white}{a}\QEDB

\end{ex}

\

A próxima proposição, cuja prova é facilmente derivada do produto de Cauchy, aborda novamente a multiplicação de funções geradoras.

\

\begin{prop} se $f_k \stackrel{\text{fgo}}{\longleftrightarrow} ((a_k)_n)_{n \in \mathbb{N}}$, para todo $k \in \{1, \ldots, K\}$.
então 

$$ \prod_{k=1}^K f_k \stackrel{\text{fgo}}{\longleftrightarrow} (\sum_{n_1 + n_2 + \ldots + n_K = n} ((a_{1})_{n_1} (a_{2})_{n_2} \ldots (a_{K})_{n_K})_{n \in \mathbb{N}}). $$
\end{prop}

\

\begin{ex} Seja $f(n, k)$ o número de maneiras que um inteiro não negativo $n$ pode ser escrito como uma soma ordenada de $k$ inteiros não negativos. Por exemplo, $f(4,2) = 5$ pois $4=4+0=3+1=2+2=1+3=0+4$. Como $ 1/(1-x) \stackrel{\text{fgo}}{\longleftrightarrow} 1$\footnote{Pois $\sum_{n=0}^\infty x^n = 1/(1-x).$}, nós temos que

$$1/(1-x)^k \stackrel{\text{fgo}}{\longleftrightarrow} (f(n,k))_{n \in \mathbb{N}}.$$

Logo, após algum trabalho algébrico, chegamos a $f(n, k) = \binom{n+k-1}{n}$.\textcolor{white}{a}\QEDB

\end{ex}

\

Agora, sejam $f \stackrel{\text{fgo}}{\longleftrightarrow} (a_n)_{n \in \mathbb{N}}$. Procuraremos a sequência gerada por $f(x)/(1 - x)$. Observe que

$$ \frac{f(x)}{(1 - x)} = (a_0 + a_1x + a_2x^2 + \ldots)(1 + x + x2 + \ldots ) = $$

$$ = a_0 + (a_0 + a_1)x + (a_0 + a_1 + a_2)x^2 + (a_0 + a_1 + a_2 + a_3)x^3 + \ldots $$
 
que indica que
 
\
 
\begin{prop} Se $f \stackrel{\text{fgo}}{\longleftrightarrow} (a_n)_{n \in \mathbb{N}}$, então

$$ \frac{f}{(1 - x)}  \stackrel{\text{fgo}}{\longleftrightarrow} (\sum_{j=0}^n a_j)_{n \in \mathbb{N}} $$

\end{prop}

\

\begin{ex} Os números harmônicos $(H_n)_{n \in \mathbb{N}}$ obedecem

$$ H_n = 1 + \frac{1}{2} + \frac{1}{3} + \ldots + \frac{1}{n} \ , \ \textit{para } n \geq 1.$$

Pela proposição B.22, a função geradora de $(H_n)_{n \in \mathbb{N}}$ é o resultado da multiplicação de $1/(1 - x)$ pela função geradora $f$ da sequência $(1/n)_{n \in \mathbb{N}}$. Como $f = \sum_{n \geq 1} x^n /n$ temos que $f' = 1/(1-x)$. Logo $f = -\log (1 - x)$. Portanto

$$ \sum_{n \geq 1} ^\infty H_n x^n = \frac{1}{1 - x} \log (\frac{1}{1 - x}).$$\textcolor{white}{a}\QEDB

\end{ex}

\begin{ex} Mostremos que os números de Fibonacci satisfazem

$$ F_0 + F_1 + F_2 + \ldots + F_n = F_{n+2} - 1 \ , \ \textit{para } n \geq 0:$$

Pela proposição B.22, se $F$ for a função geradora da sequência de Fibonacci, temos que a sequência geradora da sequência do lado esquerdo da igualdade é dada por $F/(1 - x)$. Como visto no exemplo B2, $F = x/(1 - x - x^2)$. Ainda, pela proposição B9, a função geradora da sequência do lado direito da igualdade é 

$$ \frac{F - x}{x^2} - \frac{1}{1 - x},$$

restando um trabalho algébrico para verificar a igualdade.\textcolor{white}{a}\QEDB

\end{ex}

\

\begin{ob} Essa seção abordou os conceitos básicos de funções geradoras ordinárias. Existem muitos outros tipo de funções geradoras como funções geradoras de Poisson, séries de Lambert, séries de Bell, funções geradoras exponenciais, etc.. Toda sequência com o conjunto de índices começando em 1 possui todos os tipos de funções geradoras a diferença sendo que cada problema é melhor abordado por uma ou mais funções específicas. Esse capítulo não se propõe a definir muitas funções geradoras e, consequentemente, também não abordaremos a questão de como saber qual função geradora se usar em diferentes problemas. O leitor interessado pode consultar as referências dadas no início do capítulo. Vamos, contudo, definir e dar um exemplo de funções geradoras exponenciais. Não enunciaremos porém tudo o que foi provado nessa seção para função geradoras ordinárias possui um análogo para funções geradoras exponenciais. 

\

\begin{df} O símbolo $f \stackrel{\text{fge}}{\longleftrightarrow} (a_n)_{n \in \mathbb{N}}$ significa que $f$ é a função geradora exponencial para a sequência $(a_n)_{n \in \mathbb{N}}$, i.e., que 

$$f = \sum \frac{a_n}{n!} x^n.$$

\end{df} 

\

\begin{ex} Vamos usar funções geradoras exponenciais e a proposição B.16 para calcular
$\sum (n^2 + 4n + 5)/n!$. Claramente, a resposta é o valor em $x=1$ da série 

$$ ((xD)^2 + 4(xD) + 5) e^x = (x^2 + x)e^x + 4xe^x + 5e^x = (x^2 + 5x + 5)e^x.$$

Portanto a resposta é $11e$. Mas nós realizamos um movimento ilegal. Você sabe dizer qual? Nós calculamos a função geradora para um valor de $x$. Tal operação não existe no anel das séries de potência formais. Nesse anel, as séries não assumem valores para diferentes valores de $x$. A letra $x$ é apenas um símbolo representando um pregador no varal. O que pode ser calculado para um valor de $x$ é a série de potência que converge naquele ponto, uma ideia analítica portanto. A ideia é que, se após calcularmos a função geradora usando a teoria formal nós percebermos que a série converge para uma função analítica em um disco no plano complexo, então toda a dedução formal que fizemos é também válida do ponto de vista analítico dentro desse disco. Então podemos mudar o ponto de vista se assim nos convier. A próxima seção se propõe a estudar um pouco a teoria de funções analíticas com o objetivo de provar o teorema B.8 usado no problema de contagem do capítulo 3.  \textcolor{white}{a}\QEDB \\
\textcolor{white}{a}\QEDB
\end{ex}

\end{ob}

\section{A Teoria Analítica Aplicada a Problemas de Contagem}

A teoria formal de séries de potência nos permite manipular relações de recorrência e resolver equações funcionais -- como funções diferenciais -- para séries de potências sem nos preocupar com questões de convergência. Porém, se a série convergir para uma função nós ainda conseguiremos obter informações analíticas possivelmente difíceis de serem obtidas de outra maneira. Esse capítulo se propõe a estudar a teoria analítica necessária para provar o teorema B.8. O leitor interessado em se aprofundar no tema é indicado a consultar \cite{FGA}.
\

\

\subsection{Propriedades Analíticas}

\

Nosso objetivo nessa seção é desenvolver a teoria necessária para provar o teorema B.38. Nós provaremos a existência do raio de convergência para serie de potência em $\mathbb{C}$ depois generalizaremos nosso espaço para $\mathbb{C}^n$ e provaremos em parte o teorema de preparação de Weierstrass, o teorema da função implícita e, finalmente, o teorema B.38. 

\

Considere a série 

$$ f = \sum a_n z^n,$$

onde usamos a letra $z$ para encorajar o entendimento de que estamos agora no plano complexo. Nosso desafio agora é descobrir para quais valores complexos a série $f$ converge e expressá-los em termos dos coeficientes $(a_n)_{n \in \mathbb{N}}$.

Antes, precisaremos da seguinte definição:

\

\begin{df} Seja $L$ um elemento da reta real estendida. Dizemos que $L$ é o limite superior de uma sequência  $(x_n)_{n \in \mathbb{N}}$ se

\

(a) $L$ é finito e

\
 
$\qquad\quad$ (i) para todo $\epsilon > 0$, todos com exceção de um número finito de termos da sequência satisfazem $x_n < L + \epsilon$, e\\

$\qquad\quad$ (ii) para todo $\epsilon > 0$, um número infinito de termos da sequência satisfaz $x_n > L - \epsilon$, ou

\

(b) $L = + \infty$ e, para todo $M > 0$, Existe um $n$ para o qual $x_n > M$,  ou

\

(c) $L = -\infty$ e para todo $x$, existe somente um número finito de termos que obedecem $x_n > x$.

\

Se $L$ é o limite superior de uma sequência $(x_n)_{n \in \mathbb{N}}$ então dizemos que $L = lim sup (x_n).$ 

\end{df}

\

O limite superior tem as seguinte propriedades:

\

$\bullet$ Toda sequência de números reais tem um e apenas um limite superior na reta real estendida.

\

$\bullet$ Se uma sequência tem limite $L$, então $L$ é também o limite superior da sequência.

\

$\bullet$ Se $S$ é o conjunto dos pontos de aderência de $(x_n)_{n \in \mathbb{N}}$, então $lim sup (x_n)$ é o sup do conjunto $S$.

\

Agora estamos prontos para o

\

\begin{teo} Existe um número $0 \leq R \leq +\infty$, chamado de raio de convergência da série $f$ tal que a série converge para todo $z$ com $\vert z \vert < R$ e diverge para todo $z$ com $\vert z \vert > R$. Ainda, $R$ é dado por

$$ R = \frac{1}{lim sup \vert a_n \vert ^{1/n}} \ , \  (1/0 = \infty \ e \ 1/\infty = 0).$$

\end{teo}

\textbf{Prova:} Suponha que $0 < R < \infty$. Seja $z$ tal que $\vert z \vert < R$. Nós mostraremos que a sequência converge em $z$. Seja $\epsilon > 0$ tal que 

$$\vert z \vert < \frac{R}{1 + \epsilon R}. $$

Agora, pela definição de $lim sup$, existe um $N$ para o qual para todo $n > N$ nós temos 

$$ \vert a_n \vert^{1/n} < \frac{1}{R} + \epsilon.$$

Logo, para todo $n > N$,

$$ \vert a_n z \vert \vert z \vert^n < (\vert z \vert ( \frac{1}{R} + \epsilon))^n.$$

Seja $\alpha = \vert z \vert ( \frac{1}{R} + \epsilon)$. Então, pela nossa escolha de $\epsilon$, nós temos $\alpha < 1$.  Portanto a série $\sum a_n z^n$ converge absolutamente pela comparação com os termos de uma série geométrica convergente. 
\

Agora provaremos que a série diverge para $\vert z \vert > R$. De fato, como $\vert z\vert > R$, existe $\epsilon > 0$ para o qual $ \theta = \vert (z/R) - \epsilon z\vert > 0$. Pela definição de $lim sup$, para infinitos valores de $n$ nos temos $\vert a_n\vert^{1/n} > (1/R) - \epsilon$. Logo, para tais valores de $n$,

$$ \vert a_n z^n\vert > \vert (\frac{1}{R} - \epsilon)z	\vert^n = \theta^n$$

que não tem limite superior pois $\theta > 1$. Logo, essa subsequência dos termos da série de potência não converge para zero e portanto a série diverge. Isso completa a prova para o caso $0 < R < \infty$. Os casos onde  $R = 0$ e  $R = +\infty$ são provados similarmente.\textcolor{white}{a}\QEDA

\

Nosso desafio agora é provar dois dos resultados centrais da teoria de análise complexa: o teorema de preparação de Weierstrass e o teorema da função implícita. Para tanto, precisamos desenvolver um pouco nosso linguajar:  

\begin{df} O conjunto de séries de potência formais com $n$ entradas em $\mathbb{C}$ denotado por

$$ \mathbb{C}[[X_1,\ldots, X_n]] = \mathbb{C}[[X]]$$

consiste em todas as expressões da forma

$$ P := \sum_{v\in\mathbb{N}^n} a_v X^v = \sum a_{v_1\ldots v_n} X_1^{v_1}\ldots X_n^{v_n},$$

com $a_v \in \mathbb{C}$ para todo $v \in \mathbb{N}^n$.

\end{df}

As operações de soma, multiplicação e multiplicação por escalar são  definidas de maneira natural, i.e., dado $P_1,P_2 \in \mathbb{C}[[X_1,\ldots,X_n]]$ com $P_1 = \sum_{v\in\mathbb{N}^n} a_v X^v$ e $P_2 = \sum_{v\in\mathbb{N}^n} b_v X^v$ e $\lambda \in \mathbb{C}$ temos

$$ P_1 + P_2 := \sum_{v\in\mathbb{N}^n} (a_v + b_v) X^v, $$

$$ P_1 \cdot P_2 := \sum_m (\sum_{\alpha + \beta = m } a_{\alpha} b_{\beta}) X^m \text{ e} $$

$$ \lambda P_1 := \sum_{v\in\mathbb{N}^n} \lambda a_v X^v.$$

Nós queremos introduzir o conceito de convergência no conjunto $\mathbb{C}[[X]]$. Para tanto, definiremos uma pseudonorma (i.e., uma norma que pode assumir o valor $\infty$). 

\begin{prop} Para cada $r \in \mathbb{R}_{>0}^n$, a função 

$$ \vert\vert\cdot\vert\vert_r : \mathbb{C}[[X]] \rightarrow \mathbb{R} \cup \{ \infty \}, \  \sum_{v\in\mathbb{N}^n} (a_v) X^v \mapsto \sum_{v\in\mathbb{N}^n} (\vert a_v \vert) r^v$$

é uma pseudonorma.	 

\end{prop} 

A prova da proposição B.31 segue do Lema 21.5 de \cite{Kaup}. 

\begin{df} Uma série de potência formal $P \in \mathbb{C}[[X_1,\ldots,X_n]]$ é convergente se existe um $r \in \mathbb{R}_{>0}^n$ para o qual $\vert\vert P\vert\vert < \infty$. O conjunto de séries de potência formais convergentes é denotado por  $\mathbb{C}\{X_1,\ldots, X_n\} = \mathbb{C}\{X\}.$

\end{df} 

Usaremos as notações $\mathbb{C}[[X,Y]] :=  \mathbb{C}[[X_1,\ldots,X_n,Y]]$ e $\mathbb{C}\{X,Y\} :=  \mathbb{C}\{X_1,\ldots,X_n,Y\}$, para $X \in \mathbb{C}^n$, para algum $n \in \mathbb{N}$ e $Y \in \mathbb{C}$. 

\begin{df} Uma série de potência $P \in \mathbb{C}[[X,Y]]$ é dita: 

\

i) distinguida em $Y$ com ordem $b$ se $P(0,\ldots,0,Y) = Y^b\cdot e$ para alguma unidade\footnote{Um elemento é uma unidade se possui um inverso multiplicativo.} $e \in \mathbb{C}[[Y]]$.   

\

ii) um polinômio de Weierstrass de grau $b$ em $Y$ se $P = Y^b + \sum_{j=1}^b a_j Y^{b-j}$ com $a_j \in \mathfrak{m}_{[X]}$, para cada $j \in \{1,\ldots,b\}$, onde $\mathfrak{m}_{[X]} = \{ P \in \mathbb{C}[[X]] : P(0) = 0\}$.

\end{df}

Finalmente possuímos o vocabulário para enunciar o

\begin{teo} \textbf{Teorema de Preparação de Weierstrass.} Se $P \in \mathbb{C}[[X,Y]]$ é distinguida em $Y$ com grau $b$ então existe um polinômio de Weierstrass $\omega \in \mathbb{C}[[X,Y]]$ de grau $b$ e uma unidade $e \in \mathbb{C}[[X,Y]]$ para os quais $P = \omega \cdot e$. Ainda mais, $\omega$ e $e$ são únicos e se $P \in \mathbb{C}\{X,Y\}$ então o mesmo vale para $\omega$ e $e$.  

\end{teo}

\textbf{Prova:} Para $b = 0$ o resultado é trivial, pois a serie de potência $P$ é uma unidade. Seja $b \geq 1$. No caso $n=0$ de novo o resultado segue facilmente: $Y^b$ é o polinômio de Weierstrass de grau $b$ e existe exatamente uma unidade $e \in \mathbb{C}[[Y]]$ tal que $P = eY^b.$ Para $b \geq 1$ e $n \geq 1$ nós definimos a transformação 

$$ \phi: \mathbb{C}[[X,Y]] \rightarrow \mathbb{C}[[S,T]], \ X_j \mapsto S_jT^b, \ Y \mapsto T,$$

onde $S \in \mathbb{C}^n$ e $T \in \mathbb{C}$. Para analisar as propriedades de $\phi$ nós definiremos a seguinte função peso $\gamma$: para cada monômio $aS^\sigma T^\tau$, com $a \neq 0$, definimos 

$$ \gamma(aS^\sigma T^\tau) := \tau - b\vert\sigma\vert \in \mathbb{Z},$$

e escrevemos $\gamma(P) \geq m$, para $P \in \mathbb{C}[[S,T]]$, para indicar que os únicos monômios que aparecem em $P$ tem peso pelo menos $m$. Também usamos os símbolos $\leq, <$ e $>$ de forma correspondente. Em particular, $\gamma(0)$ satisfaz todas as desigualdades. Com isso podemos enunciar o

\begin{lema} i) Im $\phi = \{ P \in \mathbb{C}[[S,T]] : \gamma(P) \geq 0 \}$.

ii) $P$ é distinguida em $Y$ com ordem $b$ se, e somente se, $\phi(P) = eT^b$ para alguma unidade $e \in \mathbb{C}[[S,T]].$

iii) $P$ é um polinômio de Weierstrass de grau $b$ em $Y$ se, e somente se, $\phi(P) = (1 + u)T^b$ para algum $u \in \mathbb{C}[[S,T]]$ tal que $\gamma(u) < 0$ (em particular, $u \in \mathfrak{m}_{[S,T]}$.

iv) Cada unidade $e \in \mathbb{C}[[S,T]]$ tem uma única decomposição da forma $e = (1 + u)\widehat{e}$ tal que $\widehat{e}$ é uma unidade, $\gamma(\widehat{e}) \geq 0$ e $\gamma(u) < 0$. Se $e$ é convergente então $\widehat{e}$ e $u$ também o são.

\end{lema}

Este lema está provado em \cite{Kaup} (lema 22.5).

\

Agora, seja $P$ distinguida em $Y$ com ordem $b$. Logo, pelo Lema B.35 ii), existe exatamente uma única unidade $e^{'} \in \mathbb{C}[[S,T]]$ para a qual $\phi(P) = e^{'}T^b$ que implica na existência da seguinte decomposição de $P$:

$$ \phi(P) = e^{'} T^b \stackrel{\text{iv}}{=} (1+u)\widehat{e}T^b \stackrel{\text{i}}{=} \phi(e)(1+u)T^b,$$

onde $e \in \mathbb{C}[[X,Y]]$ é uma unidade. Logo, pelo lema B.35 iii, a série de potência $\omega := e^{-1}P \in \mathbb{C}[[X,Y]]$ é um polinômio de Weierstrass de grau $b$. Para a unicidade e resultados complementares, indico novamente a referência \cite{Kaup}. \textcolor{white}{a}\QEDA

\

Provaremos agora, como consequência do teorema de preparação de Weierstrass, o teorema da função implícita, outro resultado central da teoria analítica. Cabe comentar, porém, que esses teoremas são de fato equivalente, i.e., o teorema de preparação de Weierstrass pode ser visto como uma consequência do teorema da função implícita. 

\begin{teo} \textbf{Teorema da Função Implícita.} Sejam $X\times Y \subset \mathbb{C}^n \times \mathbb{C}^m$ e $f: X\times Y \rightarrow \mathbb{C}^m$ holomorfa com $f = (f_1,\ldots, f_m)$ tal que, para cada $j \in \{1,\ldots, m\}$, $f_j \in \mathfrak{m}_{[X,Y]}$ e $\frac{\partial f}{\partial y}(0) \neq 0$. Então existe uma vizinhança aberta de $(0,0)$, $U\times W$, tal que existe uma função $g : U \rightarrow W$ com $g = (g_1,\ldots,g_m)$ tal que, para cada $j \in \{1,\ldots, m\}$, $g_j \in \mathfrak{m}_{[U]}$ e $f(u,w) = 0 \iff w = g(u).$

\end{teo}

\textbf{Prova:} Após uma mudança de variáveis, podemos assumir que $A := \frac{\partial f}{\partial y}(0)$ é a matriz $I_m$ (como a condição ``$f(x,A^{-1}y) = 0 \iff y = g(x)$'' implica ``$f(x,y) = 0 \iff y = (A^{-1}g(x)$'' e, pela regra da cadeia, $\frac{\partial(f A^{-1}}{\partial y}(0) = I_m$). Então, em particular, $f_m$ é distinguida com ordem 1 em  $Y_m$. Logo, pelo teorema de preparação de Weierstrass, temos uma única decomposição $f_m = e(y_m - g_0)$ com $g_0 \in \mathfrak{m}_{[X,Y']}$. É fácil de ver que podemos escolher $e = 1$. Usaremos agora um argumento de indução em $m$. O caso $m=1$ é trivial. Suponha agora que o teorema seja válido para $m - 1$
para algum $m > 1$. Mostremos que o mesmo vale para $m$: seja $h = (h_1,\ldots,h_{m-1})$ tal que, para $j \in \{ 1,\ldots, m-1\}$, $h_j \in \mathfrak{m}_{[X,Y']}$ definida por

$$ h(x,y') = f'(x,y',g_0(x,y')).$$

Logo $\frac{\partial h}{\partial y'}(0) = I_{m-1}$, e a hipótese de indução garante a existência de uma única $g' = (g_1,\ldots,g_{m-1})$ tal que, para $j \in \{ 1,\ldots, m-1\}$, $h_j \in \mathfrak{m}_{[X]}$ e

$$ h(x,y') = 0 \iff y' = g'(x).$$

Logo $g:= (g_1,\ldots,g_m)$ tem as propriedades desejadas. Com efeito, $f(x,g(x)) = 0$ pois $f_m(x,g(x)) = g_m(x) - g_0(x,g'(x)) = 0$ e $f'(x,g'(x), g_m(x)) = h(x,g'(x)) = 0$. Por outro lado, se $f(x,y) = 0$ então o fato de que $f_m(x,y) = 0$ implica que $y_m = g_0(x,y')$ e a igualdade $0 = f'(x,y',g_0(x,y')) = h(x,y')$ implica que $y' = g'(x)$. Em particular, $y_m = g_0(x,g'(x)) = g_m(x)$. \textcolor{white}{a}\QEDA

\

Estamos a um lema de estarmos aptos a provar o resultado alvo dessa seção:

\begin{lema} Seja $A(x) = \sum_{n\in \mathbb{N}} a_n x^n$ uma função analítica na região

$$ \Delta = \{ x : \vert x \vert < x_0 + \eta, \vert arg(x-x_0)\vert > \delta\}, $$

onde $x_0, \eta \in \mathbb{R}^+$ e $0 < \delta < \pi/2$. Suponha ainda que existe $\alpha \in \mathbb{R}$ para o qual

$$ A(x) = O((1-x/x_0)^{-\alpha}), \ \text{ para todo x } \in  \Delta.$$

Então 

$$ a_n = O(x_0^{-n}n^{\alpha-1}).$$

\end{lema}

O lema B.37 está provado em \cite{Drmota} (lema 1). Finalmente estamos prontos para o

\begin{teo} Seja $F(x,y)$ uma função analítica ao em uma vizinhança de (0,0) tal que $F(0,y) = 0$ e todos seus coeficientes de Taylor ao redor de zero são reais não negativos. Então existe uma única solução analítica $y = y(x)$ para a equação

$$ y = F(x,y) \eqno(1)$$

cujos coeficientes de Taylor ao redor de zero são todos reais não negativos e $y(0) = 0$. 
\

Se a região de convergência de $F(x,y)$ possuir uma solução positiva $(x_0,y_0)$ do sistema

$$ y = F(x,y) $$

$$ 1 = F_y(x,y),$$

com $F_x(x_0,y_0) \neq 0$ e $F_{yy}(x_0,y_0) \neq 0$ então $y(x)$ é analítica para $\vert x\vert < x_0$ e existem funções $g(x)$ e $h(x)$ analíticas ao redor de $x = x_0$ para as quais $y(x)$ tem uma representação da forma 

$$ y(x) = g(x) - h(x)\sqrt{1 - \frac{x}{x_0}}, \eqno(2)$$

em uma vizinhança de $x_0$. Temos também que 

$$ g(x_0) = y(x_0), \ \text{e} \ h(x_0) = \sqrt{\frac{2x_0F_x(x_0,y_0)}{F_{yy}(x_0,y_0)}}.$$

Ainda, se assumirmos que $[x^n]y(x)$\footnote{$[x^n]y(x)$ é o coeficiente de $x^n$ em $y(x)$.}
$ > 0$ para $n > n_0$ então $x_0$ é a única singularidade de $y(x)$ no círculo $\vert x\vert = x_0$ e obtemos uma expansão assintótica para $[x^n]y(x)$ da forma

$$ [x^n]y(x) = \sqrt{\frac{x_0F_x(x_0,y_0)}{2\pi F_{yy}(x_0,y_0)}} x_0^{-n}n^{-3/2}(1+O(n^{-1}). \eqno(3)$$

\end{teo}

\textbf{Prova} Primeiro mostraremos que existe uma solução (analítica) $y = y(x)$ de $y = F(x,y)$ com $y(0) = 0$. Como $F(0,y) = 0$ temos que a função

$$ y(x) \mapsto F(x,y(x))$$

é uma contração para valores pequenos de $x$. Logo as funções definida iterativamente $y_0(x) :\equiv 0$ e, para $m \geq 0$, $y_{m+1}(x) = F(x,y_m(x))$ convergem uniformemente para uma função $y(x)$ solução de (1). Temos ainda que, para todo $m \geq 0$, $y_m$ é analítica ao redor de zero e tem coeficientes de Taylor reais não negativos. Portanto seu limite uniforme goza das mesmas propriedades.
\

Seja $x_0$ o raio de convergência de $y(x)$. Então $x_0$ é uma singularidade de $y(x)$. Supondo $y(x)$ regular, a função

$$ x \mapsto F_y(x,y(x))$$

é estritamente crescente em $\mathbb{R}^+$. Note que $F_y(0,y(0)) = 0$. Supondo $F_y(x,y(x)) < 1$ temos, pelo teorema da função implícita, que $y(x)$ é regular em uma vizinhança de $x$. Logo existe um limite $x_0$ para o qual $\lim_{x\rightarrow x_0^-} y(x) = y_0$ é finito e satisfaz $F_y(x_0,y_0) = 1$. Se $y(x)$ fosse regular em $x_0$ então

$$y^{'}(x_0) = F_x(x_0,y(x)) + F_y(x_0,y(x_0))y^{'}(x_0)$$

implicaria que $F_x(x_0,y(x_0)) = 0$, o que sabemos não ser verdade. Logo $y(x)$ é singular em $x_0$ (i.e., $x_0$ é o raio de convergência de $y(x)$) e $y(x_0)$ é finito.
\

Considere agora a equação $y - F(x,y) = 0$ ao redor de $(x_0,y_0)$. Nós temos que $1 - F_y(x_0,y_0) = 0$ mas $-F_{yy}(x_0,y_0) \neq 0$. Logo, pelo teorema de preparação de Weierstrass, existem funções $H(x,y)$, $p(x)$ e $q(x)$ analíticas ao redor de $(x_0,y_0)$ que satisfazem $H(x_0,y_0) \neq 0$, $p(x_0) = q(x_0) = 0$ e

$$ y - F(x,y) = H(x,y)((y-y_0)^2 + p(x) (y - y_0) + q(x))$$

ao redor de $(x_0,y_0)$. Como $F_x(x_0,y_0) \neq 0$, nós também temos que $q_x(x_0) \neq 0$. O que significa que qualquer função analítica $y$ que satisfaça $y(x) = F(x,y(x))$ em um subconjunto de uma vizinhança de $x_0$ com $x_0$ em sua fronteira é dada por

$$ y(x) = y_0 - \frac{p(x)}{2} \pm \sqrt{\frac{p(x)^2}{4} - q(x)}.$$

Como $p(x_0) = 0$ e $q_x(x_0) \neq 0$, temos que

$$ \frac{\partial}{\partial x}\left(\frac{p(x)^2}{4} - q(x)\right)_{x = x_0} \neq 0.$$

Portanto existe uma função analítica $K(x)$ tal que $K(x_0) \neq 0$ e 

$$ \frac{p(x)^2}{4} - q(x) = K(x)(x-x_0)$$

ao redor de $x_0$. Chegamos então em uma representação de $y = y(x)$ da forma

$$ y(x) = g(x) - h(x)\sqrt{1 - \frac{x}{x_0}},\eqno(4)$$

onde $g(x)$ e $h(x)$ são analíticas ao redor de $x_0$ e satisfazem $g(x_0) = y_0$ e $h(x_0) < 0$. Para calcular $h(x_0)$, usaremos o teorema de Taylor:

$$ 0 = F(x,y(x)) = F_x(x_0,y_0)(x-x_0) + \frac{1}{2}F_{yy}(x_0,y_0)(y(x) - y_0)^2 + \ldots = $$

$$ = F_x(x_0,y_0)(x-x_0) + \frac{1}{2}F_{yy}(x_0,y_0)h(x_0)^2(1-x/x_0) + O(\vert x - x_0 \vert^{3/2}).$$

Comparando os coeficientes de $(x-x_0)$ nós obtemos

$$ h(x_0) - \sqrt{\frac{2x_0F_x(x_0,y_0)}{F_{yy}(x_0,y_0)}}.$$

Com o objetivo de aplicar o lema B.37, precisamos mostrar que $y(x)$ pode ser estendida analiticamente para uma região da forma $\Delta$. A representação (4) nos fornece tal extensão em uma vizinhança de $x_0$. Agora, suponha que $\vert x_1 \vert = x_0$ e $\vert arg(x_1)\vert \geq \delta$. Então, supondo que existe $n_0 \in \mathbb{N}$ para o qual para $y_n > 0$ para $n > n_0$, temos que $\vert y(x_1) < y(\vert x_1\vert = y(x_0)$ e portanto

$$ \vert F_y(x_1,y(x_1))\vert \leq F_y(\vert x_1 \vert , \vert y(x_1)\vert) < F_y(\vert x_1\vert, y(\vert x_1 \vert)) = F_y(x_0, y_0) = 1.$$

Logo, $F_y(x_1,y(x_1)) \neq 1$ e o teorema da função implícita mostra que existe uma solução analítica $\widehat{y}(x)$ em uma vizinhança de $x_1$. Para $\vert x \vert < x_0$, essa solução se iguala à serie de potência  $y(x)$ e, para $\vert x\vert \geq x_0$ ela fornece uma extensão analítica para a região $\Delta$. Finalmente, podemos aplicar o lema B.37 para provar (3).    \textcolor{white}{a}\QEDA

\

\chapter{Entropia de Shannon}

\

Esse apêndice se propõe a estabelecer os conceitos básicos da entropia de Shannon com ênfase aos conceitos necessários para compreender o método de entropy compression. Ele está baseado, em sua maior parte, nos primeiros capítulos de \cite{Ash} e \cite{Thomas}. 

\

Entropia é uma medida da imprevisibilidade de um pedaço de informação. Para se ter uma ideia intuitiva desse conceito, considere um lançamento de moeda. Quando a moeda é honesta, ou seja, quando a probabilidade de cara e coroa é a mesma, então a entropia desse evento (lançamento de uma moeda) assume seu valor máximo. Isso porque não temos nenhum indício de qual será o resultado do evento. Por outro lado, o lançamento de uma moeda que possui duas caras e nenhuma coroa tem entropia zero. Isso porque sabemos, de antemão, exatamente qual será o resultado do evento. 
\

O desafio então é criar uma medida para a imprevisibilidade de eventos. Isso será feito exigindo que essa medida obedeça axiomas apropriados. Estabelecidos os axiomas, provaremos que apenas uma família de funções consegue satisfazer a todos eles.

\

\section{Os Axiomas da Medida de Imprevisibilidade} 

\

Suponha que um experimento envolva a observação de uma variável aleatória discreta $X$ com imagem $\lbrace x_1,\ldots,x_M\rbrace$ e probabilidades $p_1,\ldots\, p_M$, respectivamente. Suponha ainda que $p_i > 0$, para todo $i$. Na tentativa de associar uma incerteza (imprevisibilidade) à variável $X$, nós criaremos as funções $h$ e $H$. A função $h$ denotará a incerteza associada a um evento de probabilidade $p$ (logo será definida no intervalo (0,1]). Ou seja, se o evento $\lbrace X=x_i\rbrace$ tem probabilidade $p_i$, então a incerteza removida (ou informação fornecida) ao revelar que $X$ amostrou o valor $x_i$ é $h(p_i)$. Para cada $M$, definiremos $H_M$ das $M$ variáveis $p_1,\ldots, p_M$ como $H_M(p_1,\ldots,p_M) = \sum_{i=1}^M p_i h(p_i)$. Logo $H_M(p_1,\ldots,p_M)$ é o valor esperado de incerteza removida ao se revelar o valor que $X$ assume. Em prol de uma notação menos carregada denotaremos $H_M(p_1,\ldots,p_M)$ por $H(p_1,\ldots,p_M)$ ou $H(X)$. 
\

A seguir, estabeleceremos axiomas que traduzirão o que esperamos de uma medida de incerteza. Seja $f(M)$ o valor esperado de incerteza associado a $M$ eventos equiprováveis. Ou seja, $f(M) = H(1/M\ldots, 1/M)$. Claramente, a incerteza é maior ao se escolher uma pessoa ao acaso em uma multidão do que ao se lançar uma moeda. Logo, nossa primeira exigência sobre a função $H(X)$ é

\

$\bullet \quad f(M) = H(1/M,\ldots,1/M)$ é uma função estritamente crescente de $M$.  

\

Agora, sejam $X$ e $Y$ duas variáveis aleatórias uniformes independentes com imagens $\lbrace x_1,\ldots,x_M\rbrace$ e $\lbrace y_1,\ldots,y_L\rbrace$ respectivamente. Logo o experimento envolvendo $X$ e $Y$ tem $M\cdot L$ resultados igualmente prováveis e a incerteza esperada desse experimento deve ser $f(M\cdot L)$. Se o valor de $X$ é revelado, o valor esperado de incerteza de $Y$ não deve mudar pois supusemos independência entre $X$ e $Y$. Logo, o valor esperado de incerteza associado ao experimento envolvendo $X$ e $Y$ menos o valor esperado de incerteza removida ao revelar-se o valor que $X$ assume deve ser igual ao valor esperado de incerteza de $Y$. Portanto a segunda exigência que fazemos é

\

$\bullet \quad f(ML) = f(M) + f(L)$, para $M,L \in \mathbb{N}$.

\

Considere novamente a variável aleatória  $X$ (não necessariamente uniforme) com imagem $\lbrace x_1,\ldots,x_M\rbrace$. Seja $A = \lbrace x_1,\ldots,x_r\rbrace$ e \newline $B = \lbrace x_{r+1},\ldots, x_M\rbrace$. Construiremos um experimento da seguinte maneira: primeiro selecionamos $A$ ou $B$, escolhendo $A$ com probabilidade $p_1 + \ldots + p_r$ e $B$ com probabilidade $p_{r+1} + \ldots + p_M$. Se o grupo $A$ for escolhido então, para $0 \leq i \leq r$, selecionamos $x_i$ com probabilidade $p_i/(p_1 + \ldots + p_r)$ (a probabilidade condicional de $x_i$ dado que o valor de $X$ está em $A$). Da mesma forma, se o grupo $B$ for escolhido então, para $r + 1 \leq i \leq M$, $x_i$ é escolhido com probabilidade $p_i/(p_{r+1} + \ldots + p_M)$. A próxima figura representa esse experimento:

\

\begin{center}

\begin{tikzpicture}

\node (1) at ( 0, 0) {}      ;
\node (2) at ( 2.5, 2.5) {A}     ;
\node (3) at ( 6, 6) {$x_1$} ;
\node (4) at ( 6, 4.5) {$x_2$} ;
\node (5) at ( 6, 1) {$x_r$} ;

\node (6) at ( 2.5, -2.5) {B} ;
\node (7) at ( 6, -1) {$x_{r+1}$} ;
\node (8) at ( 6, -2) {$x_{r+2}$} ;
\node (9) at ( 6, -6) {$x_M$}     ;

\draw [ line width=1.5] (0,0) -- (0.2,0.2)   node[ near start ,above = 30pt] {$p_1 + \ldots + p_r$};

\draw [ line width=1.5]  (0,0) -- (0.2,-0.2) node[ near start , below = 35pt] {$p_{r+1} + \ldots + p_M$};

\draw [ line width=1.5]  (3,3) -- (4,4) node[ near end , above = 28pt] {  {{\Large $\frac{p_1}{\sum_{i=1}^r p_i}$}}};

\draw [ line width=1.5]  (3.2,2.9) -- (5.3,4.1) node[ near end , below = 8pt] {  {{\Large $\frac{p_2}{\sum_{i=1}^r p_i} $}}};

\draw [ line width=1.5]  (3.2,2.9) -- (5.3,4.1) node[ near end , below = 5pt] {  {{\Large $\qquad\qquad\cdot$}}};

\draw [ line width=1.5]  (3.2,2.9) -- (5.3,4.1) node[ near end , below = 19pt] {  {{\Large $\qquad\qquad\cdot$}}};

\draw [ line width=1.5]  (3.2,2.9) -- (5.3,4.1) node[ near end , below = 33pt] {  {{\Large $\qquad\qquad\cdot$}}};

\draw [ line width=1.5]  (3.2,2.2) -- (5.3,1.3) node[ near start , below = 8pt] {  {{\Large $\frac{p_r}{\sum_{i=1}^r p_i}$}}};

\draw [ line width=1.5]  (3.2,-2.2) -- (5.3,-1.3) node[ near start , above = 15pt] {  {{\Large $\frac{p_{r+1}}{\sum_{i=r+1}^M p_i}$}}}; 

\draw [ line width=1.5]  (3.2,-2.4) -- (5.3,-2.1) node[ near end , below = 10pt] {  {{\Large $\frac{p_{r+2}}{\sum_{i=r+1}^M p_i}$}}}; 

\draw [ line width=1.5]  (3.2,-2.4) -- (5.3,-2.1) node[ near end , below = 20pt] { {\Large $\qquad\qquad\cdot$}}; 

\draw [ line width=1.5]  (3.2,-2.4) -- (5.3,-2.1) node[ near end , below = 34pt] { {\Large $\qquad\qquad\cdot$}};

\draw [ line width=1.5]  (3.2,-2.4) -- (5.3,-2.1) node[ near end , below = 48pt] { {\Large $\qquad\qquad\cdot$}};

\draw [ line width=1.5]  (3,-3) -- (4,-4) node[ near start , below = 25pt] {  {{\Large $\frac{p_M}{\sum_{i=r+1}^M p_i}$}}};

\draw [line width=1.5] (1) -- (2) ;
\draw [ line width=1.5] (1) -- (6); 
\draw [line width=1.5] (2) -- (3);	
\draw [line width=1.5] (2) -- (4);
\draw [line width=1.5] (2) -- (5);
\draw [line width=1.5] (6) -- (7);
\draw [line width=1.5] (6) -- (8);
\draw [line width=1.5] (6) -- (9);

\end{tikzpicture}

\end{center}

\begin{picture}(2,2)

\put(0,1){\begin{footnotesize}
$\qquad\qquad\qquad\qquad\quad$ Figura 1: Representação do experimento.  

\end{footnotesize}
}
\end{picture}

\

O evento descrito é equivalente a $X$. Com efeito, seja $Y$ a variável aleatória associada ao experimento. Seja $x_i \in A$. Logo 

$$ \textbf{P}(Y=x_i) = \textbf{P}(``A \quad é \quad escolhido \quad e \quad x_i \quad é \quad selecionado") = $$

$$ = \textbf{P}(``A \quad é \quad escolhido")\textbf{P}(``x_i \quad é \quad selecionado"\vert``A \quad é \quad escolhido") = $$

$$ = (\sum_{i=1}^r p_i) \frac{p_i}{\sum_{i=1}^r p_i} = p_i .$$

\

O mesmo raciocínio se aplica se $x_i \in B$. Logo $X$ e $Y$ têm a mesma distribuição. A incerteza associada a $X$ é $H(p_1,\ldots,p_M)$. Se for revelado qual dos grupos será escolhido no experimento, o valor esperado de incerteza que será removido é $H(p_1 + \ldots + p_r, p_{r+1} + \ldots p_M)$. Ao se escolher o grupo $A$, a incerteza restante é 

$$ H(\frac{p_1}{\sum_{i=1}^r p_i}, \ldots, \frac{p_r}{\sum_{i=1}^r p_i}).$$ 

Ao se escolher o grupo B, a incerteza
restante é

$$ H(\frac{p_{r+1}}{\sum_{i=r + 1}^M p_i}, \ldots, \frac{p_M}{\sum_{i=r+1}^M p_i}).$$

Logo o valor esperado de incerteza removida depois que o grupo for escolhido é 

$$ (p_1 + \ldots + p_r) H(\frac{p_1}{\sum_{i=1}^r p_i}, \ldots, \frac{p_r}{\sum_{i=1}^r p_i}) + $$

$$+ (p_{r+1} + \ldots + p_M) H(\frac{p_{r+1}}{\sum_{i=r + 1}^M p_i}, \ldots, \frac{p_M}{\sum_{i=r+1}^M p_i}).$$

\

É desejável que o valor esperado de incerteza relacionado ao experimento menos o valor esperado de incerteza removida ao se revelar o grupo escolhido seja igual o valor esperado de incerteza restante após o grupo ser selecionado. Logo o terceiro axioma é

$$\bullet \quad H(p_1,\ldots, p_M) = H(p_1+\ldots+p_r, p_{r+1} + \ldots + p_M)  + $$

$$+ (p_1 + \ldots + p_r) H(\frac{p_1}{\sum_{i=1}^r p_i}, \ldots, \frac{p_r}{\sum_{i=1}^r p_i}) + $$

$$+ (p_{r+1} + \ldots + p_M) H(\frac{p_{r+1}}{\sum_{i=r + 1}^M p_i}, \ldots, \frac{p_M}{\sum_{i=r+1}^M p_i}).$$ 

\

Finalmente, nós exigiremos, por conveniência matemática, que $H(p,1-p)$ seja uma função contínua em $p$. Intuitivamente isso significa uma pequena variação nas probabilidades dos valores de $X$ representam uma pequena variação da incerteza de $X$.

\

\begin{teo} A única função que satisfaz os quatro axiomas é 

$$H(p_1,\ldots,p_M)= - C \sum_{i=1}^M p_i log(p_i), $$

onde $C$ é um número positivo arbitrário e a base do logaritmo é qualquer número maior que 1.
\end{teo}

\textbf{Prova:} Não é difícil verificar que a função acima satisfaz os quatro axiomas. Provaremos que se uma função satisfaz os quatro axiomas citados então essa função pertence à família de funções citada acima. Nessa prova, os logaritmos estão em uma mesma base fixa. Isso é feito sem perda de generalidade pois mudar a base do logaritmo é equivalente a mudar a constante $C$ ($\log_a x = log_a b \cdot log_b x$). Por conveniência, dividiremos a prova em algumas partes:

\

a) É fácil provar, usando o axioma 2 e indução em $k$, que $f(M^k) = kf(M)$, para $M,k \in \mathbb{N}^*$.

\

b) $f(M) = C\log M$, para $M \in \mathbb{N}^*$ e $C \in \mathbb{R}^+$. Com efeito, seja $M = 1$. Pelo segundo axioma , $f(1) = f(1\cdot 1) = f(1) + f(1)$. Logo
$f(1) = 0$. O que está de acordo com a noção intuitiva de que um evento com probabilidade 1 não tem incerteza associada. Agora, sejam $M$ um inteiro positivo maior que 1 e $r > 0$ um número natural fixo. Existe algum número $k$ para o qual $M^k \leq 2^r < M^{k+1}$. Logo, pelo primeiro axioma, $f(M^k) \leq f(2^r) < f(M^{k+1})$ e, por a), $kf(M) \leq rf(2) < (k+1)f(M)$ ou $k/r \leq f(2)/f(M) < (k+1)/r$. O logaritmo é uma função estritamente crescente (pois sua base é maior que 1), logo $\log M^k\leq \log 2^r < \log M^{k+1}$ ou $ k/r \leq \log 2/\log M < (k+1)/r$. Como ambos $f(2)/f(M)$ e $\log 2/\log M$ estão entre $k/r$ e $(k+1)/r$, temos que

$$ \vert \frac{\log 2}{\log M} - \frac{f(2)}{f(M)} \vert < \frac{1}{r}.$$

Como $r$ foi escolhido de maneira arbitrária temos que 

$$ \frac{\log 2}{\log M} = \frac{f(2)}{f(M)} $$

ou $f(M) = C\log M$, onde $C = f(2)/\log 2$. Note que $C$ é positivo já que $f(1) = 0$ e $f$ é estritamente crescente. 

\

c) $H(p,1-p) = -C(p\log p + (1-p)\log(1-p))$, se $p \in \mathbb{Q}$. Com efeito, seja $p = r/s$, com $r,s \in \mathbb{N}^*$. Pelo terceiro axioma temos que

$$ f(s) = H(\frac{r}{s},\frac{s-r}{s}) + \frac{r}{s}f(r) + \frac{s-r}{s}f(s-r)$$

Por b, temos

$$ C\log s = H(p,1-p) + Cp\log r + C(1-p)\log(s-r).$$

Logo

$$ H(p,1-p) = -C(p\log r - \log s + (1-p)\log(s-r)) = $$

$$= -C(p\log r -p\log s + p\log s - \log s + (1-p)\log(s-r)) = $$

$$ = -C(p\log \frac{r}{s} + (1-p)\log \frac{s-r}{s}) = -C(p\log p + (1-p)\log (1-p)).$$

\

d) $H(p,1-p) = -C(p\log p + (1-p)\log(1-p))$, para todo $p$. Esse resultado segue de c) e do axioma 4. Com efeito, seja $p \in [0,1]$. Por continuidade

$$H(p,1-p) = \lim_{p^{'} \rightarrow p} H(p^{'},1-p^{'}).$$

Em particular, podemos nos aproximar de $p$ por uma sequência de números racionais. Logo

$$\lim_{p^{'} \rightarrow p} H(p^{'},1-p^{'}) = \lim_{p^{'} \rightarrow p} -C(p^{'}\log p^{'} + (1-p^{'})\log(1-p^{'})) = $$

$$ =  - C(p\log p + (1-p)\log (1-p)).$$

\

e) $H(p_1,\ldots,p_M) = -C\sum_{i=1}^M p_i\log p_i$, para $M \in \mathbb{N}^*$. Com efeito, por indução em $M$: O resultado já foi provado para $M = 1$ e $M = 2$ (lembre que $\sum_{i=1}^Mp_i = 1$). Seja $M>2$. Pelo terceiro axioma temos

$$ H(p_1,\ldots,p_M) = H(p_1 + \ldots + p_{M-1},p_M) + $$

$$ + (p_1,\ldots,p_{M-1})H(\frac{p_1}{\sum_{i=1}^{M-1}p_i},\ldots,\frac{p_{M-1}}{\sum_{i=1}^{M-1}p_i}) + p_MH(1). $$

Assumindo que a formula vale para inteiros menores que $M$ temos:

$$  H(p_1,\ldots,p_M) = -C((p_1 + \ldots + p_{M-1})\log (p_1 + \ldots + p_{M-1}) + p_M\log p_M) - $$

$$ - C(p_1 + \ldots + p_{M-1})(\frac{p_1}{\sum_{i=1}^{M-1}p_i}\log\frac{p_1}{\sum_{i=1}^{M-1}p_i} + \ldots + \frac{p_{M - 1}}{\sum_{i=1}^{M-1}p_i}\log\frac{p_{M - 1}}{\sum_{i=1}^{M-1}p_i}) = $$
\begin{samepage}
$$= -C((\sum_{i=1}^{M-1}p_i\log (\sum_{i=1}^{M-1}p_i) + p_M\log p_M) - C(\sum_{i=1}^{M-1}p_i\log p_i - $$

$$ - (\sum_{i=1}^{M-1}p_i)\log \sum_{i=1}^{M-1}p_i) =  -C\sum_{i=1}^M p_i \log p_i. $$\textcolor{white}{a}\QEDA

\

\end{samepage} 

Portanto $H(p_1,\ldots,p_M)= - C \sum_{i=1}^M p_i log(p_i)$ para algum $C$. Vimos que diferentes valores de $C$ equivalem a diferentes bases do logaritmo. Assim, diferentes unidades de entropia são definidas para cada base. Shannon ou bit para base 2, nat para base $e$ e hartley para base 10 são alguns exemplos.
\

Voltemos ao exemplo do lançamento de uma moeda. Seja $H(X) = -\sum_{i=1}^2 p_i\log_2 p_i$. Considere uma moeda com valores conhecidos para as probabilidades de cara e coroa. A situação de máxima entropia, nesse evento, é quando essas probabilidades são iguais (esse resultado vale para o caso geral, conforme provado no teorema C.11). Nesse caso temos que $H(X) = 1$. Se as probabilidades não forem iguais, então temos mais condições de prever o resultado logo a entropia deve ser menor (até o caso extremo onde uma das probabilidades é igual a 1). A próxima figura mostra a evolução da entropia e da função $h$ com relação a probabilidade de sair cara. Também constam as parcelas $p(-\log_2(p))$ e $(1 - p)(-\log_2(1 - p))$ que compõe a entropia.

\begin{center}

\includegraphics[scale=0.7]{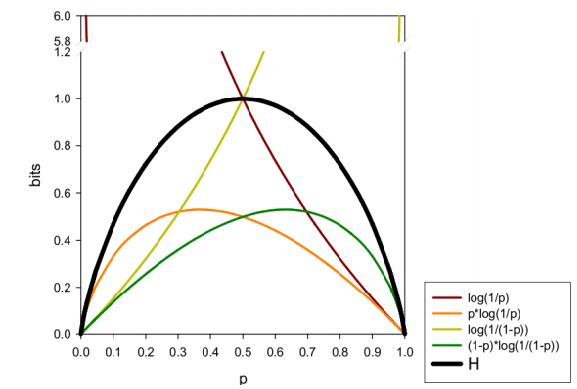}

\end{center}

\begin{picture}(10,10)

\put(1,1){$\qquad\qquad\qquad\qquad\quad$\begin{footnotesize} 
 Figura 2: Evolução da entropia.

\end{footnotesize}
}
 
\end{picture}
 
\

\begin{ob} A convenção $0 log0 = 0$ é usual na teoria da informação. Essa convenção
é facilmente justificada por continuidade, uma vez que $  xlog(x)$ tende a zero
quando $x$ tende a zero (ver Figura 2). Assim, a adição de termos com probabilidade
zero não muda a entropia.
\end{ob}

\

\begin{ob} Os axiomas para a medida de incerteza aqui apresentados são os axiomas de Shannon \cite{Shannon}. Fadiev, em \cite{Fadiev}, mostrou uma versão mais fraca desses axiomas determinando a mesma função de entropia.
\end{ob}

\section{Entropia Relativa e Propriedades Básicas}

\

Tendo em mente que o objetivo principal desse apêndice é dar ao leitor a base necessária para entender o método de entropy compression, vamos definir entropia relativa sem definir outros conceitos mais básico e/ou relacionados como entropia conjunta, entropia condicional e informação mútua. O leitor interessado pode consultar \cite{Ash} e \cite{Thomas}. 

\
 
\begin{df} A entropia relativa (ou distância de Kullback–Leibler) entre duas distribuições de probabilidade $p(x)$ e $q(x)$ é dada por

$$ D(p \vert\vert q) = \sum_{x \in \mathcal{X}} p(x)\log \frac{p(x)}{q(x)} = E_p[\log\frac{p(X)}{q(X)}]. $$
\end{df}

Na definição acima, usamos as convenções  $0 \log \frac{0} {0} = 0$, $0 \log \frac{0}{q} = 0$ e $p \log \frac{p}{0} = \infty$. Logo, se existir $x \in \mathcal{X}$ para o qual $p(x) > 0$ e $q(x)=0$, então $D(p\vert\vert q) = \infty$ . A seguir, mostraremos que a entropia relativa é não-negativa e é zero se, e somente se, $p = q$. Por essas propriedades , apesar de não ser simétrica nem respeitar a desigualdade triangular, a entropia relativa é geralmente entendida como a distância entre duas distribuições de probabilidade. 

\

\begin{ex} Sejam $\mathcal{X} = \lbrace 0,1\rbrace$ e $p$ e $q$ distribuições de probabilidade em $\mathcal{X}$ com
 
 $$ p(0) = 1-r, \quad p(1) = r, \quad q(0) =1-s\quad \textit{e} \quad q(1) = s.$$

\
 
 Então

$$ D(p\vert\vert q) = (1-r)\log \frac{1-r}{1-s} + r \log \frac{r}{s} $$
 
 e
 
$$ D(q\vert\vert p)= (1-s)\log \frac{1-s}{1-r} + s \log \frac{s}{r}. $$

\
 
 Se $r = s$, então $D(p\vert\vert q)= D(q\vert\vert p)=0$. Se $r = \frac{1}{2}$, $s = \frac{1}{4}$, então

\

 $$ D(p\vert\vert q)= \frac{1}{2} \log \frac{1/2}{3/4} + \frac{1}{2} \log \frac{1/2}{1/4} = 1 - \frac{1}{2} \log3 = 0.2075 \textit{ bit}.$$

\
 
Enquanto

\
 
$$ D(q\vert\vert p)= \frac{3}{4} \log \frac{3/4}{1/2} + \frac{1}{4} \log \frac{1/4}{1/2} = \frac{3}{4} \log3 - 1 = 0.1887 \textit{ bit}.$$\textcolor{white}{a}\QEDB
\end{ex}

Agora provaremos algumas propriedades básicas de entropia de entropia relativa. De especial importância para nós, provaremos que a entropia relativa é não-negativa. Esse fato será usado mais adiante para provar o teorema C.26 -- o resultado de maior importância para nós nessa seção. 
\

Primeiro, observamos que a entropia é sempre positiva, i.e., $H(X) \geq 0$. Isso segue diretamente do fato que $-p_i\log p_i \geq 0$, para todo $i$. 

Para provar as próximas propriedades, precisaremos das seguintes ferramentas:

\

\begin{df} Uma função $f(x)$ é convexa num intervalo $(a,b)$ se, para todo $x_1,x_2 \in (a,b)$ e $ 0 \leq \lambda \leq 1$, 

$$f(\lambda x_1 + (1 - \lambda)x_2) \leq \lambda f(x_1)+(1 - \lambda)f(x_2).$$  

Uma função $f$ é estritamente convexa se acima temos igualdade somente se $\lambda = 0$ ou $\lambda = 1$. 
\end{df}

\

\begin{df} Uma função é (estritamente) côncava se $-f$ é (estritamente) convexa. 
\end{df}
\

Intuitivamente, uma função é convexa/côncava se sempre está abaixo/acima de qualquer secante. $x^2$, $\vert x\vert$, $e^x$, $x \log x$ são exemplos de funções convexas (para $x \geq 0)$. $\log x$ e $\sqrt{x}$ são exemplos de funções côncavas (para $x \geq 0)$. Funções lineares são côncavas e convexas. 

\

\begin{teo} Seja $f$ uma função com segunda derivada não-negativa/positiva num intervalo $(a,b)$. Então $f$ é convexa/estritamente convexa nesse intervalo.
\end{teo}

\textbf{Prova:} Seja $x_0 \in (a,b)$. Considere a expansão de Taylor de $f$ ao redor de $x_0$:

$$ f(x)= f(x_0) + f^{'}(x_0)(x - x_0) + \frac{f^{''}(x^*)}{2}(x - x_0)^2,$$

Onde $x^* \in B_{\vert x - x_0\vert}(x_0)$. Por Hipótese, $f^{''}(x^*) \geq 0$. Logo, o último termo é não negativo para todo $x \in (a,b)$.  
\

Sejam $x_1,x_2 \in (a,b)$ e $\lambda \leq 1$ tais que $x_0 = \lambda x_1 +(1-\lambda)x_2$. Logo

$$ f(x_1) \geq f(x_0)+ f^{'}(x_0)(1-\lambda)(x_1 - x_2).$$

e

$$ f(x_2) \geq f(x_0) + f^{'}(x_0)\lambda (x_2 - x_1).$$ 

Multiplicando a primeira desigualdade por $\lambda$, a segunda por $1-\lambda$ e somando os resultados, obtemos

$$f(\lambda x_1 + (1-\lambda)x_2) \leq \lambda f(x_1)+(1-\lambda)f(x_2).$$

A prova para convexidade estrita segue a mesma ideia.\textcolor{white}{a}\QEDA

\

Esse teorema nos permite verificar facilmente a convexidade estrita de $x^2$, $e^x$ e $x\log x$ (para $x \geq 0$) e a concavidade estrita de $\log x$ e $\sqrt x$ (para $x \geq 0$). 
\

O próximo resultado é de suma importância para derivar as propriedades básicas da teoria da informação. Abaixo, $E[X]$ denota o valor esperado da v.a. $X$.

\begin{teo}[Desigualdade de Jensen - Caso Discreto] Seja $f$ uma função convexa e $X$ uma variável aleatória discreta. Então 

$$ E[f(X)] \geq f(E[X]).$$

Além disso, se $f$ é estritamente convexa, a igualdade na expressão acima implica que $X = E[X]$ com probabilidade 1 (i.e., $X$ é constante).
\end{teo}

\textbf{Prova:} A prova será por indução na cardinalidade de $\mathcal{X}$. Para $\vert \mathcal{X} \vert = 2$, a desigualdade toma a forma

$$ p_1 f(x_1)+ p_2 f(x_2) \geq f(p_1 x_1 + p_2 x_2)$$

que segue diretamente da definição de função convexa. Supondo que a desigualdade vale para $\vert \mathcal{X} \vert = k - 1$ concluiremos que também vale para $\vert \mathcal{X} \vert = k$. Escrevendo $p_i =  p_i/(1 - p_k)$ para $i \in \lbrace 1,2,\ldots,k-1\rbrace$ temos:

$$ E[f(X)] = \sum_{i=1}^k p_i f(x_i) = p_k f(x_k) + (1 - p_k) \sum_{i=1}^{k-1} p_i^{'} f(x_i) \geq$$

$$\geq p_k f(x_k) + (1 - p_k) f(\sum_{i=1}^{k-1} p_i^{'} x_i) \geq f(p_k x_k + (1 - p_k) \sum_{i=1}^{k-1} p_i^{'}x_i) = f(\sum_{i=1}^{k} p_i x_i),$$

onde a primeira desigualdade segue da hipótese de indução e a segunda da definição de convexidade. \textcolor{white}{a}\QEDA

\

O teorema acima pode ser estendido para distribuições contínuas por argumentos de continuidade. De posse dessa ferramenta, enfim somos capazes de demonstrar o que buscávamos:

\begin{teo} Sejam $p(x)$, $q(x)$, $x \in \mathcal{X}$, duas distribuições de probabilidade. Então

$$ D(p\vert\vert q) \geq 0,$$

com igualdade se, e somente se, $p(x)= q(x)$ para todo $x \in \mathcal{X}$. 
\end{teo}

\textbf{Prova:} Seja $A= \lbrace x : p(x) > 0\rbrace$ o suporte de $p(x)$. Então

$$ -D(p\vert\vert q) = -\sum_{x \in A} p(x)\log \frac{p(x)}{q(x)} = \sum_{x \in A} p(x)\log \frac{q(x)}{p(x)} \leq $$

$$ \leq \log \sum_{x \in A} p(x) \frac{q(x)}{p(x)} = \log \sum_{x \in A} q(x) \leq \log \sum_{x \in X} q(x) = \log 1 =0,$$

onde a primeira desigualdade segue da desigualdade de Jensen. Como $\log t$ é uma função estritamente côncava, temos igualdade na primeira desigualdade se, e somente se, $q(x)/p(x)$ é constante (i.e., $q(x)= c p(x)$ para todo $x \in \mathcal{X}$). Então
$ \sum_{x \in A} q(x)= c \sum_{x \in A} p(x)$. Temos igualdade na segunda desigualdade se, e somente se, $\sum_{x \in A} q(x)= \sum_{x \in \mathcal{X}} q(x) = 1$, o que implica que $c = 1$. Logo,  $D(p\vert\vert q)=0$ se, e somente se, $p(x)= q(x)$ para todo $x \in \mathcal{X}$.\textcolor{white}{a}\QEDA

\

A seguir, provaremos a intuitiva propriedade de entropia que diz que, fixado $\mathcal{X}$, a distribuição de probabilidade que resulta na maior entropia é a uniforme. O conceito de entropia relativa nos fornece uma prova simples de tal fato. Para uma prova alternativa veja o teorema 1.4.2 de \cite{Ash}.   

\

\begin{teo} $H(X) \leq log \vert \mathcal{X} \vert$, com igualdade se, e somente se, $X$ tem distribuição uniforme.
\end{teo}

\textbf{Prova:} Sejam $u(x) = \frac{1}{\vert \mathcal{X} \vert}$ a distribuição uniforme sobre $\mathcal{X}$ e $p(x)$ a distribuição de probabilidade de $X$. Então

$$ D(p \vert \vert u) = \sum_{x \in \mathcal{X}} p(x) \log \frac{p(x)}{u(x)} = \log|X| - H(X).$$

Logo, pela não-negatividade da entropia relativa, temos que

$$ 0 \leq D(p \vert \vert u) = \log \vert X \vert - H(X) \therefore H(X) \leq \log \vert X \vert.$$ \textcolor{white}{a}\QEDA

\

\section{Compressão de Dados}

\

Como o nome indica, comprimir dados é o ato de representar um pedaço de informação usando menos espaço (bits). Isso pode ser feito com ou sem perdas de informação. Aqui estudaremos a compressão sem perdas.
\

Esse objetivo pode ser alcançado designando símbolos curtos para frações mais frequentes do pedaço de informação (e, consequentemente, símbolos longos para frações menos frequentes). Por exemplo, no código Morse, o símbolo mais frequente é representado por um único ponto. Nosso objetivo nessa seção é descobrir até quanto e como conseguimos comprimir um pedaço de informação. 

\

\subsection{Códigos}

\

\begin{df} Um código $C$ para uma variável aleatória $X$ é uma função que leva $\mathcal{X}$, a imagem de $X$, em $D^*$, o conjunto de todas sequências finitas de símbolos de um alfabeto $D$. Para cada $x \in \mathcal{X}$, denota-se por $l(x)$  o comprimento (em número de símbolos de $D$) da palavra $C(x)$.
\end{df}

\

\begin{df} O comprimento esperado $L(C)$ de um código $C$ para uma variável aleatória $X$ com distribuição de probabilidade $p(x)$ é dado por $L(C) = \sum_{x \in \mathcal{X}} p(x)l(x)$. 
\end{df}

\

Sem perdas de generalidade, podemos assumir que o alfabeto $D$-ário (o alfabeto com $D$ símbolos) é dado por $D = \lbrace 0,1,\ldots,D-1\rbrace$. Seguem alguns exemplos de códigos:

\

\begin{ex} Seja $X$ uma variável aleatória com a seguinte distribuição de probabilidade e código: 

$$\textit{ }\textbf{P}(X = 1) = \frac{1}{2} , \quad C(1) = 0 \quad$$
$$\textbf{P}(X = 2) = \frac{1}{4} , \quad C(2) = 10 $$
$$ \textit{ }\textbf{P}(X = 3) = \frac{1}{8} , \quad C(3) = 110 $$
$$ \textit{ }\textit{ }\textbf{P}(X = 4) = \frac{1}{8} , \quad C(4) = 111. $$

\

A entropia $H(X)$ é igual a 1.75 bits assim como o comprimento esperado $L(C)$. Nota-se também que qualquer sequência de bits pode ser univocamente decodificada. Por exemplo, a sequência 0110111100110 é decodificada como 134213.\textcolor{white}{a}\QEDB

\end{ex}

\

\begin{ex} Considere agora a seguinte distribuição de probabilidade e código: 

$$ \textbf{P}(X = 1) = \frac{1}{3}, \quad C(1) = 0 $$
$$ \textit{ }\textbf{P}(X = 2) = \frac{1}{3}, \quad C(2) = 10 $$
$$ \textit{ }\textit{ }\textbf{P}(X = 3) = \frac{1}{3}, \quad C(3) = 11. $$

\

De novo temos um código univocamente decodificável. A diferença é que nesse exemplo temos $H(X) = \log 3 ( \approx 1.58)$ bits e $L(X) = 1.66$ bits.\textcolor{white}{a}\QEDB

\end{ex}

\

\begin{ex}[Código Morse] O Código Morse é um código razoavelmente eficiente para o alfabeto inglês usando um alfabeto de quatro símbolos: um ponto, um traço, um espaço de letra e um espaço de palavra. Sequências pequenas representam letras frequentes (por exemplo, um ponto representa a letra e) ao passo que sequências grandes representas letras pouco frequentes (por exemplo, Q é representado por ''traço,traço,ponto,traço”). Esse não é o código ótimo do alfabeto em quatro símbolos. Muitas palavras não são usadas pois uma palavra representando uma letra não contém espaços a não ser por um espaço de letra ao final e nenhum espaço segue outro espaço. É um problema interessante calcular o número de sequências que podem ser construídas obedecendo essas restrições. Tal problema foi resolvido por Shannon em \cite{Shannon}. O problema também está relacionado com estabelecer códigos para gravação magnética onde longas sequências de 0's são proibidas como em \cite{1} e \cite{2}.\textcolor{white}{a}\QEDB

\end{ex}

\

Os códigos são classificados de acordo com propriedades que obedecem quanto a decodificação. Vamos a elas: 

\

\begin{df} Um código é não-singular se todo elemento de $\mathcal{X}$ está associado a uma sequência diferente de $D^*$. Ou seja:

$$  x \neq x' \Rightarrow C(x)\neq C(x').$$
\end{df}

\

Não-singularidade implica em não ambiguidade na tradução de uma palavra dada. Mas, em geral, uma mensagem contém uma sequência de valores de $\mathcal{X}$. Nesse caso, nós poderíamos garantir que a mensagem será decodificada univocamente acrescentando um símbolo (exclusivo para esse fim) entre duas palavras. Claramente não é uma abordagem muito eficiente. Esse problema motiva as seguintes definições:

\

\begin{df} A extensão $C^*$ de um código $C$ é a função que leva sequências finitas de $\mathcal{X}$  em sequências finitas de $D$, definida por 

$$C(x_1,x_2, \ldots,x_n) = C(x_1)C(x_2)\ldots C(x_n),$$

onde $C(x_1)C(x_2)\ldots C(x_n)$ é a concatenação das palavras correspondentes. 
\end{df}

\

\begin{df} Um código é univocamente decodificável se sua extensão é não-singular (invertível).
\end{df}

\

Em outras palavras,  cada sequência de palavras pode ser gerada por apenas uma sequência em $\mathcal{X}$. Note, porém, que é possível que se tenha que ler a sequência inteira para descobrir quem a gerou. 
\begin{df} Um código é dito ser instantâneo se nenhuma palavra é um prefixo para outra palavra. 

\end{df}

Um código instantâneo pode ser decodificado sem a leitura de palavras futuras pois o fim de cada palavra é imediatamente reconhecido. Por exemplo, a sequência binária 01011111010 produzida pelo código do exemplo C.14 é lida como 0,10,111,110,10. Abaixo, a figura 3 ilustra a relação entre os diferentes códigos aqui mencionados. A tabela C.1 mostra um exemplo de códigos distintos usando o mesmo alfabeto para o mesmo número de eventos. Para o código não-singular a palavra 010 tem três fontes possíveis: 2, 14 e 31 (logo não se trata de um código univocamente decodificável). O código univocamente decodificável não é instantâneo (pois $C(3)$ é um prefixo de $C(4)$).

\newpage

\begin{center}

\begin{tikzpicture}

\draw[black] (2.5,0) ellipse (4.5 and 5.3) (2.5,4.7) node [anchor = north] {Todos os Códigos};

\draw[black] (2.5,-0.5) ellipse (3.5 and 4.5) (2.5,3.15) node [anchor = north] {Códigos Não-Singulares};

\draw[black] (2.5,-1) ellipse (2.5 and 3.6) (2.5,1.2) node [anchor = north] {Códigos Univocamente} (2.5, 0.8) node [anchor = north] {Decodificáveis};

\draw[black] (2.5,-2.1) ellipse (2 and 2.1)  node [anchor = center] {Códigos Instantâneos};

\draw (2.5, -6) node {Figura 3: Relação de contenção entre códigos.};

\end{tikzpicture}

\end{center}

\begin{center}

\begin{table}[h]
\centering

\begin{tabular}{|c|c|c|c|c|}

\hline 
 
 & & Não-Singular  & Univoc. Decod. & \\
                              
$X$ & Singular & mas não  & mas não  & Instantâneo \\

 & & Univoc. Decod. & Instantâneo & \\

\hline
1 & 0 & 0 & 10 & 0\\
\hline
2 & 0 & 010 & 00 & 10\\
\hline
3 & 0 & 01 & 11 & 110 \\ 
\hline
4 & 0 & 10 & 110 & 111 \\ 
\hline  

\end{tabular}

\vspace{0.2cm}
\caption{Exemplos de códigos.}

\end{table}

\end{center}

 Para verificar que se trata de um código univocamente decodificável, considere uma sequência de letras: se os primeiros dois bits forem 00 ou 10 eles podem ser decodificados imediatamente. Se os primeiros bits forem 11 é preciso ler o próximo bit. Se for lido 1, lemos $C(3)$. Se a sequência de 0's imediatamente após a dupla 11 for ímpar então lemos $C(4)$. Finalmente, se a sequência de 0's for par, lemos $C(3)$. Esse argumento pode ser repetido para ler cada palavra dessa sequência. Sardinas e Patterson em \cite{3} criaram um teste finito para univoca-decodificabilidade. O último código da tabela C.1 é claramente instantâneo.

\subsection{Desigualdade de Kraft}

\

Como é de se esperar, nosso objetivo é, dado $X$, construir um código instantâneo com o menor comprimento esperado possível. Na próximo seção indicaremos como fazer isso. Antes precisamos de alguns alicerces: 

\

\begin{teo}[Desigualdade de Kraft] Sejam $X$ uma v.a. e  $C(X) $ um código instantâneo usando um alfabeto $D$-ário. Então $l_1, l_2, \ldots, l_m$, os comprimentos das (finitas) palavras de $C(X)$, devem satisfazer

$$ \sum_i D^{-l_i} \leq 1. $$

Reciprocamente, dado um conjunto de comprimentos satisfazendo a desigualdade acima então existe um código cujas palavras têm tais comprimentos.
\end{teo}

\textbf{Prova:}  Considere a árvore $D$-ária onde cada vértice tem $D$ filhos. Assim seus galhos representam as palavras de $C(X)$. Por exemplo, os $D$ filhos da raiz representam as $D$ possibilidades para a primeira letra de uma palavra e assim por diante. Logo, o conjunto dos caminhos finitos a partir da raiz é uma representação do conjunto $D^*$. Um exemplo usando a árvore binária é mostrado na figura 4 abaixo. A condição do código ser instantâneo implica que nenhuma palavra é ancestral de outra na árvore. Seja $l_{max}$ o comprimento da maior palavra de $C(X)$. Considere todos os $D^{l_{max}}$  vértices de profundidade $l_{max}$. Eles são  ou parte de uma palavra ou descendentes de uma palavra ou nenhuma das duas. Uma palavra com tamanho $l_i$ tem $D^{l_{max}-l_i}$ descendentes no nível $l_{max}$. Os conjuntos de descendentes de palavras distintas de profundidade $l_{max}$  são disjuntos. Logo, somando em todas as palavras nós temos 

$$ \sum_i D^{l_{max}-l_i} \leq D^{l_{max}}. $$ 

Ou, como desejávamos,

$$ \sum_i D^{-l_i} \leq 1. $$

\begin{center}

\begin{tikzpicture}

\draw (-1,-7) -- (1,-5.5) (1, -5.35) node [anchor = east] {0};

\draw[dashed] (1,-5.5) -- (3,-4.55);

\draw[dashed] (1,-5.5) -- (3,-6.45);

\draw[dashed] (3,-4.55) -- (5,-3.9);

\draw[dashed] (3,-4.55) -- (5,-5.2);

\draw (-1,-7) -- (1,-8.5);

\draw (1,-8.5) -- (3,-7.55)  (3, -7.4) node [anchor = east] {10};

\draw (1,-8.5) -- (3,-9.45);

\draw[dashed] (3,-7.55) -- (5,-6.9);

\draw[dashed] (3,-7.55) -- (5,-8.2);

\draw (3,-9.45) -- (5,-8.8)  (5, -8.65) node [anchor = east] {110};

\draw (3,-9.45) -- (5,-10.1)  (5, -10.25) node [anchor = east] {111};

\draw (2.5, -11.5) node {Figura 4: Código binário respeitando a desigualdade de Kraft. };

\end{tikzpicture}

\end{center}

 Reciprocamente, seja $ L =  \lbrace l_1,l_2,\ldots,l_m\rbrace$ um conjunto que satisfaz a desigualdade de Kraft. Defina o primeiro (seguindo alguma ordem arbitrária) galho de tamanho $l_1$ como a primeira palavra de $C(X)$ e remova seus descendentes da árvore. Defina o primeiro galho remanescente de tamanho $l_2$ como a segunda palavra e assim por diante. Claramente, ao final desse processo, teremos um código instantâneo cujas palavras tem comprimento $l_1, l_2, \ldots$ e $l_m$. \textcolor{white}{a}\QEDA
 
 \

O próximo teorema é uma generalização da desigualdade de Kraft.

\

\begin{teo} Para qualquer conjunto contável de palavras de um código instantâneo usando um alfabeto $D$-áreo, seus comprimentos devem satisfazer a seguinte desigualdade:

$$  \sum_i^\infty D^{-l_i} \leq 1.$$

Reciprocamente, dados $l_1,l_2,\ldots$ satisfazendo a desigualdade acima, existe um código instantâneo cujas palavras tem esses comprimentos.  
\end{teo}

\textbf{Prova:} Sejam  $y^i = y_1y_2\ldots y_{l_i}$ a $i$-ésima palavra de $C(X)$ e $\alpha(y^i) = 0.y_1y_2\ldots y_{l_i}$ o número real dado pela expansão $D$-ária 

$$ 0.y_1y_2\ldots y_{l_i} = \sum_{j = 1}^{l_i} y_j D^{-j}. $$

O conjunto $I^i$ de todas as expansões $D$-áreas de palavras que começam com $y_1y_2\ldots y_{l_i}$ é um subconjunto de $[0.y_1y_2\ldots y_{l_i}, 0.y_1y_2\ldots y_{l_i} + \frac{1}{D^{l_i}})$ (para entender melhor essa construção veja os dois exemplos subsequentes à prova). Com efeito, a palavra que começa por $y_1y_2\ldots y_{l_i}$ e tem a maior expansão $D$-ária (em módulo) é dada por $y_1y_2\ldots y_{l_i}(D - 1)(D - 1)\ldots$ . Tal palavra obedece 

$$ 0.y_1y_2\ldots y_{l_i}(D - 1) \ldots - 0.y_1y_2\ldots y_{l_i} = \sum_{j = 1}^{\infty} y_j D^{-j}  - \sum_{j = 1}^{l_i} y_j D^{-j} = $$

$$ = \sum_{j = l{i+1}}^{\infty} (D-1) D^{-j} = \frac{(D - 1) D^{- (l_i + 1)}}{1 - \dfrac{1}{D}} = \frac{1}{D^{l_i}}.$$ 

Na verdade esses conjuntos são iguais. Com efeito, o conjunto de todas as expansões $D$-árias de palavras que começam com $y_1y_2\ldots y_{l_i}$ é o conjunto de todos os número na base $D$ entre $0.y_1y_2\ldots y_{l_i}$ e $0.y_1y_2\ldots y_{l_i} + \frac{1}{D^{l_i}}$.
\

Como o código é instantâneo, para palavras distintas os intervalos aos quais sua expansões pertencem são disjuntos. Como tais intervalos estão contidos em $[0,1]$, temos que 

$$ 1 = \vert [0,1] \vert \geq \vert \bigcup_i I^i \vert = \sum_i \vert I^i \vert = \sum_i D^{-l_i}.$$

Assim como no caso finito, podemos construir um código cujas palavras tem comprimentos $l_1,l_2,\ldots$  satisfazendo a desigualdade de Kraft. Primeiro, ordene os comprimentos de forma que $l_1 \leq l_2, \leq \ldots$.  O próximo  passo é associar os intervalos às palavras em sequência do mais próximo ao zero até o mais próximo ao 1. \begin{samepage} Por exemplo, se desejamos construir um código binário com $l_1 = 1, l_2 = 2, \ldots$ nós associamos o intervalo $[0,1/2),[1/2, 	1/4),\ldots$ às palavras correspondentes $0, 10, \ldots \textit{ }$.\textcolor{white}{a}\QEDA

\

\end{samepage}

Os seguintes exemplos foram retirados de \cite{Figueiredo}.

\

\begin{ex} Considere o caso de um código 10-ário (ou decimal). Seja $y^i = 2738$. Assim, $\alpha(y^i) = 0.2738$ (na habitual escrita em base 10) e $I^i = [0.2738, 0.2739[$, o qual é o intervalo de todos os números reais cuja escrita decimal
começa por 0.2738; por exemplo 0.273845 $\in$ [0.2738, 0.2739[.\textcolor{white}{a}\QEDB

\end{ex}

\

\begin{ex} No caso de um código binário, seja $y^i = 100101$. Nesse caso, $\alpha(y^i) = 0.100101$
(em base 2, ou seja, traduzindo para base 10, $\alpha(y^i) = 1/2 + 1/16 + 1/64 = 0.5781).$
O intervalo correspondente  é $I^i = [0.100101, 0.10011[$ (pois, em base 2, 0.100101 +
0.000001 = 0.10011), o qual contém todos os números cuja escrita em base 2 começa por
0.100101; por exemplo 0.100101101 $\in$ [0.100101, 0.10011[.\textcolor{white}{a}\QEDB
\end{ex}

\

\begin{ob} Os comprimentos das palavras em um código univocamente decodificável também obedecem a desigualdade de Kraft. Para uma prova desse fato veja \cite{Thomas}. \textcolor{white}{a}\QEDB
\end{ob}

\subsection{Códigos Ótimos}

\

Como mencionado, essa seção se dedicará a, dado $X$ e $D$, achar o código instantâneo de menor comprimento esperado. Da seção anterior, sabemos que isso é equivalente a achar os comprimentos $l_1, l_2, \ldots, l_m$ que satisfazem a desigualdade de Kraft e que minimizam $L = \sum p_i l_i$.  Esse é um típico problema de otimização: minimize L = $\sum_i p_i l_i$ sujeito a $l_1,l_2,\ldots,l_m \in \mathbb{N^*}$ e $\sum_i D^{-l_i} \leq 1$. Nós descartaremos a primeira restrição e assumiremos igualdade na segunda. Assim, podemos escrever o problema usando multiplicadores de Lagrange como a minimização de

$$ J = \sum_i p_i l_i + \lambda (\sum_i D^{-l_i}). $$ 

Diferenciando com respeito a $l_i$ obtemos 

$$\frac{\partial J}{\partial l_i} = p_i - \lambda D^{-l_i} \log_e D.$$ 

Calculando as derivadas em 0, obtemos 

$$D^{-l_i} = p_i \lambda \log_e D.$$

Substituindo esse resultado na restrição descobrimos que $\lambda = 1/\log_e D$, logo  $p_i = D^{-l_i}$ resultando nos comprimentos ótimos 

$$ l^*_i = -\log_D p_i.$$  

Essa escolha não inteira para os comprimentos resulta no comprimento esperado 

$$ L^* = \sum_i p_i l^*_i = -\sum_i p_i \log_D p_i = H_D(X).$$

 Como os comprimentos devem ser inteiros, nem sempre seremos capazes de obedecer as igualdades acima. A ideia é procurar por valores inteiros próximos aos otimais. Em vez de demonstrar por cálculo que $l^*_i  =-\log_D p_i$ é um mínimo global, verificaremos a otimalidade na prova do seguinte teorema:
 
\

\begin{teo} O comprimento esperado $L$ de um código instantâneo $D$-ário para uma v.a. $X$ é maior ou igual à entropia $H_D(X)$. Isto é

$$ L \geq H_D(X),$$

com igualdade se, e somente se, $D^{-l_i} = p_i$. 
\end{teo}

\textbf{Prova:} Podemos escrever a diferença entre o comprimento esperado e a entropia como

$$  L - H_D(X) = \sum_i p_i l_i - \sum_i p_i \log_D \frac{1}{p_i} = - \sum_i p_i \log_D D^{-l_i} + \sum_i p_i \log_D p_i.$$

Sendo  $r_i = D^{-l_i}/\sum_j D^{-l_j}$ e  $c = \sum_i D^{-l_i}$, obtemos

$$ L - H = \sum_i p_i \log_D \frac{p_i}{r_i} - \log_D c  = D(p\vert \vert r) + \log_D \frac{1}{c} \geq 0,$$  

pela não negatividade da entropia relativa (teorema C.10) e pela desigualdade de Kraft. Logo  $L \geq H$ com igualdade se, e somente se $p_i = D^{-l_i}$ (i.e., se e somente se $-\log_D p_i$ é inteiro para todo $i$). \textcolor{white}{a}\QEDA

\

\

\begin{df} Uma distribuição de probabilidade é dita ser $D$-ádica se cada probabilidade é igual a $D^{-n}$ para algum $n$. 
\end{df}

\

Logo, temos igualdade no teorema anterior se, e somente se, a distribuição de $X$ é $D$-ádica. A prova do teorema C.26 também indica um procedimento para achar um código ótimo: achar a distribuição $D$-ádica que mais se aproxima (na métrica da entropia relativa) da distribuição de $X$. Essa distribuição nos fornecerá os comprimentos das palavras do código. Já provamos, de forma construtiva, que, de posse de comprimentos satisfazendo a desigualdade de Kraft, podemos definir um código instantâneo cujas palavras têm os comprimentos dados. Assim temos um código ótimo para $X$. Ainda há um problema: a tarefa de achar a distribuição $D$-ádica mais próxima da distribuição de $X$ não é trivial. Há rotinas bem sucedidas para esse fim. O Código de Shannon-Fano mostra, além disso, que pode-se construir um código com comprimento médio tão próximo (pela direita) quanto se queira da entropia de $X$. O código de Huffman prova-se otimal, i.é., nenhum outro código instantâneo tem comprimento esperado menor que o obtido pelo código de Huffman. Para aprofundar no que acabamos de expor, o leitor pode consultar \cite{Thomas}.

\end{document}